\documentclass[11pt]{article}
\usepackage{amsmath,amssymb,amsthm,float}
\usepackage{fullpage}
\usepackage{longtable}
\usepackage{enumerate}
\usepackage{enumitem}

\usepackage{array}
\usepackage{calc}

\newtheorem{maintheorem}{Theorem}
\newtheorem{maincorollary}[maintheorem]{Corollary}
\newtheorem{theorem}{Theorem}[section]
\newtheorem{lemma}[theorem]{Lemma}
\newtheorem{proposition}[theorem]{Proposition}
\newtheorem{corollary}[theorem]{Corollary}

\theoremstyle{definition}
\newtheorem{remark}[theorem]{Remark}
\newtheorem{definition}[theorem]{Definition}

\setlist[enumerate]{topsep=0pt}

\title{Complete Reducibility in Good Characteristic}
\author{Alastair J. Litterick and Adam R. Thomas
}
\renewcommand\footnotemark{}

\begin{document}

\thispagestyle{empty}

\begin{center}
\LARGE Corrigendum: Complete Reducibility in Good Characteristic
\end{center}

Subsequent to the publication of this paper, an omission was noted by the authors. Specifically, two classes of non-$G$-cr subgroups of type $A_1$ were missed when $G = E_6$, $p = 5$ which are $\operatorname{Aut}(G)$-conjugate, leading also to a missing class of subgroups when $G = E_7$, $p = 5$.

The essential mistake is in Proposition~\ref{prop:irredA1}, with an omitted subgroup of type $A_{1}$ embedded into a Levi subgroup of type $D_5$ via the $10$-dimensional module $(2 \otimes 2^{[1]}) + 0$. This arXiv version now includes this subgroup and the non-$G$-cr subgroups arising from this. Fortunately, the qualitative results of the paper remain unchanged.

The specific additions are an additional final line in Tables \ref{E6p5tab} and \ref{E7p5tab}, and corresponding minor updates to the statements of Corollaries \ref{cor:overgroups}, \ref{cor:max} and \ref{cor:trivcentraliser}. In each case, the new subgroup classes are placed at the end of the relevant table or corollary statement.

As of October 2023, we believe this corrigendum does not affect the results of any papers citing this one. The nontrivial unipotent elements in these new non-$G$-cr subgroups $A_1$ lie in the $E_6$-class (and thus also the $E_7$-class) labelled $A_4 + A_1$ in \cite{MR1351124}. In particular, the omitted subgroups do not affect the specific reference to our results in Lemma~6.4 of [Korhonen, Mikko, \emph{Unipotent elements forcing irreducibility in linear algebraic groups}. J.\ Group Theory {\bf 21} (2018), no.~3, 365–396.], or the references in [Malle, Gunter; Testerman, Donna M. \emph{Overgroups of regular unipotent elements in simple algebraic groups.} Trans.\ Amer.\ Math.\ Soc.\ Ser.\ B {\bf 8} (2021), 788–822.].

\begin{flushright}
Alastair Litterick\\
Adam Thomas\\
October 2023.
\end{flushright}

\addtocounter{page}{-1}

\clearpage

\maketitle

\begin{abstract}
Let $G$ be a simple algebraic group of exceptional type, over an algebraically closed field of characteristic $p \ge 0$. A closed subgroup $H$ of $G$ is called $G$-completely reducible ($G$-cr) if whenever $H$ is contained in a parabolic subgroup $P$ of $G$, it is contained in a Levi subgroup of $P$. In this paper we determine the $G$-conjugacy classes of non-$G$-cr simple connected subgroups of $G$ when $p$ is good for $G$. For each such subgroup $X$, we determine the action of $X$ on the adjoint module $L(G)$ and the connected centraliser of $X$ in $G$. As a consequence we classify all non-$G$-cr connected reductive subgroups of $G$, and determine their connected centralisers. We also classify the subgroups of $G$ which are maximal among connected reductive subgroups, but not maximal among all connected subgroups.
\end{abstract}

\tableofcontents

\section{Introduction}

\vspace{-3pt}Let $G$ be a connected reductive algebraic group defined over an algebraically closed field $K$. Following Serre \cite{Ser3}, a closed subgroup $X$ of $G$ is said to be \emph{$G$-completely reducible} ($G$-cr) if whenever $X$ is contained in a parabolic subgroup $P$ of $G$, it is contained in a Levi subgroup of $P$. If $G = GL(V)$, a subgroup $X$ is $G$-cr if and only if $V$ is a completely reducible $KX$-module, and thus $G$-complete reducibility is a generalisation of the standard notion in representation theory. Similarly, $X$ is called \emph{$G$-irreducible} if $X$ is not contained in any proper parabolic subgroup of $G$, and \emph{$G$-indecomposable} if $X$ lies in no proper Levi subgroup of $G$. The $G$-cr subgroups of $G$ are precisely the $L$-irreducible subgroups of $L$, where $L$ ranges over the Levi subgroups of $G$ \cite[Corollary 3.5]{MR2178661}; the case $L = G$ yields the $G$-irreducible subgroups.

By a theorem of Borel and Tits \cite[Th\'{e}or\`{e}me 2.5]{MR0294349}, a $G$-cr subgroup is necessarily reductive. The converse is true provided that the characteristic of $K$ is zero or large relative to the root system of $G$ (cf.\ \cite{MR1635685,MR1476899,MR1329942} and \cite[Theorem 4.4]{MR2167207}); thus non-$G$-cr reductive subgroups are inherently a low positive characteristic phenomenon. In this paper we consider a weak restriction on the characteristic. Recall that $\textup{Char}(K)$ is called \emph{bad for $G$} if it is prime and divides some coefficient when the highest root in the root system of $G$ is expressed as a sum of simple roots, and \emph{good for $G$} otherwise. A number of useful subgroup structure results hold precisely when $p = \textup{Char}(K)$ is good for $G$. For instance, a result of Bate, Martin and R\"{o}hrle (Lemma \ref{lem:BMR}) states that if $p$ is good for $G$ then a closed subgroup of a subsystem subgroup $H$ is $G$-cr if and only if it is $H$-cr. Good characteristic is therefore a natural first scenario to consider when studying non-$G$-cr subgroups.

Now let $G$ be simple. When $G$ has classical type, understanding reductive subgroups of $G$ amounts to understanding the finite-dimensional representation theory of reductive groups. In this paper we consider $G$ of exceptional type, where much more explicit results can be expected. The bad primes here are $2$ and $3$ for all exceptional types, as well as $5$ for $E_{8}$. A result of Liebeck and Seitz \cite[Theorem 1]{MR1329942} states that if $p > 7$, then every closed connected reductive subgroup of $G$ is $G$-cr. In fact, if $G$ is of type $G_{2}$ or $F_{4}$, then $p > 3$ suffices, and hence for these types all connected reductive subgroups are $G$-completely reducible in good characteristic. This fails for $G$ of type $E_6$, $E_7$ or $E_8$. Indeed, it is shown in \cite[Corollary 2]{Stewart21072013} that if $p$ is good for $G$ then there exists a non-$G$-cr simple subgroup $X$ of $G$ if and only if either $X$ is of type $A_1$ with $p = 5$ or $7$, or $X$ is of type $G_2$ with $p = 7$ and $G = E_7$ or $E_8$. In this paper we classify all non-$G$-cr simple subgroups in these cases. The results are the following. The tables referenced in the statements can be found in Section \ref{sec:tables} on page \pageref{sec:tables}.

\newpage

\begin{maintheorem} \label{thm:e6}
Let $G$ be a simple algebraic group of type $E_6$ in good characteristic $p$, and let $X$ be a non-$G$-cr simple subgroup of $G$. Then $p = 5$, $X$ is of type $A_1$, and $X$ is $\textup{Aut}(G)$-conjugate to exactly one subgroup listed in Table \ref{E6p5tab}, all of which are non-$G$-cr.
\end{maintheorem}

\begin{maintheorem} \label{thm:e7}
Let $G$ be a simple algebraic group of type $E_7$ in good characteristic $p$, and let $X$ be a non-$G$-cr simple subgroup of $G$. Then either $p = 5$, $X$ is of type $A_1$, and $X$ is conjugate to exactly one subgroup listed in Table \ref{E7p5tab}; or $p = 7$, $X$ is of type $A_1$ or $G_2$, and $X$ is conjugate to exactly one subgroup listed in Table \ref{E7p7a1tab} or \ref{E7p7g2tab}. The subgroups in these tables are all non-$G$-cr.
\end{maintheorem}

\begin{maintheorem} \label{thm:e8}
Let $G$ be a simple algebraic group of type $E_8$ in good characteristic $p$, and let $X$ be a non-$G$-cr simple subgroup of $G$. Then $p = 7$, $X$ is of type $A_1$ or $G_2$, and $X$ is conjugate to exactly one subgroup listed in Table \ref{E8p7a1tab} or \ref{E8p7g2tab}, all of which are non-$G$-cr.
\end{maintheorem}

The tables in Section \ref{sec:tables} contain additional information on the non-$G$-cr simple subgroups in question. For each such subgroup $X$, we determine the $X$-module structure of the adjoint module $L(G)$ and of a non-trivial module of least dimension when $G = E_{6}$ or $E_{7}$. We also determine the connected centraliser of each non-$G$-cr simple subgroup, which allows us to extend our results to a classification of all non-$G$-cr connected reductive subgroups of $G$. In the subsequent corollaries, we are therefore able to make a number of observations concerning general reductive subgroups of $G$.

\begin{maintheorem} \label{thm:reductive}
Let $G$ be an exceptional simple algebraic group in good characteristic, and let $X$ be a non-$G$-cr connected reductive subgroup of $G$. Then one of the following holds:
\begin{enumerate}[label=\normalfont(\roman*)]
\item $X$ is simple,
\item $X$ is $\textup{Aut}(G)$-conjugate to a semisimple subgroup in Table \ref{tab:semisimple},
\item $Z(X)^{\circ} \neq 1$ and $X$ is $\textup{Aut}(G)$-conjugate to a subgroup in Table \ref{tab:reductive}.
\end{enumerate}

Each subgroup $X$ in Tables \ref{tab:semisimple} and \ref{tab:reductive} denotes a unique $\textup{Aut}(G)$-conjugacy class unless stated otherwise, and all such subgroups are non-$G$-cr.
\end{maintheorem}

The proof of Theorems \ref{thm:e6}--\ref{thm:e8} generalises the strategy developed in \cite{MR3075783}, which classifies non-$G$-cr connected reductive subgroups when $G$ has type $F_4$ with $p = 2$ or $3$. An outline is as follows. For each proper parabolic subgroup $P$ of $G$, with unipotent radical $Q$ and Levi factor $L$, we find all $L$-irreducible simple subgroups of type $A_1$ or $G_2$. This uses standard representation theory when $L$ has classical simple factors, and results of the second author \cite{Tho1,Tho2} when $L$ has exceptional simple factors. For each such $L$-irreducible subgroup $X$, we then study the cohomology set $H^1(X,Q)$, which parametrises the $Q$-conjugacy classes of complements to $Q$ in $QX$. This is the most technical part, and is discussed in detail in Section \ref{sec:preliminaries}. In each case, we either describe $H^{1}(X,Q)$ explicitly, or we determine sufficient information to limit the number of conjugacy classes of subgroups of $P$ which are isomorphic to $X$ and not conjugate to a subgroup of $L$. Finally we construct explicit representatives of each class of subgroups, for instance through known embeddings into a proper reductive subgroup as in Corollary \ref{cor:overgroups} below.

We note that the techniques described in Section \ref{sec:preliminaries} relating to $H^{1}(X,Q)$ are equally valid in bad characteristic. In addition, we still have explicit descriptions of the simple $L$-irreducible subgroups of each Levi subgroup $L$ in this case. However, results along the lines of Lemma \ref{lem:BMR}, which are intrinsic in our construction of non-$G$-cr subgroups, can fail in bad characteristic. Thus extending Theorems \ref{thm:e6}--\ref{thm:reductive} to bad characteristic presents considerable technical difficulties which do not arise here; we intend to explore this in future work.

We now present a series of consequences of Theorems \ref{thm:e6}--\ref{thm:reductive}. We begin with an observation on how simple subgroups of $G$ act on certain $G$-modules of small dimension. Recall that a module for a reductive algebraic group is called \emph{tilting} if it has a filtration by Weyl modules and a filtration by dual Weyl modules. As defined in Section \ref{sec:notation}, the notation $V^{[e]}$ denotes a twist of the module $V$ by a $p^{e}$-power Frobenius morphism. Our first corollary extends Theorems 3 and 4 of \cite{MR1973585}, which treat the case where $X$ is $G$-cr:

\begin{maincorollary} \label{cor:variations}
Let $G$ be an exceptional simple algebraic group in good characteristic $p$, and let $X$ be a simple subgroup of $G$. Then $L(G) \downarrow X$ is a direct sum of modules of the form $V_{1}^{[r_1]} \otimes \cdots \otimes V_{k}^{[r_k]}$, and exactly one of the following holds:
\begin{enumerate}[label=\normalfont(\roman*)]
\item Each $V_{i}$ is a Weyl module, a dual Weyl module or a tilting module,
\item $p = 7$ and $X = G_{2}$ is a maximal subgroup of $G = F_{4}$,
\item $p = 7$ and $X = G_{2}$ is a non-$G$-cr subgroup of $G = E_{8}$ contained in a maximal subgroup $G_{2}F_{4}$ (cf.~Theorem \ref{thm:e8}).
\end{enumerate}
\end{maincorollary}
The subgroups in (ii) and (iii) are genuine exceptions, unique up to conjugacy in $G$, and both are listed in \cite[Theorem 4]{MR1973585}, although the subgroup in (iii) was not shown to be non-$G$-cr at that time.

\begin{remark} In \cite[Theorem 4]{MR1973585}, the tensor factors in the direct summands of $L(G) \downarrow X$ arise from a factorisation $X \to X \times X \times \ldots \times X \to G$ of the inclusion of $X$ into $G$. Inspecting the first and second columns of Tables \ref{E6p5tab}--\ref{E8p7g2tab} shows that in good characteristic, the same holds for non-$G$-cr subgroups $X$; if $L(G) \downarrow X$ has a summand which is a tensor product $V_{1}^{[r]} \otimes V_{2}^{[s]}$, then $X = A_{1}$ lies in a subgroup $A_{1}A_{1}$ of $G$, such that $V_{i}$ is a Weyl module, dual Weyl module or tilting module for the $i$-th factor of this subgroup.
\end{remark}

\begin{remark} Inspection of Tables \ref{E6p5tab}--\ref{E7p7g2tab} also shows that if $G = E_{6}$ or $E_{7}$ and $X$ is non-$G$-cr then $V_{G}(\lambda_1) \downarrow X$ (resp.\ $V_{G}(\lambda_7) \downarrow X$) also satisfies part (i) of the above corollary.
\end{remark}

Next, we consider overgroups of the non-$G$-cr subgroups arising. In the following, a \emph{subsystem subgroup} of $G$ is a semisimple subgroup which is normalised by a maximal torus of $G$.

\begin{maincorollary} \label{cor:overgroups}
Let $G$ be an exceptional simple algebraic group in good characteristic $p$, and let $X$ be a non-$G$-cr connected reductive subgroup of $G$. Then exactly one of the following holds:
\begin{enumerate}[label=\normalfont(\roman*)]
\item The semisimple subgroup $X'$ is properly contained in a proper subsystem subgroup of $G$,
\item $p = 7$ and $X = A_{1}$ lies in a maximal subgroup $A_{1}G_{2}$ of $G = E_{7}$,
\item $p = 7$ and $X = G_{2}$ is maximal among proper connected reductive subgroups of $G = E_{7}$,
\item $p = 7$ and $X = G_{2}$ lies in a maximal subgroup $G_{2}F_{4}$ of $G = E_{8}$.
\item $p = 5$ and $X = A_{1}$ is maximal among proper connected reductive subgroups of $G = E_{6}$.
\end{enumerate}
Therefore a connected subgroup of $G$ which is maximal among connected reductive subgroups, is either a maximal connected subgroup of $G$, a Levi factor of a maximal parabolic subgroup of $G$, or is conjugate to a group $X$ in \textup{(iii)} or \textup{(v)}.

The subgroups $X$ in \textup{(ii)}, \textup{(iii)} and \textup{(iv)} are each unique up to conjugacy in $G$, and \textup{(v)} gives two $G$-classes of subgroups which are fused in $\operatorname{Aut}(G)$. Their embeddings are described in Table \ref{E7p7a1tab}, \ref{E7p7g2tab}, \ref{E8p7g2tab} and \ref{E6p5tab}, respectively.
\end{maincorollary}

Our next result describes certain chains $X \le M < G$ where $X$ and $M$ are reductive and $X$ is $G$-indecomposable and $M$-irreducible but not $G$-irreducible. This contrasts with Lemma \ref{lem:BMR}, which tells us that $M$ cannot be a subsystem subgroup of $G$. The notation for embeddings in Table \ref{tab:max} is defined in Section \ref{sec:notation}.

\begin{maincorollary} \label{cor:max}
Let $G$ be an exceptional simple algebraic group in good characteristic $p$. Let $X \le M$ be connected subgroups of $G$ such that $M$ is maximal among proper connected reductive subgroups of $G$, and $X$ is $M$-irr and non-$G$-cr. Then $G$, $p$, $X$ and $M$ appear in Table \ref{tab:max}, and all such chains $X \le M < G$ satisfy the given hypotheses.
\end{maincorollary}

\begin{longtable}{>{\raggedright\arraybackslash}p{0.045\textwidth - 2\tabcolsep}>{\raggedright\arraybackslash}p{0.035\textwidth - 2\tabcolsep}>{\raggedright\arraybackslash}p{0.20\textwidth - 2\tabcolsep}>{\raggedright\arraybackslash}p{0.475\textwidth-2\tabcolsep}>{\raggedright\arraybackslash}p{0.245\textwidth-\tabcolsep}@{}}
\caption{Chains $X \le M < G$ where $X$ is non-$G$-cr} \label{tab:max} \\
\hline
$G$ & $p$ & $M$ & $X$ & Embedding of $X$ in $M$ \\
\hline
$E_6$ & 5 & $A_{2}$ (2 classes)& $A_1 \hookrightarrow \bar{A}_1 A_5$ via $(1,W(5))$ (resp.\ $(1,W(5)^*)$) & via $2$ \\
$E_7$ & 5 & $A_{1}A_{1}$ & $A_1 \hookrightarrow A_2 A_5$ via $(2,W(5))$ & via $(1,1)$ \\
& 5 & $A_{2}$ & $A_1 < A_7$ via $W(7)$ & via $2$ \\
$E_{7}$ & 7 & non-$G$-cr $G_{2}$ & $X = M = G_{2}$ & -- \\
& 7 & $A_{2}$; non-$G$-cr $G_2$ & $A_1 < A_7$ via $W(7)$ & via $2$; via $6$ \\
& 7 & $A_{1}G_{2}$; $G_{2}C_{3}$ & $A_1 \hookrightarrow A_1 G_2$ via $(1,6)$ & via $(1,6)$; via $(6,5)$ \\
$E_8$ & 7 & $B_{2}$ & $A_1 < A_8$ via $W(8)$ & via $4$ \\
& 7 & $G_{2}F_{4}$ & $G_{2} \hookrightarrow G_{2}G_{2} \textup{ max } G_{2}F_{4}$ via $(10,10)$ & -- \\
$E_6$ & 5 & non-$G$-cr $A_1$  & $X = M = A_1$ & -- \\
& & (2 classes) & \\
 \hline
\end{longtable}

Recall from \cite{MR1973585} that a simple subgroup $X$ of an exceptional simple algebraic group $G$ is called \emph{restricted} if either $X = A_{1}$ and each high weight of a composition factor of $L(G) \downarrow X$ is at most $2p - 2$, or $X \neq A_{1}$ and whenever a high weight of a composition factor of $L(G) \downarrow X$ is expressed as a sum of fundamental dominant weights, each coefficient is at most $p - 1$. A semisimple subgroup is then called \emph{restricted} if each of its simple factors is restricted. Our next corollary complements Theorem 1.1 of \cite{MR1779618}, which tells us that a restricted subgroup of type $A_{1}$ (called a ``good $A_{1}$'' there) is $G$-completely reducible.

\begin{maincorollary} \label{cor:restricted}
Let $G$ be an exceptional simple algebraic group in good characteristic $p$, and let $X$ be a restricted semisimple subgroup of $G$. Then either $X$ is $G$-cr, or $p = 7$ and $X$ has type $G_{2}$ or $A_{1}G_{2}$. In the latter cases, every such non-$G$-cr subgroup $X$ is restricted.
\end{maincorollary}

The connected centraliser of a $G$-cr subgroup is reductive \cite[Lemma 3.12]{MR2178661}, and it follows that a $G$-cr subgroup has trivial connected centraliser if and only if it is $G$-irreducible. Theorems \ref{thm:e6}--\ref{thm:reductive} allow us to classify those connected reductive subgroups of $G$ that have trivial connected centraliser but are not $G$-irreducible.   

\begin{maincorollary} \label{cor:trivcentraliser}
Let $G$ be an exceptional simple algebraic group in good characteristic $p$, and let $X$ be a semisimple subgroup of $G$. Then $C_{G}(X)^\circ = 1$ if and only if $X$ is either $G$-irreducible or conjugate to one of the non-$G$-cr subgroups in Table \ref{tab:centraliser}.
\end{maincorollary}

\begin{longtable}{>{\raggedright\arraybackslash}p{0.06\textwidth - 2\tabcolsep}>{\raggedright\arraybackslash}p{0.06\textwidth - 2\tabcolsep}>{\raggedright\arraybackslash}p{0.65\textwidth-\tabcolsep}@{}}
\caption{Non-$G$-cr semisimple $X$ with $C_{G}(X)^{\circ} = 1$.} \label{tab:centraliser} \\
\hline
$G$ & $p$ & Non-$G$-cr subgroup \\
\hline
$E_6$ & 5 & $A_1 \hookrightarrow \bar{A}_1 A_5$ via $(1^{[r]},W(5)^{[s]})$ or $(1^{[r]},(W(5)^*)^{[s]})$ $(rs=0)$ \\
& 5 & $\bar{A}_1 A_1 < \bar{A}_1 A_5$ where $A_1 < A_5$ via $W(5)$ or $W(5)^*$ \\
$E_7$ & 5 & $A_1 \hookrightarrow A_2 A_5$ via $(2^{[r]},W(5)^{[s]})$ $(rs=0)$ \\
& 5 & $A_1 A_1 < A_2 A_5$ where $A_1 < A_2$ via 2 and $A_1 < A_5$ via $W(5)$ \\
& 5 & $A_2 A_1 < A_2 A_5$ where $A_1 < A_5$ via $W(5)$ \\
& 7 & $A_1 \hookrightarrow A_1 G_2$ via $(1,6)$ \\
$E_8$ & 7 & $A_1 < A_8$ via $W(8)$ \\
& 7 & $A_1 \hookrightarrow \bar{A}_1 A_1 G_2 < \bar{A}_1 E_7$ via $(1^{[r]},1^{[s]},6^{[s]})$ $(r \neq s$; $rs=0)$ \\
& 7 & $\bar{A}_1 G_2 < \bar{A}_1 E_7$ where $G_2 < E_7$ is non-$E_7$-cr (cf.~\ref{thm:e7})\\
& 7 & $G_2 \hookrightarrow G_2 G_2 < G_2 F_4$ via $(10,10)$ \\
$E_6$ & 5 & $A_1$ in no proper reductive overgroup (2 classes; cf.~Theorem \ref{thm:e6}) \\
\hline
\end{longtable}

\begin{remark}
Theorems \ref{thm:e6}--\ref{thm:reductive} in fact classify all non-$G$-cr connected subgroups in good characteristic with reductive centraliser. In particular, in good characteristic every connected reductive subgroup of $G = E_{6}$ has a reductive centraliser, but this is not true for $G = E_{7}$ and $E_{8}$.
\end{remark}

\begin{remark}
For a reductive subgroup $X$ of $G$, let $S$ be a maximal torus of $C_G(X)$ and let $H = C_G(S)$. Then $H$ is minimal among Levi subgroups of $G$ containing $X$ (possibly $H = G$). We observe that for each non-$G$-cr subgroup $X$ appearing in Theorems \ref{thm:e6}--\ref{thm:reductive}, the reductive subgroup $C_G(H')$ is a complement to the unipotent radical of $C_G(X)^{\circ}$. We are not aware of a general reason for this phenomenon. Note that $S$ is also a maximal torus of $C_G(H')$, hence $C_G(X)^{\circ}/R_{u}(C_G(X)^{\circ})$ and $C_G(H')$ necessarily have the same rank.
\end{remark}

\begin{remark}
In the corresponding scenario with $G$ of classical type, we do not expect to be able to classify non-$G$-cr subgroups $X$ with $C_{G}(X)^{\circ} = 1$. For instance, if $V$ is a non-trivial irreducible $X$-module with $H^{1}(X,V) \neq 0$, then $X$ occurs as a non-$G$-cr subgroup in a maximal parabolic subgroup of $G = SL(V \oplus K)$. The corresponding unipotent radical is isomorphic to $V$ as an algebraic $X$-group, hence $C_G(X)^{\circ}$ contains no unipotent elements, and is therefore a torus. Since $X$ does not lie in a proper Levi subgroup of $G$, it cannot centralise a non-trivial torus, and so $C_{G}(X)^{\circ} = 1$.
\end{remark}

Lastly, we make an observation on the number of $G$-conjugacy classes of semisimple subgroups. In principle, if $P$ is a parabolic subgroup of $G$ with Levi decomposition $P = QL$, and $X \le L$ is an $L$-irreducible reductive subgroup, then the conjugacy classes of complements to $Q$ in $QX$ can depend on the cardinality of the underlying algebraically closed field $K$. For instance if $\textup{dim}\, H^{1}(X,V) = n \ge 2$ for some irreducible $X$-module $V$, then $X$ occurs as a subgroup of $SL(V \oplus K)$ with image in a maximal parabolic subgroup $P = QL$, where $L \cong GL(V)$ and $Q \cong V$ as algebraic $X$-groups. The quotient of $H^{1}(X,Q)$ by the action of $Z(L)$ is a projective variety over $K$ of dimension $n - 1$, parametrising conjugacy classes of complements to $Q$ in $QX$. We observe that this phenomenon does not occur for exceptional groups in good characteristic:

\begin{maincorollary} \label{cor:fincomps}
Let $G$ be an exceptional simple algebraic group in good characteristic $p$, let $P$ be a parabolic subgroup of $G$ with Levi decomposition $P = QL$, and let $X \le L$ be an $L$-irreducible simple subgroup. Then complements to $Q$ in $QX$ fall into finitely many $G$-conjugacy classes of subgroups.
\end{maincorollary}

As a consequence of \cite[Theorem 3]{MR2043006}, every Levi subgroup $L$ of $G$ has either finitely many $L$-irreducible simple subgroups of a given isomorphism type, up to $L$-conjugacy, or a countably infinite number, depending on a choice of a finite number of field twists. The following is thus immediate from Corollary \ref{cor:fincomps}.

\begin{maincorollary} \label{cor:countable}
Let $G$ be an exceptional simple algebraic group, over an algebraically closed field $K$ of characteristic $p$ which is good for $G$. Then $G$ has countably many conjugacy classes of connected reductive subgroups.

In particular, the classes of non-$G$-cr simple subgroups occurring are as follows: Countably infinitely many classes of subgroups $A_{1}$ when $(G,p) = (E_6,5)$, $(E_7,5)$ or $(E_8,7)$; two classes of subgroups $A_{1}$ and one class of subgroups $G_{2}$ when $(G,p) = (E_7,7)$; and two classes of subgroups $G_{2}$ when $(G,p) = (E_8,7)$.
\end{maincorollary}

Note that the corresponding result does not hold for exceptional groups in bad characteristic, for instance if $K$ has characteristic $2$ or $3$ then by \cite[Theorem 1(B)]{MR3075783} the group $F_{4}(K)$ has a series of pairwise non-conjugate non-$G$-cr subgroups of type $A_{1}$, parametrised by closed points of a positive-dimensional variety over $K$.


\section{Notation} \label{sec:notation}

Throughout, $G$ denotes a simple algebraic group over an algebraically closed field $K$ of characteristic $p > 0$, where $p$ is good for $G$. Subgroups of $G$ are taken to be closed, and homomorphisms are taken to be morphisms of varieties. For us, ``simple'' and ``semisimple'' subgroups of $G$ will always refer to connected subgroups. In addition, conjugation will always be a left action.

Let $\Phi$ be set of roots of $G$, with respect to a fixed maximal torus $T \le G$. Let $\Pi = \{ \alpha_1, \ldots, \alpha_l \}$ be a base of simple roots corresponding to a choice of Borel subgroup containing $T$, and let $W(G) = N_G(T) / T$ be the Weyl group of $G$. Let $\{ \lambda_1, \ldots, \lambda_{l}\}$ be the set of fundamental dominant weight of $G$. We use the Bourbaki ordering on nodes of the Dynkin diagram, cf.\ \cite[Ch.\ VI, Planches I-IX]{MR0240238}. We sometimes use $a_1 a_2 \ldots a_l$ to denote either a root $a_1 \alpha_1 + a_2 \alpha_2 + \ldots + a_l \alpha_l$ or a dominant weight $a_1 \lambda_1 + a_2 \lambda_2 + \ldots + a_l \lambda_l$; context will prevent ambiguity.

For a root $\alpha$, we make the following definitions. The Weyl group $W(G)$ acts on the left on $\mathbb{Z}\Phi \otimes \mathbb{R}$, and we let $s_{\alpha}$ denote the reflection in the hyperplane perpendicular to $\alpha$. We let $U_{\alpha} = \{ x_{\alpha}(t) \, : \, t \in K \}$ denote the $T$-root subgroup of $G$ corresponding to $\alpha$, and for $t \in K^{\ast}$ we let $n_{\alpha}(t) = x_{\alpha}(t)x_{-\alpha}(-t^{-1})x_{\alpha}(t) \in N_{G}(T)$, so that $n_{\alpha}(t)$ maps to $s_{\alpha} \in W(G)$ under $N_{G}(T) \to W(G)$ (cf.\ \cite[\S 6.4]{MR0407163}). We set $n_{\alpha} = n_{\alpha}(1)$. Furthermore we let $h_{\alpha}(t) = n_{\alpha}(t)n_{\alpha}(-1) \in T$. If $\alpha = \alpha_{i}$ is a simple root, we set $s_{i} = s_{\alpha_{i}}$, and similarly $n_{i} = n_{\alpha_i}$ and $h_{i}(t) = h_{\alpha_i}(t)$.

For a dominant weight $\lambda$, we denote by $V_G(\lambda)$ (or just $\lambda$) the irreducible $G$-module of highest weight $\lambda$. Similarly, the Weyl module of highest weight $\lambda$ is denoted $W(\lambda) = W_G(\lambda)$, and the tilting module of highest weight $\lambda$ is denoted by $T(\lambda)$. The dual of a $G$-module $V$ is $V^{*} = \textup{Hom}_{K}(V,K)$. If $Y = Y_1 Y_2 \ldots Y_k$ is a commuting product of simple algebraic groups, then $(V_1, \ldots, V_k)$ denotes the $Y$-module $V_1 \otimes \cdots \otimes V_k$  where $V_i$ is a $Y_i$-module for each $i$.  Let $F\, :\, G \rightarrow G$ be the Frobenius endomorphism of $G$ which acts by sending the root element $x_{\alpha}(c)$ to $x_{\alpha}(c^{p})$, and let $V$ be a $G$-module afforded by a representation $\rho\,:\, G \rightarrow GL(V)$. Then $V^{[r]}$ denotes the module afforded by the representation $\rho^{[r]} \stackrel{\textup{def}}{=} \rho \circ F^r$. Let $M_1, \ldots, M_k$ be $X$-modules and $m_1, \ldots, m_k$ be positive integers. Then $M_1^{m_1} / \ldots / M_k^{m_k}$ denotes an $X$-module having the same composition factors as $M_1^{m_1} \oplus \ldots \oplus M_k^{m_k}$. Furthermore, $V = M_1 | \ldots | M_k$ denotes an $X$-module with a series $V = V_{1} > V_{2} > \ldots > V_{k+1} = 0$ of submodules such that $\textup{soc}(V/V_{i+1}) = V_{i}/V_{i+1} \cong M_{i}$ for $1 \le i \le k$. Finally, we let $L(G)$ denote the adjoint module for $G$, which is always irreducible for $G$ of exceptional type when $p$ is good for $G$. For $G$ of type $E_{6}$ we set $V_{27} = V_{G}(\lambda_1)$, and for $G$ of type $E_{7}$ we set $V_{56} = V_{G}(\lambda_7)$.

The notation $\bar{A}_{1}$ denotes a subgroup $A_{1}$ of $G$ which is generated by long root subgroups; we use this to help keep track of simple factors in semisimple groups. If the root system of $G$ contains short roots, then $\tilde{A}_{1}$ and $\tilde{A}_{2}$ are used to denote subgroups generated by short root subgroups.

Let $J = \{\alpha_{j_1}, \alpha_{j_2}, \ldots, \alpha_{j_r}\} \subseteq \Pi$ and define $\Phi_J = \Phi \cap \mathbb{Z}J$. Then the standard parabolic subgroup corresponding to $J$ is $P = P_{j_1 j_2 \ldots j_r} = \left< T, U_\alpha \, : \, \alpha \in \Phi_J \cup \Phi^+ \right>$. The Levi decomposition of $P$ is $P = Q L$ where $Q = Q_{j_1 j_2 \ldots j_r} = R_{u}(P) = \left< U_\alpha \, : \, \alpha \in \Phi^+ \setminus \Phi_J \right>$, and $L = L_{j_1 j_2 \ldots j_r} = \left< T, U_\alpha \, : \, \alpha \in \Phi_J \right>$ is the standard Levi complement. If the semisimple subgroup $L'$ has Lie type $\mathbf{X}$, for brevity we refer to $P$ as an ``$\mathbf{X}$-parabolic subgroup of $G$'', so for instance the parabolic subgroup of $G = E_{6}$ corresponding to the roots $\{\alpha_2,\alpha_3,\alpha_4,\alpha_5\}$ is a $D_{4}$-parabolic subgroup.

For a standard Levi subgroup $L$, we use the following notation. If $L_{0}$ is a simple component of $L$ then the simple roots for $L_{0}$ are a subset of $\Pi$, say $\Psi = \{\alpha_{1}',\ldots,\alpha_{m}'\}$. Order $\Pi$ according to $\alpha_{i} < \alpha_{j}$ if $i < j$. If $L_{0}$ has Lie type $A_{m}$, the embedding is chosen such that $\alpha_{1}'$ is the least simple root of $G$ contained in $\Psi$. If $L_{0}$ has type $E_{6}$ or $E_{7}$, then $\alpha_{i}' = \alpha_{i}$ for all $i$. If $L_{0}$ has type $D_{4}$ then $(\alpha_{1}',\alpha_{2}',\alpha_{3}',\alpha_{4}') = (\alpha_{2},\alpha_{4},\alpha_{3},\alpha_{5})$. If $L_{0}$ has type $D_{5}$ then $(\alpha_{1}', \alpha_{2}',\alpha_{3}',\alpha_{4}',\alpha_{5}') = (\alpha_{1},\alpha_{3},\alpha_{4},\alpha_{5},\alpha_{2})$ or $(\alpha_{6},\alpha_{5},\alpha_{4},\alpha_{3},\alpha_{2})$. If $L_{0}$ has type $D_{6}$ then $(\alpha_{1}',\alpha_{2}',\alpha_{3}',\alpha_{4}',\alpha_{5}',\alpha_{6}') = (\alpha_{7},\alpha_{6},\alpha_{5},\alpha_{4},\alpha_{3},\alpha_{2})$. Finally, if $L$ has multiple components of the same type, then these components are ordered according to the position of their least simple root ``$\alpha_{1}'$'' as an element of $\Pi$. For instance, if $G = E_{7}$ and $L = L_{12567}$ is a Levi subgroup of type $A_{1}A_{1}A_{3}$, then the first $A_{1}$ factor corresponds to $\alpha_1$, and the second to $\alpha_{2}$.

For $i \geq 1$ we define the subgroups
\[ Q(i) = \left< U_\gamma \, : \, \gamma = \sum_{\alpha \in \Pi}^{} c_{\alpha} \alpha, \sum_{\alpha \in \Pi \setminus J}^{} c_{\alpha} \geq i \right>. \vspace{-4pt} \]
The \emph{$i$-th level of $Q$} is $V_{i} \stackrel{\textup{def}}{=} Q(i)/Q(i+1)$, which is central in $Q/Q(i+1)$. By \cite[Theorem 2]{MR1047327} each level of $Q$ has the structure of a completely reducible $L$-module. The \emph{level} of a root $\gamma = \sum_{\alpha \in \Pi}c_{\alpha} \alpha$ is the sum $\sum_{\alpha \in \Pi \setminus J} c_{\alpha}$, and the \emph{height} of $\gamma$ is $\sum_{\alpha \in \Pi}c_{\alpha}$.

When $G$ has type $E_7$ or $E_8$, we will need to distinguish between certain isomorphic subsystem subgroups of $G$. In $E_7$ there are two conjugacy classes of Levi subgroups of type $A_5$, with representatives $A_{5} = L_{24567}$ and $A_{5}' = L_{34567}$, where our notation follows that of \cite[Table 8.2]{MR1329942}. These subgroups have connected centralisers $C_{G}(A_{5})^{\circ} = A_{2}$ and $C_{G}(A_{5}')^{\circ} = A_{1}T_{1}$ where $T_{1}$ is a $1$-dimensional torus. Furthermore $A_5'$ is contained in a subgroup $E_6$ whereas $A_5$ is not. In $E_8$ there are two conjugacy classes of subgroups $A_7$, with representatives the Levi subgroup $L_{1345678}$, which we denote by $A_{7}$, and a subgroup $A_{7}'$ which is a subsystem subgroup of a Levi subgroup $E_{7}$. Then $C_{G}(A_{7})^{\circ}$ is a $1$-dimensional torus, and $C_{G}(A_{7}')^{\circ} = A_{1}$.

Next, suppose that $G$ is of classical type or of type $G_2$, and that $X$ has type $A_1$ or $G_2$. Then the notation ``$X < G$ via $M$'' denotes an embedding $X \to G$ such that $V_{G}(\lambda_1) \downarrow X \cong M$. Since $p \neq 2$ throughout, this determines the image of $X$ up to $G$-conjugacy, unless $G$ has type $D_{n}$, in which case there are potentially two $G$-conjugacy classes with this action and we use this notation to refer to both simultaneously, specifying additional information when appropriate.

Now let $Y = H_1 H_2 \ldots H_k$ be a commuting product of simple subgroups $H_{i}$ all having the same type, $A_{1}$ or $G_{2}$. Then the simply-connected cover of $Y$ is $\hat{Y} \cong H \times H \ldots \times H$ ($k$ terms), where $H$ is simply-connected of type $A_1$ or $G_2$. We have a natural isogeny $\hat{Y} \to Y$, and a \emph{diagonal subgroup of $Y$} is the image in $Y$ of a subgroup of the form $\{ (\phi_1(h), \ldots, \phi_k(h)) \, : \, h \in H \}$, where each $\phi_i$ is an endomorphism of $H$. By \cite[Chapter 12]{steinberglectureson}, an endomorphism of $H$ is a product of an inner, graph and field morphism. Since $A_{1}$ and $G_{2}$ have no non-trivial graph automorphisms, to specify a diagonal subgroup of $Y$ up to $Y$-conjugacy it suffices to specify non-negative integers $r_{1},\ldots,r_{k}$, and we then take $\phi_{i} = F^{r_{i}}$ for each $i$. A diagonal subgroup $X$ of $Y$ is thus denoted by ``$X \hookrightarrow Y$ via $(\lambda_1^{[r_1]}, \lambda_1^{[r_2]},  \ldots, \lambda_1^{[r_k]})$''.

Taking this further, let $X$ have type $A_1$ or $G_2$, and let $Y = H_1 H_2 \ldots H_k$ ($k > 1)$ be semisimple, where each simple factor $H_{i}$ is classical or of type $G_2$. Then ``$X \hookrightarrow Y$ via $(M_1^{[r_1]}, \ldots, M_k^{[r_k]})$'' denotes a diagonal subgroup of $X_1 \ldots X_k$, with field twists $r_1, \ldots, r_k$, where each $X_i$ has the same type as $X$, and $X_i < H_i$ via $M_i$.

When discussing the centraliser of a subgroup $X$ of $G$, we will use the notation $U_{i}$ to denote an $i$-dimensional unipotent group, and $T_{j}$ to denote a $j$-dimensional torus. For instance $C_G(X)^{\circ} = U_{5} T_{1}$ means that $C_{G}(X)^{\circ}$ has a $5$-dimensional unipotent radical, with corresponding quotient a 1-dimensional torus.


\section{Preliminary Results} \label{sec:preliminaries}

\subsection{Exhibiting non-$G$-cr subgroups}

\vspace{-3pt}In this section we present preliminary results required for the proofs of Theorems \ref{thm:e6}--\ref{thm:e8} and their corollaries. The first of these limits the isomorphism types of non-$G$-cr simple subgroups occurring.

\begin{lemma} \label{lem:subtypes}
Let $G$ be an exceptional algebraic group in good characteristic $p$, and suppose $G$ has a non-$G$-cr simple subgroup of type $X$. Then $(G,X,p)$ is one of $(E_6,A_1,5)$, $(E_7,A_1,5 \text{ or } 7)$, $(E_7,G_2,7)$ or $(E_8,A_1 \text{ or } G_2,7)$.
\end{lemma}

\vspace{-10pt}\proof This follows immediately from \cite[Theorem 1]{MR1329942} and \cite[Theorem 1]{Stewart21072013}, the latter result ruling out the possibility $(E_7,A_2,5)$. \qed

Recall that our basic strategy for finding non-$G$-cr subgroups of $G$ is to iterate through the parabolic subgroups $P = QL$ of $G$, letting $X$ be an $L$-irreducible subgroup $A_{1}$ or $G_{2}$ of $L$, and then studying complements to $Q$ in the semidirect product $QX$. We derive an upper bound for the number of such complements (up to $G$-conjugacy), and it then remains to exhibit an explicit representative of each conjugacy class of subgroups. The next result allows us to quickly find such representatives.

Recall that a simple group $H$ of classical type is related by isogenies to a special linear, orthogonal or symplectic group, say $Cl(V)$. For $p \neq 2$, a subgroup $X$ of $Cl(V)$ is $Cl(V)$-cr if and only if $V$ is a completely reducible $X$-module (cf.\ Lemma \ref{lem:irredclass} below). Since complete reducibility is well-behaved with respect to isogenies \cite[Lemma 2.12]{MR2178661}, a subgroup of a simple algebraic group $H$ of classical type is $H$-completely reducible if and only if the corresponding subgroup of $Cl(V)$ acts completely reducibly on $V$. We now appeal to the following:

\begin{lemma}[{\cite[Theorem 3.26]{MR2178661}}] \label{lem:BMR}
Let $G$ be a simple algebraic group in good characteristic and let $M$ be a subsystem subgroup of $G$. Then a subgroup $H$ of $M$ is $G$-cr if and only if it is $M$-cr. 
\end{lemma}
Thus a non-completely reducible $X$-module $V$ gives rise to an embedding of $X$ into a classical group $H$ with non-$H$-cr image. If $H$ is a subsystem subgroup of $G$ this then gives rise to a non-$G$-cr subgroup of $G$.

In employing Lemma \ref{lem:BMR}, we will often refer to \emph{maximal connected} subgroups of $G$, that is, subgroups which are maximal among connected subgroups of $G$. The following result determines all reductive such subgroups up to $G$-conjugacy.

\begin{lemma}[{\cite[Corollary 2]{MR2044850}}] \label{lem:maximalexcep}
Let $G$ be a simple exceptional algebraic group in characteristic $p = 5$ or $7$. Let $M$ be a reductive, maximal connected subgroup of $G$. Then $M$ is $G$-conjugate to a subgroup in the following table, where each isomorphism type $X$ denotes a unique $G$-conjugacy class of subgroups. \vspace{-10pt}
\begin{longtable}{cl}
\\ \hline $G$ & $X$ \\ \hline
$G_2$ & $A_2$, $A_1 \tilde{A}_1$, $A_1$ $(p = 7)$ \\
$F_4$ & $B_4$, $A_1 C_3$, $A_1 G_2$, $A_2 \tilde{A}_2$, $G_2$ $(p = 7)$ \\
$E_6$ & $A_1 A_5$, $A_2^3$, $F_4$, $C_4$, $A_2 G_2$, $G_2$ $(p = 5$, two classes$)$, $A_2$ $($two classes$)$ \\
$E_7$ & $A_1 D_6$, $A_7$, $A_2 A_5$, $G_2 C_3$, $A_1 F_4$, $A_1 G_2$, $A_1 A_1$, $A_2$ \\
$E_8$ & $D_8$, $A_8$, $A_1 E_7$, $A_2 E_6$, $A_4^2$, $G_2 F_4$, $B_2$, $A_1 A_2$ \\ \hline
\end{longtable}

\end{lemma}

For $L$ a Levi subgroup of an exceptional algebraic group, the following results yield all $L$-irreducible subgroups of type $A_{1}$ or $G_{2}$.

\begin{lemma}[{\cite[Lemma 2.2]{MR2043006}}] \label{lem:irredclass}
Let $G = Cl(V)$ be a classical simple algebraic group in characteristic $p \neq 2$, and let $X$ be a $G$-irreducible semisimple subgroup of $G$. Then one of the following holds:
\begin{enumerate}[label=\normalfont(\roman*)]
\item $G = A_n$ and $X$ is irreducible on $V$. 
\item $G = B_n, C_n$ or $D_n$ and $V \downarrow X = V_1 \perp \ldots \perp V_k$ with the $V_i$ all non-degenerate, irreducible and inequivalent as $X$-modules.  
\end{enumerate}
\end{lemma}

\begin{lemma} \label{lem:A1modules}
Let $p > 2$ and let $X$ be simple of type $A_1$. Let $V = V_{X}(n)$ be a restricted irreducible $X$-module (so $n < p$). Then $X$ preserves a non-degenerate symplectic form on $V$ if $n$ is odd, and $X$ preserves a non-degenerate orthogonal form on $V$ if $n$ is even.  
\end{lemma}

\begin{proposition} \label{prop:irredA1}
Let $L_{0}$ be a simple factor of a proper Levi subgroup of $G = E_{6}$ or $E_{7}$ with $p = 5$, or of $G = E_{7}$ or $E_{8}$ with $p = 7$. The table below lists all simple $L_{0}$-irreducible subgroups $X$ of type $A_{1}$, up to $L_{0}$-conjugacy.
\end{proposition}

\begin{longtable}{>{\raggedright\arraybackslash}p{0.07\textwidth - 2\tabcolsep}>{\raggedright\arraybackslash}p{0.07\textwidth - 2\tabcolsep}>{\raggedright\arraybackslash}p{0.86\textwidth-\tabcolsep}@{}}

\caption{$L_0$-irreducible subgroups of type $A_{1}$} \\

\hline

$L_0$ & $p$ & Embedding of $X$ of type $A_1$ \\ \hline  

$A_1$ & $5,7$ & $X = L_0$ \\

$A_2$ & $5,7$ & $V_{A_2}(\lambda_1) \downarrow X = 2$ \\

$A_3$ & $5,7$ & $V_{A_3}(\lambda_1) \downarrow X = 3$ \\

& $5,7$ & $V_{A_3}(\lambda_1) \downarrow X = 1 \otimes 1^{[r]}$ $(r > 0)$ \\

$A_4$ & $5,7$ & $V_{A_4}(\lambda_1) \downarrow X = 4$ \\

$A_5$ & $7$ & $V_{A_5}(\lambda_1) \downarrow X = 5$ \\

& $5,7$ & $V_{A_5}(\lambda_1) \downarrow X = 2^{[r]} \otimes 1^{[s]}$ $(r \neq s$; $rs=0)$ \\

$A_6$ & $7$ & $V_{A_6}(\lambda_1) \downarrow X = 6$ \\

$D_4$ & $7$ & $V_{D_4}(\lambda_1) \downarrow X = 6 + 0$ \\

& $5,7$ & $V_{D_4}(\lambda_1) \downarrow X = 4^{[r]} + 2^{[s]}$ $(rs=0)$ \\

& $5,7$ & $V_{D_4}(\lambda_1) \downarrow X = 3^{[r]} \otimes 1^{[s]}$ $(r \neq s$; $rs=0)$ (two classes) \\

& $5,7$ & $V_{D_4}(\lambda_1) \downarrow X = 2^{[r]} + 1^{[s]} \otimes 1^{[t]} + 0$ ($s < t$; $rst=0$) \\

& $5,7$ & $V_{D_4}(\lambda_1) \downarrow X = 1 \otimes 1^{[r]} + 1^{[s]} \otimes 1^{[t]}$ ($r \neq 0$; $s \neq t$; $\{s,t\} \neq \{0,r\}$) \\

$D_5$ & $7$ & $V_{D_5}(\lambda_1) \downarrow X = 6^{[r]} + 2^{[s]}$ $(rs=0)$ \\

 & $5,7$ & $V_{D_5}(\lambda_1) \downarrow X = 4 + 4^{[r]}$ $(r > 0)$ \\

 & $5,7$ & $V_{D_5}(\lambda_1) \downarrow X = 4^{[r]} + 1^{[s]} \otimes 1^{[t]} + 0$ $(s<t$; $rs = 0)$ \\

 & $5,7$ & $V_{D_5}(\lambda_1) \downarrow X = 2^{[r]} + 2^{[s]} + 1^{[t]} \otimes 1^{[u]}$ $(r<s;$ $t \neq u;$ $rtu=0)$ \\

 & $5,7$ & $V_{D_5}(\lambda_1) \downarrow X = 2 + 2^{[r]} + 2^{[s]} + 0$ $(0 < r < s)$ \\

$D_6$ & $7$ & $V_{D_6}(\lambda_1) \downarrow X = 5^{[r]} \otimes 1^{[s]}$ $(r \neq s$; $rs=0)$ (two classes) \\

& $7$ & $V_{D_6}(\lambda_1) \downarrow X = 6^{[r]} + 4^{[s]}$ $(rs=0)$ \\

& $7$ & $V_{D_6}(\lambda_1) \downarrow X = 6^{[r]} + 1^{[s]} \otimes 1^{[t]} + 0$ $(s < t;$ $rs=0)$ \\

& $5,7$ & $V_{D_6}(\lambda_1) \downarrow X = 2^{[r]} \otimes 1^{[s]} \otimes 1^{[t]}$ ($r,s,t$ distinct; $s < t$; $rst=0$) (two classes) \\

& $5,7$ & $V_{D_6}(\lambda_1) \downarrow X = 2^{[r]} \otimes 2^{[s]} + 2^{[t]}$ $(r<s;$ $rt=0)$ \\

& $5,7$ & $V_{D_6}(\lambda_1) \downarrow X = 4^{[r]} + 2^{[s]} + 1^{[t]} \otimes 1^{[u]}$ $(t<u;$ $rst=0)$ \\

& $5,7$ & $V_{D_6}(\lambda_1) \downarrow X = 3^{[r]} \otimes 1^{[s]} + 1^{[t]} \otimes 1^{[u]}$ $(r \neq s; t \neq u;$ $rstu=0)$  \\

& $5,7$ & $V_{D_6}(\lambda_1) \downarrow X = 4^{[r]} + 2^{[s]} + 2^{[t]} + 0$ $(s<t;$ $rs=0)$ \\

& $5,7$ & $V_{D_6}(\lambda_1) \downarrow X = 3^{[r]} \otimes 1^{[s]} + 2^{[t]} + 0$ $(r \neq s;$ $rst=0)$  \\

& $5,7$ & $V_{D_6}(\lambda_1) \downarrow X = 2 + 2^{[r]} + 2^{[s]} + 2^{[t]}$ $(0 < r < s < t)$ \\

& $5,7$ & $V_{D_6}(\lambda_1) \downarrow X = 2^{[r]} + 1^{[s]} \otimes 1^{[t]} + 1^{[u]} \otimes 1^{[v]} + 0$ \\ & & $(s < t; u < v; s \leq u$; if $s=u$ then $t < v$; $rs=0)$ \\

& $5,7$ & $V_{D_6}(\lambda_1) \downarrow X = 1 \otimes 1^{[r]}  + 1^{[s]} \otimes 1^{[t]} + 1^{[u]} \otimes 1^{[v]}$ \\
& & ($r \neq 0$; $s \neq t$, $u \neq v$; $s \le \textup{min}\{t,u,v\}$; $\{0,r\}$, $\{s,t\}$, $\{u,v\}$ pairwise distinct.) \\

$D_7$ & $7$ & $V_{D_7}(\lambda_1) \downarrow X = 6 + 6^{[r]}$ $(r > 0)$ \\

& $7$ & $V_{D_7}(\lambda_1) \downarrow X = 6^{[r]} + 2^{[s]} + 2^{[t]} + 0$ $(s<t;$ $rs=0)$ \\

& $7$ & $V_{D_7}(\lambda_1) \downarrow X = 6^{[r]} + 2^{[s]} + 1^{[t]} \otimes 1^{[u]}$ $(t<u;$ $rst = 0)$ \\

& $7$ & $V_{D_7}(\lambda_1) \downarrow X = 4^{[r]} + 4^{[s]} + 1^{[t]} \otimes 1^{[u]}$ $(r<s;$ $t<u;$ $rt=0)$  \\

& $7$ & $V_{D_7}(\lambda_1) \downarrow X = 4^{[r]} + 4^{[s]} + 2^{[t]} + 0$ $(r<s;$ $rt=0)$  \\

& $7$ & $V_{D_7}(\lambda_1) \downarrow X = 4^{[r]} + 3^{[s]} \otimes  1^{[t]} + 0$ $(s \neq t;$ $rst=0)$  \\

& $7$ & $V_{D_7}(\lambda_1) \downarrow X = 4^{[r]} + 2^{[s]} + 2^{[t]} + 2^{[u]}$ $(s<t<u;$ $rs=0)$ \\

& $7$ & $V_{D_7}(\lambda_1) \downarrow X = 3^{[r]} \otimes 1^{[s]} + 2^{[t]} + 2^{[u]}$ $(r \neq s;$ $ t<u;$ $ rst=0)$ \\

& $7$ & $V_{D_7}(\lambda_1) \downarrow X = 4^{[r]} + 1^{[s]} \otimes 1^{[t]} + 1^{[u]} \otimes 1^{[v]} + 0$ \\
& & $(s<t$; $u<v$; $s \le u$; if $s = u$ then $t < v$; $rsu=0)$ \\

& $7$ & $V_{D_7}(\lambda_1) \downarrow X = 2^{[r]} + 2^{[s]} + 2^{[t]} + 1^{[u]} \otimes 1^{[v]} + 0$ $(r<s<t;$ $ u<v;$ $ ru=0)$ \\

& $7$ & $V_{D_7}(\lambda_1) \downarrow X = 2^{[r]} + 2^{[s]} + 1^{[t]} \otimes 1^{[u]} + 1^{[v]} \otimes 1^{[w]}$ \\
& & ($r<s$; $t<u$; $v<w$; $t \le v$; if $t = v$ then $u < w$; $rtv=0$) \\

$E_6$ & $7$ & $X \hookrightarrow \bar{A}_1 A_5$ via $(1^{[r]},5^{[s]})$ $(rs=0)$ \\

& $7$ & $X \hookrightarrow A_2 G_2$ via $(2^{[r]},6^{[s]})$ $(r \neq s;$ $ rs=0)$ \\

& $5,7$ & $X \hookrightarrow \bar{A}_1 A_5$ via $(1^{[r]},2^{[s]} \otimes 1^{[t]})$ $(s \neq t;$ $ rst=0)$ \\

& $5,7$ & $X \hookrightarrow A_2^3$ via $(2,2^{[r]}, 2^{[s]})$ $(0<r<s)$ \\

$E_7$ & $7$ & $X < \bar{A}_1 D_6$ see \cite[Theorem 4]{Tho2} for explicit classes   \\

& $7$ & $X < A_1 A_1$ via $(1^{[r]},1^{[s]})$ $(r \neq s;$ $ rs=0)$ \\

& $7$ & $X < A_1 G_2$ via $(1^{[r]},6^{[s]})$ $(r \neq s;$ $ rs=0)$ \\

& $7$ & $X < G_2 C_3$ via $(6^{[r]},5^{[s]})$ $(r \neq s;$ $rs=0)$ \\

& $7$ & $X < G_2 C_3$ via $(6^{[r]},2^{[s]} \otimes 1^{[t]})$ $(s \notin \{r,t\}$; $rst=0)$ \\ 

$D_5$ & $5,7$ & $V_{D_5}(\lambda_1) \downarrow X = 2 \otimes 2^{[r]} + 0$ $(r > 0)$ \\ \hline
\end{longtable}

\begin{proof}
For $L_0 = E_6$ and $E_7$, this follows from \cite[Theorems 3 and 4]{Tho2}. For $L_0$ of type $A_{n}$, this follows from Lemmas \ref{lem:irredclass} and \ref{lem:A1modules}. For $L_{0}$ of type $D_{n}$, Lemmas \ref{lem:irredclass} and \ref{lem:A1modules} imply that the given subgroups of type $A_{1}$ are distinct and unique up to conjugacy as a subgroup of $GO_{2n}(K) = SO_{2n}(K){\left<\tau\right>}$, where $\tau$ induces a graph automorphism on $SO_{2n}(K)$. It remains to determine whether such a conjugacy class of subgroups splits into two classes of subgroups of $SO_{2n}(K)$. This is equivalent to determining whether $N_{GO_{2n}(K)}(X) \le SO_{2n}(K)$. Since $X = A_{1}$ has no non-trivial outer algebraic automorphisms, there are precisely two possibilities: Either $X$ lies in the centraliser of some element of $GO_{2n}(K) \setminus SO_{2n}(K)$, which by \cite[Table 4.3.1]{MR1490581} is the stabiliser $B_{k}B_{n-k-1}$ ($n/2 \le k \le n - 1$) of a direct-sum decomposition of the natural orthogonal module; or $X$ is irreducible on the natural module, and is not normalised by any element in $GO_{2n}(K) \setminus SO_{2n}(K)$, and so gives rise to two classes of subgroups in $SO_{2n}(K)$ by the orbit-stabiliser theorem. 
\end{proof}

\newpage

\begin{proposition} \label{prop:irredG2}
Let $L_{0}$ be a simple factor of a proper Levi subgroup of $G = E_{7}$ or $E_{8}$ with $p = 7$. The table below lists all simple $L_{0}$-irreducible subgroups $X$ of type $G_{2}$, up to $L_{0}$-conjugacy.
\end{proposition}

\begin{longtable}{>{\raggedright\arraybackslash}p{0.07\textwidth - 2\tabcolsep}>{\raggedright\arraybackslash}p{0.60\textwidth-\tabcolsep}@{}}
\caption{$L_0$-irreducible subgroups of type $G_{2}$} \\
\hline $L_0$ & Embedding of $X$ of type $G_2$ \\ \hline
$A_6$ & $X < A_6$ via $10$ \\
$D_4$ & $X < D_4$ via $10 + 00$ \\
$D_7$ & $X < D_7$ via $01$ (two classes) \\
$E_6$ & $V_{27} \downarrow X = 20 + 00$, $X$ maximal in $F_4$ \\ \hline
\end{longtable}

\begin{proof}
For $L_0$ exceptional, the result follows from \cite[Theorems 2 and 3]{Tho1}. For $L_0$ classical, this is similar to the previous proposition, using Lemma \ref{lem:irredclass} and \cite{MR1901354}, and easier since $G_{2}$ has far fewer modules of low dimension.
\end{proof}

\subsection{Complements and non-abelian cohomology} \label{sec:nonab}

In this section, we describe a method for classifying complements in parabolic subgroups by approximating certain non-abelian cohomology sets. The technique was pioneered by D.\ Stewart in his Ph.\ D.\ thesis, and our definitions and strategy are taken from \cite[Section 3.2]{MR3075783}.

Recall that if $X$ and $Q$ are algebraic groups over $K$, with a morphism $X \times Q \to Q$ giving an action of $X$ on $Q$, then complements to $Q$ in the semidirect product $QX$ correspond bijectively with rational \emph{1-cocycles}, which are variety morphisms $\phi \, : \, X \to Q$ such that $\phi(xy) = \phi(x)(x.\phi(y))$ for all $x,y \in X$. Here, a complement $X'$ to $Q$ is a closed subgroup of $QX$ satisfying (i) $QX' = QX$, (ii) $Q \cap X' = 1$, and (iii) $L(Q) \cap L(X') = \{0\}$ (cf.\ \cite[4.3.1]{MR2753264}). By \cite[Lemma 3.6.1]{MR3105754}, a subgroup satisfying (i) and (ii) automatically satisfies (iii) when $X$ is connected reductive, $Q$ is unipotent and $p \neq 2$. 

Two cocycles $\phi$, $\psi$ are \emph{cohomologous} if there exists $q \in Q$ such that $\phi(x) = q^{-1}\psi(x)(x.q)$ for all $x \in X$. This defines an equivalence relation on the set $Z^{1}(X,Q)$ of 1-cocycles, and the corresponding quotient is called the \emph{cohomology set} $H^{1}(X,Q)$, which parametrises complements up to $Q$-conjugacy. The set $H^{1}(X,Q)$ has a distinguished element, denoted $[0]$, which is the class of the map sending every element of $X$ to the identity of $Q$. In general, $H^{1}(X,Q)$ is only a pointed set, but if $Q$ is a $KX$-module, then $H^{1}(X,Q)$ is a $K$-vector space in a natural way.

In our calculations $X$ will always be simple, and $Q$ will always be connected and unipotent, with a filtration by $X$-stable connected subgroups $Q = Q(1) \rhd Q(2) \rhd \ldots \rhd Q(r+1) = 1$ such that each section $V_{i} \stackrel{\textup{def}}{=} Q(i)/Q(i+1)$ is a rational $KX$-module which is central in $Q/Q(i+1)$. This allows us to study $H^{1}(X,Q)$ in terms of the vector-space direct sum
\[ \mathbb{V} = \mathbb{V}_{X,Q} \stackrel{\textup{def}}{=} \bigoplus_{i = 1}^{r} H^{1}(X,V_{i}) \]
and then appeal to the representation theory of the simple group $X$.

We first recall some results from non-abelian cohomology, using \cite[\S I.5]{MR1867431} as a standard reference. If $R$ is an $X$-stable central subgroup of $Q$, then the short exact sequence
\[ 0 \to R \to Q \to Q/R \to 0 \]
gives rise to a long exact sequence of pointed sets:
\label{eqn:exact}
\begin{align*}
0 \to R^{X} \to Q^{X} \to (Q/R)^{X} \to H^{1}(X,R) \to H^{1}(X,Q) \to H^{1}(X,Q/R) \to H^{2}(X,R)
\end{align*}
where $-^{X}$ denotes a fixed-point subgroup, and where $H^{2}(X,R)$ is the usual second abelian cohomology group, defined for example in \cite[\S II.4.2]{MR2015057}. Since $R$ is central in $Q$, the group $H^{1}(X,R)$ acts on $H^{1}(X,Q)$ on the left; for $[\phi] \in H^{1}(X,R)$ and $[\psi] \in H^{1}(X,Q)$ we have $[\phi].[\psi] = [\phi \cdot \psi]$, where $\phi \cdot \psi(x) \stackrel{\textup{def}}{=} \phi(x)\psi(x)$ for all $x \in X$. Since we shall use this in Lemma \ref{lem:rho} below, we note that the condition for a group action becomes $[\phi].([\phi'].[\psi]) = [\phi + \phi'].[\psi]$ if $H^{1}(X,R)$ is written additively.

\begin{lemma}[{\cite[{\S I.5.7, Proposition 42}]{MR1867431}}] \label{lem:diff_lifts}
With $X$, $Q$ and $R$ as above, two elements of $H^{1}(X,Q)$ have the same image in $H^{1}(X,Q/R)$ if and only if they lie in the same $H^{1}(X,R)$-orbit.
\end{lemma}

\begin{definition}[{cf.\ \cite[Definition 3.2.5]{MR3075783}}]
For $i = 1,\ldots,r$ we define a partial map $\rho_{i} \, : \, \mathbb{V} \to H^{1}(X,Q/Q(i+1))$ as follows. For $i = 1$, set $\rho_{1}([\gamma_{1}],\ldots,[\gamma_{r}]) = [\gamma_{1}]$. For $i > 1$, if $\rho_{i-1}([\gamma_{1}],\ldots,[\gamma_{r}])$ is defined and lifts to some element $[\Gamma]$ under the natural map $H^{1}(X,Q/Q(i+1)) \to H^{1}(X,Q/Q(i))$, then set $\rho_{i}([\gamma_{1}],\ldots,[\gamma_{r}]) = [\gamma_{i}].[\Gamma]$, otherwise declare $\rho_{i}$ undefined at $([\gamma_{1}],\ldots,[\gamma_{r}])$.

We set $\rho = \rho_{r}$. Note that $\rho_{i}([\gamma_{1}],\ldots,[\gamma_{r}])$ depends only on the first $i$ coordinates of its argument, but does depend on the choice of lifts made. We pick lifts to be consistent with the convention that $\rho_{i}([0],[0],\ldots,[0]) = [0]$ for all $i$.
\label{def:rho}
\end{definition}

\begin{lemma}[{cf.\ \cite[Proposition 3.2.6]{MR3075783}}] \label{lem:rho}
Each partial map $\rho_{i} \, : \, \mathbb{V} \to H^{1}(X,Q/Q(i+1))$ is surjective.
\end{lemma}
\proof Proceed by induction on $i$. The result is trivial for $i = 1$. If $i > 1$, then a class $[\phi] \in H^{1}(X,Q/Q(i+1))$ gives a class $[\phi'] \in H^{1}(X,Q/Q(i))$ under the natural map. By induction, there exists $\mathbf{v} = ([\gamma_1],[\gamma_2],\ldots,[\gamma_r]) \in \mathbb{V}$ such that $\rho_{i-1}(\mathbf{v}) = [\phi']$. Since $[\phi']$ has a lift in $H^{1}(X,Q/Q(i+1))$, namely $[\phi]$, the class $\rho_{i}(\mathbf{v})$ is defined. Now, $\rho_{i}(\mathbf{v})$ and $[\phi]$ have the same image in $H^{1}(X,Q/Q(i))$, hence by Lemma \ref{lem:diff_lifts} there exists $[\delta] \in H^{1}(X,V_{i})$ such that $[\phi] = [\delta].\rho_{i}(\mathbf{v})$. But by definition of $\rho_{i}$, the right-hand side is equal to $[\delta].([\gamma_{i}].[\Gamma]) = [\delta + \gamma_{i}].[\Gamma]$, where $[\Gamma]$ is the chosen fixed lift of $[\phi']$ to $H^{1}(X,Q/Q(i+1))$. Then we have $[\delta + \gamma_{i}].[\Gamma] = \rho_{i}([\gamma_{1}],[\gamma_{2}],\ldots,[\gamma_{i-1}],[\delta + \gamma_{i}],[\gamma_{i+1}],\ldots,[\gamma_{r}])$, hence $[\phi]$ lies in the image of $\rho_{i}$. \qed

As an immediate consequence, we obtain the following.
\begin{corollary}
If $H^{1}(X,V_{i}) = 0$ for each $i$, then $H^{1}(X,Q) = 0$. In this case, every complement to $Q$ in $QX$ is $Q$-conjugate to $X$.
\label{cor:nonzeroH1}
\end{corollary}

Since the map $\rho$ is a surjection from a subset of $\mathbb{V}$ to $H^{1}(X,Q)$, choosing a basis of $\mathbb{V}$ allows us to parametrise conjugacy classes of complements to $Q$ in $QX$ by certain ordered $m$-tuples $(k_1,\ldots,k_m)$ of elements of $K$, for some $m \ge 0$. We denote by $X_{[k_1,k_2,\ldots,k_m]}$ a fixed complement to $Q$ in $QX$ corresponding to $(k_1,k_2,\ldots,k_m) \in \mathbb{V}$, when $\rho$ is defined at this point.

We now consider the question of when $\rho(\mathbf{v}) = \rho(\mathbf{w})$ for $\mathbf{v}$, $\mathbf{w} \in \mathbb{V}$. If $Q = R_u(P)$ for a parabolic subgroup $P$ of $G$, the Chevalley commutator relations give us precise information about the $X$-group structure of $Q$, which in turn lets us derive information about the set $H^{1}(X,Q)$.

For each $i$ and $j$, we have the containment $[Q(i),Q(j)] \subseteq Q(i+j)$, where $Q(i) = 1$ for all $i > r$. Thus for each $v \in V_{j}$, we get a map $Q(i) \to Q(i+j)$, sending $q \in Q(i)$ to $[\hat{v},q]$, where $\hat{v}$ is a fixed lift of $v$ to $Q(j)$. Composition with the quotient $Q(i+j) \to V_{i+j}$ gives a map $c_{v} = c_{v,i} \, : \, Q(i) \to V_{i+j}$ which is independent of the choice of lift $\hat{v}$. From standard properties of commutators,
\begin{align*}
c_{v}(uw) &= [\hat{v},uw]Q(i+j+1)\\
&= [\hat{v},u]\left({}^{u}[\hat{v},w]\right)Q(i+j+1)
\end{align*}
and since $V_{i+j}$ is central in $Q/Q(i+j+1)$, the above is equal to $c_{v}(u)c_{v}(w)$, so $c_{v}$ is a group homomorphism. In addition, if $v \in V_{j}$ is fixed by $X$, then we claim $c_{v}$ is $X$-equivariant. Indeed, if $u \in Q(i)$, then for all $x \in X$ we have
\begin{align*}
c_{v}({}^{x}u) &= [\hat{v},{}^{x}u]Q(i+j+1)\\
&= [({}^{x}\hat{v})v',{}^{x}u]Q(i+j+1)
\end{align*}
for some $v' \in Q(j+1)$. Since $v'$ commutes with $Q(i)$ modulo $Q(i+j+1)$, the above equals ${}^{x}[\hat{v},u]Q(i+j+1) = {}^{x}c_{v}(u)$. Moreover, since $Q(i+1)$ lies in the kernel of $c_{v,i}$, this descends to a homomorphism $V_{i} \to V_{i+j}$ of $X$-modules.

Our next proposition is rather technical, so we take a moment to discuss its use in proving our main theorems. Informally, it states that if $U$ and $W$ are summands occurring in levels of $Q$, with $H^{1}(X,U) \cong H^{1}(X,W) \cong K$, and if $x \in Q^{X}$ induces an isomorphism $U \to W$, then the partial maps $\rho_{i}$ are surjective on the subset of $\mathbb{V}$ which gives a cocycle on at least one of $U$ or $W$. If elements of $\mathbb{V} = H^{1}(X,U) \oplus H^{1}(X,W) \oplus \ldots$ are given by tuples $(k_1,k_2,\ldots)$, then this allows us to assume that $k_1 k_2 = 0$.

\begin{proposition} \label{prop:trivs}
With $X$ and $Q$ as above, suppose that $v \in V_{j}^{X}$ for some $j < r$, and let $n$ be minimal such that $c_{v,n} \, : \, V_{n} \to V_{n+j}$ is non-zero. Suppose that the subspace $\left<v\right>$ lifts to a $1$-dimensional subgroup of $(Q/Q(n+j+1))^{X}$. Let $U$ and $W$ be direct summands of $V_{n}$ and $V_{n+j}$, respectively, such that $c_{v,n} \, : \, V_{n} \to V_{n+j}$ restricts to an isomorphism $U \to W$. Suppose also that $H^{1}(X,U) \cong K$.

Pick a basis $\{e_{l} \, : \, l = 1,\ldots,m\}$ of $\mathbb{V}$, where each vector $e_{l}$ lies in the image of some inclusion $H^{1}(X,V_{i}) \to \mathbb{V}$, and where $e_{1}$ and $e_{2}$ respectively lie in the images of $H^{1}(X,U)$ and $H^{1}(X,W)$. Finally, define
\[ \mathbb{V}_{0} = \left\{ \sum_{l=1}^{m} t_{l}e_{l} \, : \, t_{l} \in K,\ t_{1} t_{2} = 0 \right\}. \]
Then for each $i \ge n + j$, the restricted partial map $\rho_{i} \, : \, \mathbb{V}_{0} \to H^{1}(X,Q/Q(i+1))$ is surjective.
\end{proposition}

\proof Note first that by our choice of basis, if $i \neq n,\ n + j$ then the image of $H^{1}(X,V_{i})$ is spanned by some vectors $e_{l}$ with $l \notin \{1,2\}$. In particular, if $w \in \mathbb{V}_{0}$ and if $w'$ lies in the image of $H^{1}(X,V_{i})$ with $i \neq n,\ n + j$, then $w + w' \in \mathbb{V}_{0}$. We will make use of this shortly.

By Lemma \ref{lem:rho}, if $1 < i \le r$ then every element of $H^{1}(X,Q/Q(i+1))$ is of the form $\rho_{i}(\mathbf{v}) = [\phi].\rho_{i-1}'(\mathbf{v})$ for some $[\phi] \in H^{1}(X,V_{i})$ and some $\mathbf{v} \in \mathbb{V}$, where $\rho_{i-1}'(\mathbf{v}) \in H^{1}(X,Q/Q(i+1))$ is the fixed lift of $\rho_{i-1}(\mathbf{v}) \in H^{1}(X,Q/Q(i))$.

For an induction on $i$, suppose that $r > i-1 \ge n + j$, and suppose that the restriction $\rho_{i-1} \, : \, \mathbb{V}_{0} \to H^{1}(X,Q/Q(i))$ is surjective. Let $\rho_{i}(\mathbf{v}) \in H^{1}(X,Q/Q(i+1))$, and write $\rho_{i-1}(\mathbf{v}) = \rho_{i-1}(\mathbf{w})$ for some $\mathbf{w} \in \mathbb{V}_{0}$. Then $\rho_{i-1}'(\mathbf{v}) = \rho_{i-1}'(\mathbf{w})$ since the choice of lift is fixed, and thus $\rho_{i}(\mathbf{v}) = [\phi].\rho_{i-1}'(\mathbf{v}) = [\phi].\rho_{i-1}'(\mathbf{w})$ for some $[\phi] \in H^{1}(X,V_{i})$. Moreover, by definition of $\rho_{i}$ we have $[\phi].\rho_{i-1}'(\mathbf{w}) = \rho_{i}([\phi] + \mathbf{w})$, where on the right-hand side we have identified $[\phi]$ with its image under the natural inclusion $H^{1}(X,V_{i}) \to \mathbb{V}$. By the observation in the first paragraph above, we have $[\phi] + \mathbf{w} \in \mathbb{V}_{0}$, and the restriction of $\rho_{i}$ to $\mathbb{V}_{0}$ is surjective, as required.

It remains to prove that $\rho_{n+j} \, : \, \mathbb{V}_{0} \to H^{1}(X,Q/Q(n+j+1))$ is surjective. For ease of notation we may replace $Q$ with $Q/Q(n+j+1)$ without loss of generality, so that $\rho = \rho_{r} = \rho_{n+j}$.

Suppose that $\rho(\mathbf{v})$ is defined, where $\mathbf{v} = \sum_{l=1}^{m} t_{l}e_{l}$ and $t_{1}t_{2} \neq 0$. By minimality of $n$, conjugation by a lift $\hat{v} \in Q^{X}$ of $v \in V_{j}^{X}$ sends a complement corresponding to $[\phi] \in H^{1}(X,Q)$ to a complement corresponding to $[c_{v,n} \circ \phi].[\phi]$, where $[c_{v,n} \circ \phi] \in H^{1}(X,Q(n+j))$. Therefore these cohomology classes correspond to the same $Q$-conjugacy class of complements to $Q$ in $QX$, and hence are equal. Since $W$ is an $X$-module direct summand of $Q(n+j)$, there is a natural projection $Q(n+j) \to W$, and the image of $[c_{v,n} \circ \phi]$ under the induced map $H^{1}(X,Q(n+j)) \to H^{1}(X,W)$ is a multiple of $e_{2}$. Moreover this multiple must be non-zero, since by hypothesis the map $c_{v,n}$ induces an isomorphism $U \to W$ and therefore induces an isomorphism $H^{1}(X,U) \to H^{1}(X,V)$. Scaling $e_{2}$ if necessary, we can assume that that this image is $-t_{2} e_{2}$. So if $[\phi] = \rho(\mathbf{v})$, then $[\phi] = [c_{v,n} \circ \phi].[ \phi ] = [-t_{2} e_{2} ] . \rho(\mathbf{v}) = \rho(\mathbf{v} - t_{2}e_{2})$, and the vector $\mathbf{v} - t_{2}e_{2}$ lies in $\mathbb{V}_{0}$, which proves the desired result. \qed

As a special case, we obtain the following.
\begin{corollary} \label{cor:trivs}
With $X$ and $Q$ as above, suppose that $\mathbb{V} = H^{1}(X,U) \oplus H^{1}(X,W)$, where $U$ and $W$ are direct summands in levels of $Q$ with $H^{1}(X,U) \cong H^{1}(X,W) \cong K$. Fix a basis of $\mathbb{V}$ consisting of non-zero elements from $H^{1}(X,U)$ and $H^{1}(X,W)$.

Suppose that $Q^{X}$ contains a $1$-dimensional subgroup inducing isomorphisms $U \to W$. Then the restriction $\rho \, : \, \{ (k_1,k_2) \in \mathbb{V} \, : \, k_1 k_2 = 0\} \to H^{1}(X,Q)$ is surjective.
\end{corollary}

\subsection{Representations and abelian cohomology}

In light of the preceding section, we now wish to describe $H^{1}(X,V)$ for various $X$-modules $V$, when $X$ is simple of type $A_{1}$ or $G_{2}$. It will also be useful to understand $H^{2}(X,V)$ in some cases, since if $X$ lies in a parabolic subgroup $P = QL$ and $H^{2}(X,V_{i}) = 0$ for some level $V_{i}$, then it follows from the long exact sequence that every element of $H^{1}(X,Q/Q(i))$ lifts to an element of $H^{1}(X,Q/Q(i+1))$.

\begin{lemma}[{\cite[Corollary 3.9]{MR684821}}] \label{lem:h1fora1}
Let $X$ be simple of type $A_{1}$ and let $M$ be an irreducible $X$-module. Then $H^{1}(X,M) \neq 0$ if and only if $M$ is a Frobenius twist of $(p-2) \otimes 1^{[1]}$; in this case $H^1(X,M) \cong K$.
\end{lemma}

\begin{lemma}[{\cite[Theorem 1]{MR2557160}}] \label{lem:h2fora1}
Let $X$ be simple of type $A_1$ and let $M$ be an irreducible $X$-module. Then $H^2(X,M) \neq 0$ if and only if $M$ is a Frobenius twist of $V_{X}(r)$ where r is $2p$, $2p^2 - 2p - 2$ or $2p-2+(2p-2)p^e$ $(e > 1)$; in this case $H^{2}(X,M) \cong K$.
\end{lemma}

The next result is a special case of results due to Cline, Parshall, Scott and van der Kallen \cite[Corollaries 3.9, 3.10]{MR0439856}. Recall that $W_{X}(\lambda)$ denotes the Weyl module for $X$ of highest weight $\lambda$, and additionally let $M(\lambda)$ denote the unique maximal submodule of $W_X(\lambda)$.

\begin{lemma} \label{lem:cps}
Let $X$ be a simple algebraic group. Then for any dominant weight $\lambda$ there are isomorphisms $H^{2}(X,V_{X}(\lambda)) \cong H^{1}(X,M(\lambda)^{\ast})$ and $H^{1}(X,V_{X}(\lambda)) \cong H^{0}(X,M(\lambda)^{\ast})$.
\end{lemma}
It follows that $H^{2}(X,V_{X}(\lambda_{i})) = 0$ for any simple algebraic group $X$ and every fundamental dominant weight $\lambda_{i}$. We now give a description of certain Weyl modules and tilting modules for $X$ of type $G_{2}$, to which we will refer later on.

\begin{lemma} \label{lem:G2mods}
Let $X$ be simple of type $G_2$ in characteristic $7$, and let $\lambda$ be one of $10$, $20$, $11$, $01$, $30$. Then the Weyl modules $W(\lambda)$ and tilting modules $T(\lambda)$ have the following structure:
\begin{enumerate}[label=\normalfont(\roman*)]
\item $W(10) = T(10) = 10$,
\item $W(01) = T(01) = 01$,
\item\label{G2mod20} $W(20) = 20|00$, $T(20) = 00 | 20 | 00$,
\item $W(11) = 11 | 20$,  $T(11) = 20 | (11 + 00) | 20$,
\item $W(30) = T(30) = 30$.
\end{enumerate}
\end{lemma}

\proof The composition factors of the Weyl modules are well-known, see for instance \cite{MR1901354}. The submodule structure of these, as well as the composition factors and submodule structure of each $T(\lambda)$, follows easily from the fact that $W(\lambda)/\textup{soc}(W(\lambda)) \cong V_{X}(\lambda)$ for each $\lambda$, and that $T(\lambda)$ admits both a filtration by Weyl modules and a filtration by duals of Weyl modules. \qed

In Sections \ref{sec:e6}--\ref{sec:e8} we will make implicit use of the following results on tilting modules:
\begin{lemma} \label{lem:tiltingnoH1}
Let $X$ be an algebraic group and $\lambda$ be a dominant weight for $X$. Then:
\begin{enumerate}[label=\normalfont(\roman*)]
\item A direct summand of a tilting module is tilting,
\item The tensor product of two tilting modules is tilting,
\item If $p > r$ then the $r$-th symmetric and alternating powers of a tilting module are tilting,
\item $H^1(X,T(\lambda)) = 0$.
\end{enumerate}
\end{lemma}
\proof Parts (i), (ii) and (iv) are well-known, see for instance \cite[\S{E.1}, E.2, E.7]{MR2015057}. For part (iii), if $p > r$ then the symmetric power $S^{r}(V)$ and alternating power $\bigwedge^{r}(V)$ can each be realised as the image of a projection operator $V^{\otimes r} \to V^{\otimes r}$ (see for instance \cite[\S 11.5, Proposition 40]{MR2286236}), and therefore as a direct summand of $V^{\otimes r}$. Now apply parts (i) and (ii). \qed

\begin{lemma} \label{lem:tiltingform}
Let $X$ be a simple algebraic group of type $A_{1}$, and let $V = T(n)$ be a tilting module for $X$, where $n > 0$. Then $V$ supports a non-degenerate $X$-invariant bilinear form, which is symmetric if $n$ is even, and skew-symmetric otherwise. 
\end{lemma}

\proof The result holds for $n = 1$ since $T(1)$ is the natural $2$-dimensional module for $X = SL_{2}(K) = Sp_{2}(K)$. Now let $n > 1$, and for an induction, assume that the result holds for all integers $1 \le m < n$. If $n$ is even then let $W = T(n/2) \otimes T(n/2)$, otherwise let $W = T(\frac{n-1}{2}) \otimes T(\frac{n+1}{2})$. Then $W$ supports the non-degenerate tensor product form, which is symmetric in the first case and skew-symmetric in the second. Moreover, $T(n)$ occurs as a direct summand of $W$, and is the unique indecomposable summand containing a vector of weight $n$. Hence the natural isomorphism $W/T(n)^{\perp} \to T(n)^{\ast}$, and the fact that $T(n) \cong T(n)^{\ast}$, shows that $T(n)^{\perp} \cap T(n) = 0$. Thus $T(n)$ is a non-degenerate subspace of $W$, and the result follows. \qed

\begin{lemma} \label{lem:twoclasses}
Let $X$ be a simple algebraic group of type $A_{1}$ or $G_{2}$ with an indecomposable orthogonal module $V$ of dimension $2n$. Then there are precisely two conjugacy classes of subgroups of $D_{2n}$ isomorphic to $X$ acting on $V_{D_{2n}}(\lambda_1)$ via $V$. 
\end{lemma}

\proof Since neither $A_1$ nor $G_2$ have any non-trivial outer algebraic automorphisms, the result follows from the discussion in the proof of Proposition \ref{prop:irredA1}. 

\subsection{Explicit cohomology calculations for $A_1$}

When $X$ has type $A_{1}$, we sometimes adopt a computational approach to studying $H^{1}(X,Q)$. Lemma \ref{lem:a1rels} below lets us verify that a given subset of $QX$ indeed generates a subgroup of type $A_{1}$. Lemma \ref{lem:a1cocycle} gives us an explicit formula for a cocycle $X \to V = 1^{[1]} \otimes (p - 2)$ when $p > 2$, which restricts the possibilities for a cocycle $X \to Q/[Q,Q]$.

\begin{lemma}[{\cite[12.1.1 and Remark p.198]{MR0407163}}] \label{lem:a1rels}
Let $K$ be any field and let $X$ be a group generated by $\{x_{+}(t),\, x_{-}(t) \, : \, t \in K \}$, with relations
\begin{enumerate}[label=\normalfont(\roman*)]
\item $x_{\pm}(t_{1})x_{\pm}(t_{2}) = x_{\pm}(t_{1} + t_{2})$,
\item $h(t)h(u) = h(tu)$,
\item $n(t)x_{+}(t_1)n(t)^{-1} = x_{-}(-t^{-2}t_{1})$,
\end{enumerate}
for all $t_{1}$, $t_{2} \in K$ and all $t$, $u \in K^{*}$, where $n(t) = x_{+}(t)x_{-}(-t^{-1})x_{+}(t)$ and $h(t) = n(t)n(-1)$.

Then $X$ is abstractly isomorphic to $SL_{2}(K)$ or $PSL_{2}(K)$.
\end{lemma}

Now let $X$ be as in Lemma \ref{lem:a1rels}, where $K$ has characteristic $p > 2$. Let $T = \{h(t) \, : \, t \in K^{*}\}$ be a maximal torus of $X$, with corresponding root subgroups $U_{\pm} = \{x_{\pm}(t) \, : \, t \in K\}$. Let $K[v_{1},v_{2}]$ be a polynomial algebra in two variables, with $X$ acting as algebra automorphisms induced via:
\[ \begin{array}{rlcrl}
x_{+}(t).v_{1} =& v_{1}, &\quad& x_{+}(t).v_{2} =& v_{2} + tv_{1},\\
x_{-}(t).v_{1} =& v_{1} + tv_{2}, &\quad& x_{-}(t).v_{2} =& v_{2}.
\end{array} \]
View $(p-2)$ as the $KX$-submodule of all homogeneous polynomials of degree $p-2$ in $v_{1}$ and $v_{2}$, and view $1^{[1]}$ as the submodule spanned by $f_{1} \stackrel{\textup{def}}{=} v_{1}^{p}$ and $f_{-1} \stackrel{\textup{def}}{=} v_{2}^{p}$. Then the basis
\[ \left\{v_{1}^{r-1}v_{2}^{p-r-1} \otimes f_{\pm 1} \, : \, 0 < r \le p - 1 \right\}\]
of $V = (p - 2) \otimes 1^{[1]}$ consists of weight vectors. For each $0 < r \le p - 1$ let $e_{2r}$ denote the vector $v_{1}^{r-1}v_{2}^{p-r-1} \otimes f_{1}$ and let $e_{-2r}$ denote $v_{1}^{p-r-1}v_{2}^{r-1} \otimes f_{-1}$. Clearly each $e_{\pm2r}$ has weight $\pm 2r$, and the elements $x_{\pm}(t)$ act as follows:
\begin{align*}
x_{+}(t).(v_{1}^{a}v_{2}^{p-2-a}) &= v_{1}^{a}(tv_{1} + v_{2})^{p-2-a} = \sum_{j = 0}^{p-2-a} \binom{p-2-a}{j}t^{p-2-a-j}v_{1}^{p-2-j}v_{2}^{j}\\
x_{-}(t).(v_{1}^{a}v_{2}^{p-2-a}) &= (v_{1} + tv_{2})^{a}v_{2}^{p-2-a} = \sum_{j = 0}^{a} \binom{a}{j}t^{a-j}v_{1}^{j}v_{2}^{p-j-2}
\end{align*}
and from this, for any $r > 1$ it follows that
\label{eqn:xplus}
\begin{align*}
x_{+}(t).e_{2r} &= \left(\sum_{j = 0}^{p-r-1} \binom{p-r-1}{j}t^{p-r-1-j}v_{1}^{p-2-j}v_{2}^{j}\right) \otimes f_{1}\\
&= \sum_{m = r}^{p-1} \binom{p-r-1}{m-r}t^{m-r}e_{2m}.
\end{align*}

Similarly, we find that
\begin{align*}
x_{+}(t).e_{-2r} &= \left(\sum_{j = 0}^{r-1} \binom{r-1}{j}t^{r-1-j}v_{1}^{p-j-2}v_{2}^{j}\right) \otimes \left(f_{-1} + t^{p}f_{1}\right)\\
&= \left( \sum_{s = 1}^{r} \binom{r-1}{s-1} t^{r-s}e_{-2s}\right)
+ \left( \sum_{s = p-r}^{p-1} \binom{r-1}{s+r-p} t^{r+s} e_{2s} \right)
\end{align*}
with similar expressions for $x_{-}(t).e_{\pm2r}$. As expected, it also follows that
\begin{align*}
n(t).e_{\pm 2r} &= (-t^{\mp 2})^{r}e_{\mp 2r},\\
h(t).e_{\pm 2r} &= t^{\pm 2r}e_{\pm 2r}.
\end{align*}

\begin{lemma} \label{lem:a1cocycle}
Let $X$ be simple of type $A_1$ and $V = (p-2) \otimes 1^{[1]}$. With the notation above, for each $k \in K$ define $\gamma_{k} \, : \, U_{+} \to V$ by
\begin{align*}
\gamma_{k}(x_{+}(t)) &= k\sum_{r = 1}^{p-1}  \binom{p-1}{r}t^{r}e_{2r}.
\end{align*}
Then each $\gamma_{k}$ extends to a rational cocycle $X \to V$, with
\begin{align*}
\gamma_{k}(x_{-}(t)) &= k\sum_{r = 1}^{p-1}  \binom{p-1}{r}t^{r}e_{-2r}
\end{align*}
for all $t \in K$. Furthermore, two such cocycles $\gamma_{k}$ and $\gamma_{l}$ are cohomologous if and only if $k = l$, and every cocycle $X \to V$ is cohomologous to some $\gamma_{k}$.
\end{lemma}

\proof It is clear that $k \mapsto \gamma_{k}$ is a linear map between $K$ and the vector space $Z^{1}(U_{+},V)$ of $1$-cocycles $U_{+} \to V$. Therefore, once we have shown that each $\gamma_{k}$ extends to a cocycle $X \to V$, and that $\gamma_{k}$ and $\gamma_{l}$ are cohomologous if and only if $k = l$, it follows that this map induces an isomorphism $K \to H^{1}(X,V)$, so that every cocycle $X \to V$ is cohomologous to exactly one such $\gamma_{k}$.

Clearly, $\gamma_{k}$ is a rational map on $U_{+}$. To prove that $\gamma_{k}$ extends to a cocycle on $X$, we check firstly that the cocycle condition holds on $U_{+}$, and that $\gamma_{k}$ is a coboundary if and only if $k = 0$, so that $\gamma_{k}$ and $\gamma_{l}$ are cohomologous if and only if $k = l$. We then check that $h(u).\gamma_{k}(x_{+}(t)) = \gamma_{k}({}^{h(u)}x_{+}(t))$ for all $t \in K$ and $u \in K^{*}$, so that $[\gamma_{k}] \in H^{1}(U_{+},V)^{T}$. By \cite[Lemma 1.1 and Theorem 2.1]{MR0439856}, we have vector-space isomorphisms $H^{1}(U_{+},V)^{T} \to H^{1}(U_{+}T,V) \leftarrow H^{1}(X,V)$, and we deduce that each $\gamma_{k}$ extends to a cocycle $X \to V$, again with $[\gamma_{k}] = [\gamma_{l}]$ if and only if $k = l$. Finally, since $x_{-}(-t) = {}^{n(-1)}x_{+}(t)$ and $n(-1)$ sends each vector $e_{2r}$ to $(-1)^{r}e_{-2r}$, it follows immediately that $\gamma_{k}$ has the stated form on $U_{-}$.

So consider $\gamma_{k}(x_{+}(t_{1})) + x_{+}(t_{1}).\gamma_{k}(x_{+}(t_{2}))$. Substituting the expression above for $x_{+}(t_{1}).e_{2s}$, this becomes
\begin{align*}
 &k\left(\left(\sum_{r = 1}^{p-1} \binom{p-1}{r}t_{1}^{r}e_{2r}\right) + \left(\sum_{s = 1}^{p-1}\binom{p-1}{s}\sum_{r = s}^{p - 1}\binom{p-s-1}{r-s}t_{1}^{r-s}t_{2}^{s}e_{2r} \right)\right)\\
=\ &k \left(\left(\sum_{r = 1}^{p-1} \binom{p-1}{r}t_{1}^{r}e_{2r}\right) + \left(\sum_{r = 1}^{p-1}\sum_{s= 1}^{r}\binom{p-1}{s}\binom{p-s-1}{r-s}t_{1}^{r-s}t_{2}^{s}e_{2r} \right)\right)\\
=\ &k\left(\left(\sum_{r = 1}^{p-1} \binom{p-1}{r}t_{1}^{r}e_{2r}\right) + \left(\sum_{r = 1}^{p-1}\sum_{s=1}^{r}\binom{p-1}{r}\binom{r}{s}t_{1}^{r-s}t_{2}^{s}e_{2r} \right)\right)\\
=\ &k\left(\sum_{r = 1}^{p-1} \left(\binom{p-1}{r} (t_{1} + t_{2})^{r} \right) e_{2r}\right)\\
=\ &\gamma_{k}(x_{+}(t_{1} + t_{2}))
\end{align*}
and so $\gamma_{k} \,:\, U_{+} \to V$ is a cocycle. Here we have used the following identity, which holds for all $1 \le s \le r \le p - 1$:
\[ \binom{p-1}{s} \binom{p-s-1}{r-s} = \binom{p-1}{r} \binom{r}{s}. \]

Next, fix $0 \neq k \in K$ and suppose that $\gamma_{k}$ is a coboundary on $U_{+}$, so $\gamma_{k}(x_{+}(t)) = (x_{+}(t).v) - v$ for some $v \in V$. Express $v = \sum_{r = 1}^{p-1}\left(c_{r} e_{2r} + d_{r}e_{-2r}\right)$. It is easily shown that if $d_{r} \neq 0$ for some $r > 1$, then $e_{-2(r-1)}$ has non-zero coefficient in $(x_{+}(t).v) - v$, contradicting the definition of $\gamma_{k}$. So $d_{r} = 0$ for $r > 1$. It then follows that the coefficient of $e_{2}$ in $(x_{+}(t).v) - v$ is identically zero for all $t$ and $k$, while by definition the coefficient of $e_{2}$ in $\gamma_{k}(x_{+}(t))$ is $kt(p-1)$, a contradiction. Thus $\gamma_{k}$ is a coboundary if and only if $k = 0$.

It remains to show that $h(u).\gamma_{k}(x_{+}(t)) = \gamma_{k}({}^{h(u)}x_{+}(t))$. This is straightforward:
\begin{align*}
h(u).\gamma_{k}(x_{+}(t)) &= k\sum_{r = 1}^{p-1}  \binom{p-1}{r}t^{r} h(u).e_{2r}\\
&= k\sum_{r = 1}^{p-1}  \binom{p-1}{r}(u^{2}t)^{r}e_{2r}\\
&= \gamma_{k}(x_{+}(u^{2}t))\\
&= \gamma_{k}({}^{h(u)}x_{+}(t))
\end{align*}
as required. \qed

An entirely similar proof to the above yields the following.
\begin{corollary} \label{cor:a1cocycle}
If $X$ is simple of type $A_{1}$, with root elements $x_{\pm}(t)$ and if $V^{[s]}$ denotes the Frobenius twist of $V = (p-2) \otimes 1^{[1]}$, then $V^{[s]}$ has a basis $\{e_{2-2p},e_{4-2p},\ldots,e_{-2},e_{2},e_{4},\ldots,e_{2p-2}\}$, where $e_{2i}$ has weight $2i(p^s)$ for each $i$, such that every cocycle $X \to V^{[s]}$ is cohomologous to exactly one cocycle $\gamma_{k}$, where
\begin{align*}
\gamma_{k}(x_{+}(t)) &= k\sum_{r = 1}^{p-1}  \binom{p-1}{r}t^{rp^{s}}e_{2r},\\
\gamma_{k}(x_{-}(t)) &= k\sum_{r = 1}^{p-1}  \binom{p-1}{r}t^{rp^{s}}e_{-2r}.
\end{align*}

\end{corollary}

\begin{remark} \label{rem:magma}
We will apply the above results in the context of a simple group $X$ of type $A_{1}$ acting on a group $U$, where $X$ and $U$ are generated by products of root elements in a simple algebraic group $G$. We therefore perform a number of intricate calculations involving products and commutators of root elements of the simple algebraic group $G$. These calculations can in principle be carried out by hand, however we have made use of the computational algebra package \textsc{Magma} \cite{MR1484478} for both speed and accuracy. We have therefore taken our structure constants for computing commutators in $G$ to be consistent with those found in \textsc{Magma}.
\end{remark}

\subsection{From $Q$-conjugacy to $G$-conjugacy}
If $P = QL$ is a parabolic subgroup of $G$ and $X$ is a connected subgroup of the Levi factor $L$, then we now have the necessary tools to study $H^{1}(X,Q)$, but it remains to consider how the corresponding conjugacy classes of subgroups fuse in $G$, and also to consider conjugacy between subgroups in non-conjugate parabolic subgroups of $G$.

Consider first the non-trivial torus $Z(L)$. This centralises $X$ and normalises each root subgroup of $G$. Moreover the action of $Z(L)$ on each such root subgroup in $Q$ is non-trivial since $C_{G}(Z(L)) = L$. This fuses together various classes of complements to $Q$ in $QX$:

\begin{lemma} \label{lem:torus}
Let $G$ be a simple algebraic group over an algebraically closed field $K$ and let $P$ be a parabolic subgroup of $G$ with Levi decomposition $P = QL$. Let $V  = \bigoplus_{i = 1}^{r} V_{i}$ be the sum of the levels of $Q$, let $X$ be a subgroup of $L$, and suppose that $V = \left(\bigoplus_{i = 1}^{m} M_{i} \right) \oplus \left(\bigoplus_{i = 1}^{n} W_{i} \right)$ as $X$-modules, where each $M_{i}$ and each $W_{i}$ is indecomposable, and $H^{1}(X,M_{i}) \cong K$, $H^{1}(X,W_{i}) = 0$.

If for each $i$, the action of $C_{Z(L)}\left(\bigoplus_{j \neq i} M_{j}\right)$ on $M_{i}$ is non-trivial, then complements to $Q$ in $QX$ fall into at most $2^{m}$ classes of subgroups in $G$.
\end{lemma}

\proof Recall that $\mathbb{V} = \mathbb{V}_{X,Q} = \sum_{i = 1}^{r} H^{1}(X,V_{i})$, and fix a basis of $\mathbb{V}$ consisting of a non-zero element from each of the $m$ spaces $H^{1}(X,M_{i})$ ($i = 1,\ldots,m)$. Recall also the surjective map $\rho \, : \, \mathbb{V} \to H^{1}(X,Q)$ from Definition \ref{def:rho}. If $\rho(k_1,k_2,\ldots,k_m) \in H^{1}(X,Q)$, then the hypotheses imply that whenever $k_{i} \neq 0$ for some $i$, we may assume that $k_{i} = 1$ by replacing the corresponding complement with a $Z(L)$-conjugate, without changing the other $k_{j}$. Hence we may take each $k_{j}$ to be either $0$ or $1$. \qed

Next, let $w$ be an element of the Weyl group $W(G)$ and let $\dot{w}$ be a fixed preimage of $w$ in $N_G(T)$. Let $I$ and $J$ be subsets of the simple roots $\Pi$ of $G$, and let $P_{I} = Q_{I}L_{I}$ and $P_{J} = Q_{J}L_{J}$ be standard parabolic subgroups of $G$, where $L_{I}$ and $L_{J}$ are standard Levi subgroups. Let $X$ be a subgroup of $L_{I}$ and suppose that:
\begin{enumerate}[label=\normalfont(\roman*)]
\item \label{list1} $w(I) = J$,
\item \label{list2} For some subgroup $R_{I}$ of $Q_{I}$ generated by root subgroups of $G$, the inclusion $R_{I} \to Q_{I}$ induces a bijection $H^{1}(X,R_{I}) \to H^{1}(X,Q_{I})$,
\item \label{list3} For each root subgroup $U_{\alpha} \le R_{I}$ we have $U_{w(\alpha)} \le Q_{J}$.
\end{enumerate}
Then from \ref{list1} it follows that ${}^{\dot{w}}X$ is a subgroup of $L_{J}$. By \ref{list2}, every complement to $Q_{I}$ in $Q_{I}X$ is $Q_{I}$-conjugate to a subgroup of $R_{I}X$, and by \ref{list3} these are therefore $G$-conjugate to a subgroup of $Q_{J}({}^{\dot{w}}X) \le P_{J}$.

Conjugation by elements of $W(G)$ can also fuse different subgroup classes within a single parabolic. Suppose that $\dot{w} \in N_G(T)$ normalises $X$, and suppose also that $w$ stabilises a certain set of roots, such that the corresponding root subgroups generate an $X$-stable normal subgroup $R$ of $Q$. Then if $\phi$ is a cocycle $X \to R$ and $X_{\phi}$ is the corresponding complement to $R$ in $RX$, then conjugation by $\dot{w}$ sends $X_{\phi}$ to another complement to $R$ in $RX$.

Rephrasing this in terms of $\mathbb{V}$, if elements are represented by $m$-tuples $(k_1,\ldots,k_m)$ with respect to an appropriate basis, then conjugation by $\dot{w}$ induces a permutation on the indices. We will encounter our first instance of this, and the first non-trivial application of Lemma \ref{lem:torus}, in Section \ref{sec:E6D4}.

Two parabolic subgroups of $G$ are called \emph{associated} if their Levi subgroups are $G$-conjugate to one another. The following lemma shows that up to association, there is a unique minimal parabolic subgroup of $G$ containing a given non-$G$-cr subgroup $X$. This prevents double-counting of subgroups during our classification in Sections \ref{sec:e6}--\ref{sec:e8}.

\begin{lemma} \label{lem:associated}
Let $X$ be a closed subgroup of $G$, and let $P_{1}$ and $P_{2}$ be minimal among parabolic subgroups of $G$ containing $X$. Then $P_{1}$ and $P_{2}$ are associated.
\end{lemma}

\proof Let $I$ and $J$ be subsets of the simple roots of $G$ such that the standard parabolic subgroups $P_{I}$ and $P_{J}$ are respectively conjugate to $P_{1}$ and $P_{2}$. It is a standard result \cite[Propositions 2.8.2, 2.8.3]{MR794307} that $P_{1} \cap P_{2}$ is contained in a conjugate of the standard parabolic subgroup $P_{I \cap w(J)}$, for some element $w$ of the Weyl group. From the minimality of $P_{1}$ it follows that $I \cap w(J) = I$, hence $w(J) \supseteq I$. By symmetry, there exists an element $w'$ of the Weyl group such that $w'(I) \supseteq J$. Hence $I$ and $J$ have the same size, so $w(J) = I$, and if $\dot{w} \in N_{G}(T)$ is a lift of $w$ to an element of $N_{G}(T)$, then ${}^{\dot{w}}L_{J} = L_{I}$. \qed


\section{Proof of Theorem \ref{thm:e6}: $G=E_6$, $p=5$} \label{sec:e6}

In this section, we show that when $p = 5$, each non-$G$-cr subgroup of $G$ is conjugate to a subgroup listed in Table \ref{E6p5tab}. Since every subgroup in Table \ref{E6p5tab} is indeed non-$G$-cr by Lemma \ref{lem:BMR}, this proves Theorem \ref{thm:e6}.

Let $P = QL$ be a parabolic subgroup of $G$, such that $P$ contains a non-$G$-cr subgroup $X$ necessarily of type $A_{1}$, and further assume that $P$ is minimal among parabolic subgroups of $G$ containing $X$. Then the image of $X$ in $L$ is an $L'$-irreducible subgroup of $L'$. Moreover, by Corollary \ref{cor:nonzeroH1}, there is some level $M$ of $Q$ such that $H^{1}(X,M \downarrow X) \neq 0$. The following lemma classifies the possible occurrences of this scenario.

\begin{lemma} \label{lem:e6badlevi}
Let $L$ be a Levi subgroup of $G$ containing an $L$-irreducible subgroup $X$ of type $A_1$. If there exists a parabolic subgroup $P$ of $G$ with Levi factor $L$ and unipotent radical $Q$, such that $H^{1}(X,M \downarrow X) \neq 0$ for some level $M$ of $Q$, then $X$ and the type of $L'$ appear in Table \ref{tab:e6p5badX}.
\end{lemma}
\begin{table}[htbp]
\caption{$L'$-irreducible $X$ with $H^1(X,M\downarrow X) \neq 0$\label{tab:e6p5badX}} \vspace{10pt}
\centering
\begin{tabular}{>{\raggedright\arraybackslash}p{0.09\textwidth - 2\tabcolsep}>{\raggedright\arraybackslash}p{0.70\textwidth-\tabcolsep}@{}}
\hline $L'$ & Embedding of $X$ \\ \hline
$D_5$ & $X < D_5$ via $4^{[r]} + 1^{[r+1]} \otimes 1^{[s]} + 0$ $(rs = 0;$ $r+1 \neq s)$ \\
$D_4$ & $X < D_4$ via $4 + 2^{[1]}$, \\
 & $X < D_4$ via $3 \otimes 1^{[1]}$ (two $L'$-classes) \\
$A_1 A_3$ & $X \hookrightarrow A_1 A_3$ via $(1^{[1]}, 3)$ \\
$A_1^{2}A_2$ & $X \hookrightarrow A_1^2 A_2$ via $(1,1^{[1]},2)$ or $(1^{[1]},1,2)$ \\
$D_5$ & $X < D_5$ via $2 \otimes 2^{[1]} + 0$ \\ \hline
\end{tabular}
\end{table}

\proof Let $P = QL$ be a parabolic subgroup of $G$. The action of $L'$ on the levels of $Q$ is straightforward to determine, as described in \cite{MR1047327}. Now Proposition \ref{prop:irredA1} gives all $L'$-irreducible subgroups of type $A_1$, and it is straightforward to determine the action of each such subgroup $X$ on the levels of $Q$, for instance using the tables of Section \ref{tabs:misc}. Checking each level $M$ against Lemma \ref{lem:h1fora1} tells us whether or not $H^{1}(X,M \downarrow X) = 0$. Whenever we find a module $M$ such that $H^{1}(X,M \downarrow X) \neq 0$, the full description of the action of $X$ on the levels of $Q$ is given in the relevant section below. So let us illustrate the process with an example where $H^{1}(X,M \downarrow X) = 0$ for all levels $M$ of $Q$.

Let $L$ be the unique standard Levi subgroup of $G$ such that $L' = A_{5}$. Since $p = 5$, the only $L'$-irreducible subgroup $X$ of type $A_{1}$ acts on the natural module as $2^{[r]} \otimes 1^{[s]}$ where $rs = 0$ and $r \neq s$. In $Q$, there are two levels. As an $L'$-module, $Q/Q(2)$ is irreducible with high weight $\lambda_3$, and is generated by the image of the root subgroup $U_{\alpha_2}$, and $Q(2)$ is irreducible with high weight $0$, and is generated by $U_{122321}$. Using Table \ref{tab:extpowers} and the fact that $V_{A_5}(\lambda_3) = \bigwedge^{3} (V_{A_5}(\lambda_1)$), it follows that $V_{A_5}(\lambda_3) \downarrow X = (4^{[r]} \otimes 1^{[s]}) + (2^{[r]} \otimes 1^{[s]}) + 3^{[s]}$. Lemma \ref{lem:h1fora1} shows that no indecomposable summand of such a module has a nonvanishing first cohomology group. \qed

The remainder of this section is as follows. There is a subsection for each Lie type of $L'$ in Table~\ref{tab:e6p5badX}. For each type, we enumerate the standard parabolics with Levi factor of that type and let $X$ vary over the $L'$-irreducible subgroups in Table \ref{tab:e6p5badX}. Using the tools from Section \ref{sec:preliminaries}, we bound the number of possible classes of complements to $Q$ in $QX$, as $Q$ varies over unipotent radicals of these parabolics, and we thus determine a global bound on the number of non-$G$-cr subgroups arising from such a parabolic. Finally, we exhibit representatives of each possible conjugacy class and calculate their connected centralisers. Lemma \ref{lem:associated} shows that two such subgroups lying in non-associated parabolics cannot be $G$-conjugate. Table \ref{E6p5tab} contains precisely these representatives.

\subsection{$L' = D_5$} \label{sec:E6D5}
The two standard $D_5$-parabolic subgroups of $G$ are $P_{12345}$ and $P_{23456}$. By Lemma \ref{lem:e6badlevi}, the only embeddings of $A_{1}$ into $D_5$ that we need to consider are $X = A_1 < D_5$ via $4^{[r]} + 1^{[r+1]} \otimes 1^{[s]} + 0$, where $rs = 0$ and $r+1 \neq s$.

Consider first $P = P_{12345} = QL$. Then $Q$ has a single level, and is a $16$-dimensional irreducible module for $L'$ of high weight $\lambda_4$. From the action of $X$ on the natural $D_{5}$-module, we see that $X$ lies in a subgroup $A_{1}^{2}B_{2}$ of $L'$, via $(1^{[s]},1^{[r+1]},4^{[r]})$. In turn, this subgroup lies in a subgroup $A_{1}^{2}A_{3}$ of $L'$. Now, the restriction of a spin module for $D_{n}$ $(n \ge 4)$ to a subgroup $B_{n-1}$, or from $B_{n}$ or $D_{n}$ to a proper Levi subgroup of the same type, is itself a sum of spin modules \cite[Prop.\ 2.7]{MR1329942}. It follows that $V_{D_{5}}(\lambda_4) \downarrow A_{1}^{2}A_{3} = (1,0,\lambda_1) + (0,1,\lambda_3)$, hence $V_{D_5}(\lambda_4) \downarrow A_{1}^{2}B_{2} = V_{D_5}(\lambda_5) \downarrow A_{1}^{2}B_{2} = (1,0,\lambda_2) + (0,1,\lambda_2)$, and so $V_{D_5}(\lambda_4) \downarrow X = V_{D_5}(\lambda_{5}) \downarrow X = (3^{[r]} \otimes 1^{[r+1]}) + (3^{[r]} \otimes 1^{[s]})$. By Lemma \ref{lem:h1fora1}, $\mathbb{V} \cong K$, since $r+1 \neq s$. The 1-dimensional torus $Z(L)$ acts non-trivially on each summand of $Q$, hence by Lemma \ref{lem:torus} there is at most one $G$-conjugacy class of non-$G$-cr complements to $Q$ in $QX$. An entirely similar argument applies to the parabolic subgroup $P_{23456}$, yielding at most one $G$-conjugacy class of non-$G$-cr complements to $Q_{23456}$ in $Q_{23456}X$.

Let $Y$ and $Z$ be subgroups of the subsystem subgroup $A_1 A_5$ of $G$, respectively embedding via $(1^{[s]},W(5)^{[r]})$ and $(1^{[s]},(W(5)^{[r]})^{*})$ $(rs=0;$ $r+1 \neq s)$. Note that $Y$ and $Z$ are $\text{Aut}(G)$-conjugate, since a graph automorphism of $E_6$ induces a graph automorphism of $A_1 A_5$. Both $Y$ and $Z$ are non-$A_1 A_5$-cr, hence by Lemma \ref{lem:BMR} they are non-$G$-cr. We claim that $Y$ and $Z$ each lie in a $D_5$-parabolic subgroup of $G$, with irreducible image in a Levi factor. Their actions on the natural $A_1 A_5$-modules imply that $Y$ and $Z$ are each contained in an $A_1^2 A_3$-parabolic subgroup of $A_1 A_5$, and it follows that $Y$ and $Z$ each lie in a parabolic subgroup of $G$ whose Levi factor contains a subgroup of type $A_1^{2} A_3$. The only such subgroups are $D_5$-parabolic subgroups, hence $Y$ and $Z$ each lie in a $D_5$-parabolic.

The action of $Y$ on $V_{27}$ is given in Table \ref{E6p5tab}, and this determines the action of $Z$ since the outer automorphism of $G$ swaps the $G$-modules $V_{27}$ and $V_{27}^{\ast}$, and swaps the subgroup classes of $Y$ and $Z$. Now, if $Y$ and $Z$ were $G$-conjugate, then the submodule latties of $V_{27} \downarrow Y$ and $V_{27} \downarrow Z$ would be identical. But this is not the case, since $Y$ has an $8$-dimensional submodule $1^{[r]} \otimes 3^{[s]} \subset 1^{[r]} \otimes W(5)^{[s]}$, while $Z$ does not. Finally, the $G$-conjugacy class of a reductive sugroup determines its high weights on each $G$-module, up to multiplying by a power of $p$. Therefore if $r'$, $s'$ are non-negative integers with $r's' = 0$ and $(r',s') \neq (r,s)$, then the corresponding non-$G$-cr subgroups arising are not conjugate to $Y$ or to $Z$. Thus for each $(r,s)$, there are precisely two classes of non-$G$-cr subgroups $A_{1}$ contained in a $D_{5}$-parabolic with irreducible image in a Levi factor, with representatives $Y$ and $Z$.

Inspecting Table \ref{E6p5tab}, we see that $C_{L(G)}(Y) = C_{L(G)}(Z) = \{0\}$, and thus $C_{G}(Y)^{\circ} = C_{G}(Z)^{\circ} = 1$.

\begin{remark}
The condition $r + 1 \neq s$ is only necessary to ensure that $X$ is $D_5$-irreducible. If instead $(r,s)=(0,1)$, the subgroups $Y$ and $Z$ described above are still non-$G$-cr and non-conjugate. The images of these subgroups in $D_5$ now lie in a Levi subgroup $D_4$, and hence these provide two non-conjugate, non-$G$-cr subgroups lying in a $D_4$-parabolic subgroup of $G$.
\end{remark}

\subsection{$L' = D_4$} \label{sec:E6D4}

Let $P = P_{2345} = QL$ be the unique standard $D_4$-parabolic subgroup of $G$. Let $X$, $Y$, $Z$ be representatives of the three $L'$-conjugacy classes of $L'$-irreducible subgroups, with $V_{D_4}(\lambda_1) \downarrow X \cong V_{D_4}(\lambda_3) \downarrow Y = V_{D_4}(\lambda_4) \downarrow Z \cong 4 + 2^{[1]}$, so that the remaining 8-dimensional modules $V_{D_4}(\lambda_i)$ restrict as $3 \otimes 1^{[1]}$. Now, $Q$ has two levels: $Q/Q(2) \cong V_{D_4}(\lambda_3) + V_{D_4}(\lambda_4)$, the factors respectively generated as a $D_4$-module by the images of the root groups $U_{\alpha_1}$ and $U_{\alpha_6}$, and $Q(2) \cong V_{D_4}(\lambda_1)$, generated by $U_{101111}$.

The action of $X$ on $Q$ is as follows:
\begin{align*}
Q/Q(2) \downarrow X &= (3 \otimes 1^{[1]}) + (3 \otimes 1^{[1]}),\\
Q(2) \downarrow X &= 4 + 2^{[1]}.
\end{align*}
And thus $\mathbb{V}_{X,Q} \cong K^{2}$ by Lemma \ref{lem:h1fora1}. In the second level, $H^{2}(X,2^{[1]}) \cong K$ by Lemma \ref{lem:h2fora1}. This means that not all pairs $(k,l) \in \mathbb{V}_{X,Q}$ necessarily give rise to an element of $H^{1}(X,Q)$. We will show that the condition $kl = 0$ is necessary.

For this, we now describe $X$ explicitly in terms of the root groups of $G$, which allows us to identify the weight vectors of $X$ in its action on each level of $Q$, and hence apply Lemma \ref{lem:a1cocycle} to give an explicit description of cocycles $X \to Q$. Using the module decomposition $V_{L'}(\lambda_1) \downarrow X = 4 + 2^{[1]}$, we can identify root elements $x_{\pm}(t)$ of $X$:
\begin{align*}
x_{+}(t) &= x_{\alpha_3}(3t)x_{001100}(2t^2)x_{001110}(t^3)x_{\alpha_4}(2t)x_{000110}(4t^2) x_{\alpha_5}(t) x_{010110}(2t^5) x_{011100}(t^5),\\
x_{-}(t) &= x_{-\alpha_3}(t)x_{-001100}(t^2)x_{-001110}(t^3)x_{-\alpha_4}(2t)x_{-000110}(3t^2) x_{-\alpha_5}(3t) x_{-010110}(3t^5) x_{-011100}(t^5).
\end{align*}
Furthermore, a maximal torus of $X$ is given by $T_{X} = \{h(t) \, : \, t \in K^{\ast}\}$ with $h(t)$ as defined in Lemma \ref{lem:a1rels}. Multiplying out the above elements gives the following formula for $h(t)$:
\begin{align*}
h(t) &= h_{2}(t^{10})h_{3}(t^{8})h_{4}(t^{14})h_{5}(t^{8}).
\end{align*}
We can now check directly that each non-trivial element of the form $x_{111100}(t)x_{010111}(u)$ has weight $2$ under the action of $T_{X}$. Next, when $p = 5$ the formula for $x_{+}(t).e_{2}$ on page \pageref{eqn:xplus} becomes
\begin{align*}
x_{+}(t).e_{2} = e_{2} + 3te_{4} + 3t^{2}e_{6} + t^{3}e_{8}.
\end{align*}
In $Q/Q(2)$ we can therefore let $e_{2} = x_{111100}(1)Q(2)$ (resp.\ $x_{010111}(1)Q(2)$) and calculate the conjugate ${}^{x_{+}(t)}e_{2}$, and then equate coefficients with the equation above to find that $e_{4} = x_{111110}(3)$, $e_{6} = x_{111210}(3)$, $e_{8} = x_{112210}(4)$ (resp.\ $x_{011111}(1)$, $x_{011211}(1)$, $x_{011221}(1)$). Hence by Lemma \ref{lem:a1cocycle}, a general cocycle $X \to Q/Q(2)$ is cohomologous to exactly one map $\gamma_{k,l}$ as follows:
\[ \gamma_{k,l} \ : \ \left\{
\begin{array}{ll}x_{+}(c) \mapsto & x_{111100}(kc) x_{11110}(2kc^2) x_{111210}(3kc^3) x_{112210}(kc^4) \\
	& \times \ x_{010111}(lc) x_{011111}(4lc^2) x_{011211}(lc^3) x_{011221}(4lc^4) Q(2).
\end{array}
\right. \]
Now if $\phi \in Z^{1}(X,Q)$ has image $\gamma_{k,l}$ under composition with $Q \twoheadrightarrow Q/Q(2)$, then it follows that for each $x \in X$ we have $\phi(x) = \gamma_{k,l}(x)q$ for some $q \in Q(2) = Z(Q)$ depending on $x$. From the identity $xq = xqx^{-1}x = ({}^{x}q)x$, for all $q \in Z(Q)$, $x \in X$, we then have
\begin{align*}
(\phi(x_{+}(c))x_{+}(c))^{5} &= (\gamma_{k,l}(x_{+}(c))qx_{+}(c))^{5} \\
&= q ({}^{x_{+}(c)}q)({}^{x_{+}(c)^{2}}q)({}^{x_{+}(c)^{3}}q)({}^{x_{+}(c)^{4}}q)( \gamma_{k,l}(x_{+}(c))x_{+}(c))^{5}.
\end{align*}
Since $x_{+}(c)$ induces a linear transformation of order $5$ on the vector space $Q(2)$, call it $T$, in additive notation the first five terms in the above product are $(1 + T + T^{2} + T^{3} + T^{4})(q)$, which is identically zero. So the above is equal to $(\gamma_{k,l}(x_{+}(c))x_{+}(c))^{5}$. Substituting, and using \textsc{Magma} to simplify calculations, we find that this is equal to $x_{122321}(-klc^5)$. But since $\phi$ is a cocycle, the element $\phi(x_{+}(c))x_{+}(c)$ is a positive root element in a complement to $Q$ in $QX$, and in particular its order divides $5$. Thus $x_{122321}(-klc^5) = 1$ for all $c \in K$, hence $kl = 0$ as claimed.

The $2$-dimensional torus $Z(L)$ consists of elements of the form
\[ h(t,u) \stackrel{\textup{def}}{=} h_{1}(t^{-2}u^{2}) h_{2}(t) h_{3}(u) h_{4}(t^{2}) h_{5}(t^{3}u^{-1})h_{6}(t^{4}u^{-2}) \]
for $t$, $u \in K^{\ast}$. Now ${}^{h(t,u)}x_{\alpha_1}(c) = x_{\alpha_1}(t^{4} u^{-3} c)$ and ${}^{h(t,u)}x_{\alpha_6}(c) = x_{\alpha_6}(t^{-5}u^{3}c)$. Thus the action of $Z(L)$ on $Q$ satisfies Lemma \ref{lem:torus}, and each complement to $X$ in $QX$ is $G$-conjugate to one of the complements $X = X_{[0,0]}$, $X_{[1,0]}$ or $X_{[0,1]}$ (recall the notation from Section \ref{sec:nonab}). So in particular, there exist at most two $G$-conjugacy classes of non-$G$-cr complements to $Q$ in $QX$.

Next consider the action of $Y$ on $Q$:
\begin{align*}
Q/Q(2) \downarrow Y &= 4 + 2^{[1]} + (3 \otimes 1^{[1]}),\\
Q(2) \downarrow Y &= 3 \otimes 1^{[1]}
\end{align*}
and thus $\mathbb{V}_{Y,Q} \cong K^{2}$ by Lemma \ref{lem:h1fora1}. Again, by Lemma \ref{lem:torus} each complement to $Q$ in $QY$ is $G$-conjugate to one of $Y = Y_{[0,0]}$, $Y_{[1,0]}$, $Y_{[0,1]}$ or $Y_{[1,1]}$. By an identical argument, each complement to $Q$ in $QZ$ is $G$-conjugate to one of $Z = Z_{[0,0]}$, $Z_{[1,0]}$, $Z_{[0,1]}$ or $Z_{[1,1]}$. We now exhibit $G$-conjugacies between these various subgroup classes, and then construct a representative of each possible distinct class that remains.

We claim that $Y_{[1,0]}$, $Y_{[0,1]}$ and $X_{[0,1]}$ are $G$-conjugate to one another, and that $Z_{[1,0]}$, $Z_{[0,1]}$ and $X_{[1,0]}$ are $G$-conjugate to one another. Indeed, recalling notation from Section \ref{sec:notation}, the element $w_1 = n_1 n_3 n_4 n_2 n_5 n_4 n_3 n_1 \in N_{G}(T)$ acts as an outer involution in $D_4.S_3$, preserving $Y$ whilst swapping $X$ and $Z$. Furthermore, $w_1$ swaps the root subgroups $U_{\alpha_6}$ and $U_{101111}$, hence swaps the two $Y$-modules generated by these elements, and also swaps the $X$-invariant and $Z$-invariant subgroups generated by these. It follows that $Y_{[1,0]}$ is $G$-conjugate to $Y_{[0,1]}$, and $X_{[1,0]}$ is $G$-conjugate to $Z_{[0,1]}$. Similarly, the element $n_6 n_5 n_4 n_2 n_3 n_4 n_5 n_6$ stabilises $Z$ and swaps $X$ and $Y$, and also swaps the root subgroups $U_{\alpha_1}$ and $U_{101111}$, and therefore $Z_{[1,0]}$ is $G$-conjugate to $Z_{[0,1]}$ and $X_{[0,1]}$ is $G$-conjugate to $Y_{[0,1]}$. This proves the claim, and we now have at most four $G$-conjugacy classes of non-$G$-cr subgroups $A_1$ of $P$, namely $X_{[1,0]}$, $X_{[0,1]}$, $Y_{[1,1]}$ and $Z_{[1,1]}$.

Consider the following four $A_1$ subgroups of $G$: Two subgroups $A_1 \hookrightarrow A_1 A_5$ via $(1^{[1]},W(5))$ and via $(1^{[1]},W(5)^*)$, respectively; and two classes of subgroups $A_1 < D_5$ via $T(8)$ (these exist by Lemma \ref{lem:twoclasses}). Each of these four subgroups is non-$G$-cr by Lemma \ref{lem:BMR}. The two subgroups of $A_{1}A_{5}$ are interchanged by a graph automorphism of $A_{1}A_{5}$, hence by a graph automorphism of $E_{6}$, and similarly the two subgroups of $D_{5}$ are exchanged by a graph automorphism of $E_{6}$. The third and fourth lines of Table \ref{E6p5tab} give the actions of one subgroup of $A_{1}A_{5}$ and one subgroup of $D_{5}$ on the modules $V_{27}$ and $L(G)$, which also determines the actions of the other two subgroups. In particular, we see that no two of these four subgroups have the same lattice of submodules on $V_{27}$, and hence all four are pairwise non-conjugate in $G$.

The subgroups in $A_1 A_5$ each lie in an $A_1^2 A_3$-parabolic subgroup of $A_1 A_5$, which is contained in a $D_5$-parabolic subgroup of $G$. The images of these subgroups in the Levi factor $D_{5}$ then act on $V_{D_5}(\lambda_1)$ as $4 + 1^{[1]} \otimes 1^{[1]} + 0 = 4 + 2^{[1]} + 0^2$, and so these images are $D_4$-irreducible subgroups of a Levi subgroup $D_4$.

Similarly, from the module structure of $T(8) = 0|(3 \otimes 1^{[1]})|0$ it follows that the two subgroups in $D_5$ each lie in a $D_4$-parabolic subgroup of $D_5$, and hence of $G$. Thus all four of these subgroups lie in a $D_{4}$-parabolic subgroup of $G$.

The connected centraliser of the non-$G$-cr subgroups contained in $A_1 A_5$ is trivial, as their fixed-point space on $L(G)$ is zero. The connected centraliser of each subgroup $X = A_1 < D_5$ via $T(8)$ is $T_1$. Indeed $\text{dim}(C_G(X)^\circ) \le 1$ by Table \ref{E6p5tab}, and $T_1 = C_G(D_5)^\circ \leq C_G(X)^\circ$.

\subsection{$L' = A_1 A_3$}

The four standard $A_{1} A_{3}$-parabolic subgroups of $G$ are $P_{1346}$, $P_{2346}$, $P_{1456}$ and $P_{1245}$. For each, let $X = A_1 \hookrightarrow A_1 A_3$ via $(1^{[1]},3)$. 

First, consider $P_{1346}$. Then $Q = Q_{1346}$ has four levels and $X$ acts as follows.
\begin{align*}
Q/Q(2) \downarrow X &= (3 \otimes 1^{[1]}) + 3,\\
Q(2)/Q(3) \downarrow X &= (4 \otimes 1^{[1]}) + 1^{[1]}, \\
Q(3)/Q(4) \downarrow X &= 3, \\
Q(4) \downarrow X &= 0,
\end{align*}
where the module of highest weight $3 \otimes 1^{[1]}$ is generated by $U_{\alpha_5} Q(2)$. By Lemma \ref{lem:h1fora1}, $\mathbb{V}_{X, Q_{1346}} \cong K$. Applying Lemma \ref{lem:torus}, we have at most one class of non-$G$-cr complements to $Q_{1346}$ in $Q_{1346}X$. The same argument applied to $P_{2346}$ allows us to conclude that there is at most one class of non-$G$-cr complements to $Q_{2346}$ in $Q_{2346}X$, with the module of highest weight $3 \otimes 1^{[1]}$ in the action of $X$ on $Q_{2346}$ generated by $U_{101110} Q_{2346}(3)$. Conjugation by the element $n_1 n_3 n_4 n_2 \in N_{G}(T)$ sends $L_{1346}$ to $L_{2346}$, and also sends $U_{\alpha_5}$ to $U_{101110}$. It follows that any non-$G$-cr $A_{1}$ subgroup of $P_{1346}$ with irreducible image in $L_{1346}$ is $G$-conjugate to a subgroup of $P_{2346}$. 

Similarly, considering the parabolic subgroups $P_{1456}$ and $P_{1245}$ leads to at most one $G$-conjugacy class of non-$G$-cr subgroups, where now conjugation by the element $n_6 n_5 n_4 n_2 \in N_{G}(T)$ sends $L_{1456}$ to $L_{1245}$ and sends $U_{\alpha_3}$ to $U_{001111}$. 

Consider the subgroups $Y$, $Z = A_1 < A_5$ via $W(5)$ and $W(5)^{*}$, respectively. Each of these is non-$G$-cr by Lemma \ref{lem:BMR}, and they are contained in an $A_1 A_3$-parabolic subgroup of $G$ since they are contained in an $A_1 A_3$-parabolic subgroup of $A_5$. They are exchanged by an outer automorphism of $A_{5}$, and are therefore $\textup{Aut}(G)$-conjugate. Thus the first line of Table \ref{E6p5tab} determines the action of each on $V_{27}$ and $L(G)$, and in particular we see that the subgroups have different submodule lattices on $V_{27}$, and are therefore not $G$-conjugate. Hence these are representatives of the above conjugacy classes of non-$G$-cr subgroups in $A_1 A_3$-parabolic subgroups of $G$. 

From Table \ref{E6p5tab} we see that $C_{G}(Y)$ and $C_{G}(Z)$ are each at most $3$-dimensional, and so $C_{G}(Y)^{\circ} = C_{G}(Z)^{\circ} = C_G(A_5)^\circ = A_{1}$.

\subsection{$L' = A_1^{2} A_2$}

The five standard $A_{1}^{2} A_{2}$-parabolic subgroups of $G$ are $P_{1235}$, $P_{1236}$, $P_{1256}$, $P_{2356}$ and $P_{1246}$. For each, we need to consider subgroups $X$ and $Y$, respectively embedded in $A_1^2 A_2$ via $(1^{[1]},1,2)$ and $(1,1^{[1]},2)$.

First, consider $P_{1235}$. Then $Q = Q_{1235}$ has four levels, and
\begin{align*}
Q/Q(2) \downarrow L_{1235} &= (1,1,01) + (0,1,00),\\
Q(2)/Q(3) \downarrow L_{1235} &= (1,0,01) + (0,0,10), \\
Q(3)/Q(4) \downarrow L_{1235} &= (0,1,10), \\
Q(4) \downarrow L_{1235} &= (1,0,00). 
\end{align*}
Applying Lemma \ref{lem:h1fora1}, it follows that $\mathbb{V}_{X,Q} \cong \mathbb{V}_{Y,Q} \cong K$. Applying Lemma \ref{lem:torus} yields at most one class of non-$G$-cr complements to $Q$ in $QX$ and at most one class of non-$G$-cr complements to $Q$ in $QY$.

With entirely similar calculations, we find that there are at most two $G$-conjugacy classes of non-$G$-cr simple subgroups of type $A_{1}$ in each of the other four standard parabolic subgroups having irreducible image in the corresponding Levi factor. Moreover, for each such parabolic $P$, we find an element of $N_{G}(T)$, given in the table below, whose image in the Weyl group sends the roots of the standard Levi factor to those of $L_{1235}$ and sends $\alpha$ to $\alpha_{4}$, where the image of $U_{\alpha}$ generates the $X$-module of high weight $3 \otimes 1^{[1]}$ in the appropriate level. This element therefore conjugates each non-$G$-cr subgroups $A_1$ of $P$, with irreducible image in $L$, to a subgroup of $P_{1235}$.
\begin{center}
\begin{tabular}{c|c|c}
$P$ & Root $\alpha$ &  Element of $N_{G}(T)$ \\ \hline 
$P_{1236}$ & $000110$ & $n_6 n_5$ \\
$P_{1256}$ & $001100$ & $n_{6} n_{5} n_{4} n_{2} n_{3} n_{1} n_{4} n_{3} n_{5} n_{4} n_{2} n_{6} n_{5} n_{4}$ \\
$P_{2356}$ & $\alpha_4$ & $n_{6} n_{5} n_{4} n_{2} n_{3} n_{1} n_{4} n_{3} n_{5} n_{4} n_{2} n_{6} n_{5} n_{4} n_{3} n_{1}$ \\
$P_{1246}$ & $001110$ & $n_{6} n_{5} n_{4} n_{2} n_{3} n_{1} n_{4} n_{3} $ \\
\end{tabular}
\end{center}
Note that the element above sending $L_{2356}$ to $L_{1235}$ sends $U_{\alpha_2}$ to $U_{\alpha_5}$ and $U_{\alpha_3}$ to $U_{\alpha_2}$, hence sends the subgroup embedded in $L_{2356}$ via $(1^{[1]},1,2)$ to the subgroup of $L_{1235}$ embedded via $(1,1^{[1]},2)$, and vice versa.

Therefore there are at most two $G$-conjugacy classes of non-$G$-cr subgroups arising in this case. By Lemma \ref{lem:twoclasses}, there are two conjugacy classes of subgroups embedded in $D_5$ via $T(6)$. Let $Z_1$ and $Z_2$ be representatives of these. As $D_{5}$ is a Levi subgroup of $G$, both $Z_{1}$ and $Z_{2}$ are non-$G$-cr by Lemma \ref{lem:BMR}. They are each contained in an $A_1^2 A_2$-parabolic subgroup of $D_5$, hence in such a parabolic of $G$, and they are exchanged by an outer automorphism of $D_{5}$, hence are $\textup{Aut}(G)$-conjugate. Therefore the second line of Table \ref{E6p5tab} determines their actions on $V_{27}$, and we see that they are not $G$-conjugate.

Since $\textup{dim}(C_G(Z_1))$ and $\textup{dim}(C_G(Z_2)) \le 1$, we find that $C_G(Z_1)^{\circ} = C_G(Z_2)^{\circ} = C_G(D_5)^{\circ} = T_1$.

\subsection{$L' = D_5$ (Previously missing case)} \label{sec:E6D5again}

Recall that the two standard $D_5$-parabolic subgroups of $G$ are $P_{12345}$ and $P_{23456}$, and let $X = A_1 < D_5$ via $2 \otimes 2^{[1]} + 0$. In each case the unipotent radical is abelian, and as given in Table~\ref{tab:spinD5}, this is a spin module $\lambda_4$ on which $X$ acts as $3 \otimes 1^{[1]} + 3^{[1]} \otimes 1$. This has $1$-dimensional cohomology, giving rise to one non-$G$-cr subgroup of type $A_1$ in each case, after considering the action of $Z(L)^{\circ}$. Call these non-$G$-cr subgroups $Y$ and $Z$. The actions of $Y$ and $Z$ are given in Table~\ref{E6p5tab}, and since the actions on $V_{27}$ are inequivalent, we see that $Y$ and $Z$ are not $G$-conjugate. For both subgroups, $V_{27}$ has indecomposable summand dimensions $10$ and $17$; and $L(G)$ has summand dimensions $40$, $20$, $15$ and $3$. This is incompatible with the summand dimensions of all the reductive maximal connected subgroups of $G$; hence $Y$ and $Z$ are maximal among reductive subgroups of $G$. 

Inspecting Table~\ref{E6p5tab}, we see that $C_{L(G)}(Y) = C_{L(G)}(Z) = \{0\}$, and thus $C_{G}(Y)^{\circ} = C_{G}(Z)^{\circ} = 1$.


\section{Proof of Theorem \ref{thm:e7}: $G = E_7$, $p = 5$} \label{sec:e7p5}

In this section we prove Theorem \ref{thm:e7} in the case $p=5$. Again, our starting point is to determine those parabolics $P = QL$ of $G$ and $L$-irreducible subgroups $X$ of $L$, such that $H^{1}(X,Q)$ may be non-zero.

\begin{lemma} \label{lem:e7p5badlevi}
Let $L$ be a Levi subgroup of $G$ containing an $L$-irreducible subgroup $X$ of type $A_1$. If there exists a parabolic subgroup $P$ of $G$ with Levi factor $L$ and unipotent radical $Q$, such that $H^{1}(X,M \downarrow X) \neq 0$ for some level $M$ of $Q$, then $X$ and the type of $L'$ appear in Table \ref{tab:e7p5badX}.
\end{lemma}
\begin{longtable}{>{\raggedright\arraybackslash}p{0.14\textwidth - 2\tabcolsep}>{\raggedright\arraybackslash}p{0.86\textwidth-\tabcolsep}@{}}
\caption{$L'$-irreducible $X$ with $H^1(X,M\downarrow X) \neq 0$\label{tab:e7p5badX}} \\
$L'$ & Embedding of $X$ \\ \hline
$E_6$ & $X \hookrightarrow A_1 A_5 < E_6$ via $(1,2 \otimes 1^{[1]})$ \\
$A_1 D_5$ & $X \hookrightarrow A_1 D_5$ via $(1^{[r]},2^{[r]} + 2^{[r+1]} + 1^{[r]} \otimes 1^{[s]})$ $(rs=0$; $r \neq s)$ \\
$A_1 A_2 A_3$ & $X \hookrightarrow A_1 A_2 A_3$ via $(1^{[r+1]},2^{[s]},3^{[r]})$ $(rs=0)$ \\
$D_5$ & $X < D_5$ via $4^{[r]} + 1^{[r+1]} \otimes 1^{[s]} + 0$ $(rs = 0;$ $r+1 \neq s)$ \\
$A_1 D_4$ & $X \hookrightarrow A_1 D_4$ via $(1^{[r]},3^{[s]} \otimes 1^{[s+1]})$ $(rs=0)$ (two $L'$-classes) \\
$A_2 A_3$ & $X \hookrightarrow A_2 A_3$ via $(2,1 \otimes 1^{[1]})$ \\
$A_1^2 A_3$ & $X \hookrightarrow A_1^2 A_3$ via $(1^{[r]},1^{[s+1]},3^{[s]})$ or $(1^{[s+1]},1^{[r]},3^{[s]})$ ($rs=0$), \\
& $X \hookrightarrow A_1^2 A_3$ via $(1,1,1 \otimes 1^{[1]})$ \\
$A_1^3 A_2$ & $X \hookrightarrow A_1^3 A_2$ via $(1^{[r]},1^{[s]},1^{[t]},2^{[u]})$ \newline ($rstu=0$; $u=i=j-1$ where $i,j \in \{r,s,t\}$) \\
$D_4$ & $X < D_4$ via $4 + 2^{[1]}$, \\
 & $X < D_4$ via $3 \otimes 1^{[1]}$ (two $L'$-classes) \\
$A_1 A_3$ & $X \hookrightarrow A_1 A_3$ via $(1^{[1]},3)$ \\
$A_1^2 A_2$ & $X \hookrightarrow A_1^2 A_2$ via $(1,1^{[1]},2)$ or $(1^{[1]},1,2)$ \\ 
$D_5$ & $X < D_5$ via $2 \otimes 2^{[1]} + 0$ \\ \hline
\end{longtable}

\proof As in the proof of Lemma \ref{lem:e6badlevi}, this is straightforward. For each parabolic $P = QL$ and each $L$-irreducible subgroup $X$ of type $A_1$, we systematically check the action of $X$ on each level of $Q$ using Proposition \ref{prop:irredA1}. \qed

\subsection{$L'=E_6$} \label{sec:E7E6p5}

Let $P = P_{123456} = QL$ be the unique standard $E_6$-parabolic subgroup of $G$. Let $X \cong A_1 \hookrightarrow  A_1 A_5 < E_6$ via $(1,2 \otimes 1^{[1]})$. Then $Q$ is abelian, and carries the structure of a rational $KL$-module, with $Q \downarrow X = 3 \otimes 1^{[1]} + 1 \otimes 1^{[1]} + 2 \otimes 2^{[1]} + 4 + 0$. By Lemma \ref{lem:h1fora1}, $\mathbb{V} = H^1(X,Q) \cong K$. By Lemma \ref{lem:torus}, there is at most one class of non-$G$-cr complements to $Q$ in $QX$.

Let $Y \cong A_1 < A_7$ via $W(7)$. Then $Y$ is non-$A_7$-cr, lying in an $A_1 A_5$-parabolic subgroup of $A_7$. By Lemma \ref{lem:BMR}, $Y$ is non-$G$-cr. The projection of $Y$ to $A_1 A_5$ embeds via $(1, 2 \otimes 1^{[1]})$. The only subgroups $A_{1} A_{5}$ of $G$ which lie in a subgroup $A_{7}$ are those in the conjugacy class of $A_1 A_5' < E_6$. Moreover, the only Levi subgroup of $G$ that contains $A_1 A_5'$ is $E_6$. Therefore, the $A_1 A_5$-parabolic subgroup of $A_{7}$ must lie in an $E_6$-parabolic subgroup of $G$. Hence $Y$ lies in an $E_6$-parabolic subgroup of $G$ and the projection of $Y$ to $E_6$ is $E_6$-irreducible by Proposition \ref{prop:irredA1}. Hence $Y$ is a representative of the $G$-conjugacy class above.  

From Table \ref{E7p5tab}, we see that $\text{dim}(C_{G}(Y)^\circ) \leq 1$. Since $Q$ is abelian and $X$ centralises a 1-dimensional subgroup $U_{1}$ of $Q$, we deduce that $U_{1} \le C_G(QX)^{\circ} \le C_G(Y)^{\circ}$, and so $C_G(Y)^\circ = U_1$.

\subsection{$L'=A_1 D_5$}  \label{sec:E7A1D5p5}

Let $P = P_{123457} = QL$ be the unique standard $A_{1} D_{5}$-parabolic subgroup of $G$. Let $X \cong A_1 \hookrightarrow A_1 D_5$ via $(1^{[r]}, 1^{[r]} \otimes 1^{[s]} + 2^{[r]} + 2^{[r+1]} )$ $(rs=0, r \neq s)$. Then the projection of $X$ to $D_5$ is contained in a subgroup $\bar{A}_1^2 A_1^2 \leq \bar{A}_1^2 A_3$. The two $A_{1}$ factors of this latter group are generated by the root subgroups $U_{\pm \alpha_{2}}$ and $U_{\pm \alpha_{5}}$, respectively, and the $A_{3}$ factor is generated by $U_{\pm \alpha_{1}}$, $U_{\pm \alpha_{3}}$ and $U_{\pm 1122100}$. This allows us to pick the following root elements of $X$:
\begin{align*}
x_{+}(t) &= x_{\alpha_{7}}(t^{5^r})x_{\alpha_2}(t^{5^r})x_{\alpha_{5}}(t^{5^s})x_{-1122100}(t^{5^r})x_{\alpha_{1}}(t^{5^r}) x_{-1112100}(t^{5^{r+1}})x_{1010000}(t^{5^{r+1}}),\\
x_{-}(t) &= x_{-\alpha_{7}}(t^{5^r})x_{-\alpha_{2}}(t^{5^r})x_{-\alpha_{5}}(t^{5^s})x_{1122100}(t^{5^r})x_{-\alpha_{1}}(t^{5^r}) x_{1112100}(t^{5^{r+1}})x_{-1010000}(t^{5^{r+1}}),
\end{align*}
and the actions of $X$ on the levels of $Q$ are:
\begin{align*}
Q/Q(2) \downarrow X &= 3^{[r]} \otimes 1^{[r+1]} + 1^{[r]} \otimes 1^{[r+1]} + 2^{[r]} \otimes 1^{[r+1]} \otimes 1^{[s]}, \\
Q(2) \downarrow X &= 2^{[r]} + 2^{[r+1]} + 1^{[r]} \otimes 1^{[s]},
\end{align*}
where the summand $3^{[r]} \otimes 1^{[r+1]}$ of $Q/Q(2)$ is generated as an $X$-module by the image of the root subgroup $U_{1111111}$. Using Corollary \ref{cor:a1cocycle}, we find that a complete set of representatives of $H^{1}(X,Q/Q(2))$ is given by
\[ \begin{array}{rl}
\gamma_{k} \, : \, x_{+}(t) \mapsto& x_{1111111}(kt^{4(5^r)})x_{0111111}(3kt^{3(5^r)})x_{1011111}(3kt^{3(5^r)})x_{1111110}(2kt^{3(5^r)})\\
& \times \ x_{0011111}(2kt^{2(5^r)})x_{1011110}(3kt^{2(5^r)})x_{0111110}(3kt^{2(5^r)})x_{0011110}(kt^{5^r})Q(2).
\end{array} \]
Now, suppose that $\gamma_{k}$ lifts to a cocycle $\phi \, : \, X \to Q$. Then $\phi(x_{+}(t))$ is a product of the above root elements with an element $q_{k,t} \in Q(2)$. If we let $x_{k}(t) = \phi(x_{+}(t))x_{+}(t)$, we find that
\[ x_{k}(t)^5 = x_{1122221}(2k^{2}t^{5^{r+1}}). \]
On the other hand, $x_{k}(t)$ is a root element in a group of type $A_{1}$, hence has order dividing $5$; it follows that $k = 0$. Thus the map $H^{1}(X,Q) \to H^{1}(X,Q/Q(2))$ is zero, so $H^{1}(X,Q) = H^{1}(X,Q(2)) = 0$ and all complements to $Q$ in $QX$ are conjugate to $X$, hence are $G$-cr.

\subsection{$L' = A_1 A_2 A_3$} \label{sec:E7A1A2A3p5}

Let $P = P_{123567} = QL$ be the unique standard $A_{1} A_{2} A_{3}$-parabolic subgroup of $G$. Let $X \cong A_1 \hookrightarrow A_1 A_2 A_3$ via $(1^{[r+1]},2^{[s]},3^{[r]})$ $(rs=0)$. This determines root elements of $X$ as in Section \ref{sec:E7A1D5p5}, and the actions of $X$ on the levels of $Q$ are as follows: 
\begin{align*}
Q/Q(2) \downarrow X &= 3^{[r]} \otimes 2^{[s]} \otimes 1^{[r+1]}, \\
Q(2)/Q(3) \downarrow X &= 4^{[r]} \otimes 2^{[s]} + 2^{[s]}, \\
Q(3)/Q(4) \downarrow X &= 3^{[r]} \otimes 1^{[r+1]}, \\
Q(4)/Q(5) \downarrow X &= 2^{[s]}. 
\end{align*}

Thus $H^{1}(X,Q(3)/Q(4)) \cong K$ for all $r$, $s$, while $H^{1}(X,Q/Q(2)) \cong K$ if $(r,s) = (0,1)$, and $H^{1}(X,Q/Q(2)) = 0$ otherwise. Moreover, by Lemma \ref{lem:h2fora1} if $s > 0$ then there are direct summands in levels $2$ and $4$ with non-vanishing second cohomology group. With similar calculations to Section \ref{sec:E7A1D5p5}, we find that every cocycle $X \to Q(3)/Q(4)$ lifts to a cocycle $X \to Q(3)$, while if $(r,s) = (0,1)$, a cocycle $X \to Q/Q(2)$ lifts to a cocycle $X \to Q/Q(3)$ if and only if it is a coboundary. This shows that $\mathbb{V} \cong H^{1}(X,Q(3)/Q(4)) \cong K$, for any choice of $r$ and $s$.

Therefore, for each $r$ and $s$, elements of $H^{1}(X,Q)$ are parametrised by elements of $K$. Applying Lemma \ref{lem:torus}, for each possible $(r,s)$ we have at most one conjugacy class of non-$G$-cr complements to $Q$ in $QX$.     

Let $Y \cong A_1 \hookrightarrow A_2 A_5$ via $(2^{[s]},W(5)^{[r]})$ where $rs=0$. By Lemma \ref{lem:BMR}, $Y$ is non-$G$-cr. Moreover, $Y$ is contained in a $A_1 A_2 A_3$-parabolic subgroup of $A_2 A_5$, and hence in such a parabolic subgroup of $G$. Thus there is exactly one conjugacy class of non-$G$-cr complements to $Q$ in $QX$ for each pair $(r,s)$ with $rs = 0$.

From Table \ref{E7p5tab}, we find that $\text{dim}(C_{L(G)}(Y)) = 0$ and hence $C_G(Y)^\circ =1$. 

\subsection{$L'=D_5$} \label{E7p5D5}

The two standard $D_5$-parabolic subgroups of $G$ are $P_{12345}$ and $P_{23456}$. Let $X \cong A_1 < D_5$ via $4^{[r]} + 1^{[r+1]} \otimes 1^{[s]} + 0$ $(rs = 0;$ $r+1 \neq s)$. Let $P = P_{12345} = QL$. The action of $X$ each level of $Q$ is as follows: 
\begin{align*}
Q/Q(2) \downarrow X &= 3^{[r]} \otimes 1^{[r+1]} + 3^{[r]}  \otimes 1^{[s]} + 0, \\
Q(2)/Q(3) \downarrow X &= 3^{[r]} \otimes 1^{[r+1]} + 3^{[r]}  \otimes 1^{[s]}, \\
Q(3) \downarrow X &= 4^{[r]} + 1^{[r+1]} \otimes 1^{[s]} + 0.
\end{align*}
Let $U$ and $V$ be the respective summands of $Q/Q(2)$ and $Q(2)/Q(3)$ which are isomorphic to $3^{[r]} \otimes 1^{[r+1]}$. By Lemma \ref{lem:h1fora1}, $\mathbb{V}_{X,Q} = H^{1}(X,U) \oplus H^{1}(X,V) \cong K^{2}$. The modules $U$ and $V$ are respectively generated by the images of the root subgroups $U_{\alpha_{6}}$ and $U_{\alpha_{6} + \alpha_{7}}$. The root group $U_{\alpha_{7}}$ commutes with $L'$, hence gives a $1$-dimensional subgroup of $Q^{X}$, and does not commute with $U_{\alpha_{6}}$. Hence, as described in Section \ref{sec:nonab}, conjugation by $v \stackrel{\textup{def}}{=} x_{\alpha_{7}}(1)$ induces a non-trivial homomorphism of $L'$-modules $c_{v,1} \, : \, Q/Q(2) \to Q(2)/Q(3)$, which restricts to an isomorphism $U \to V$.

We now apply Corollary \ref{cor:trivs}. This tells us that, with respect to a basis of $\mathbb{V}$ consisting of a non-zero element from each of $H^{1}(X,U)$ and $H^{1}(X,V)$, elements of $H^{1}(X,Q)$ are parametrised by pairs $(k_1,k_2) \in \mathbb{V}$ with $k_1 k_2 = 0$.

Next, the element $n_{7}$ of $N_{G}(T)$ normalises each root subgroup in $L'$, and swaps $U_{\alpha_{6}}$ with $U_{\alpha_{6} + \alpha_{7}}$. It follows that two complements to $Q$ in $QX$ corresponding to $(k_1,k_2)$ and $(k_2,k_1)$ are conjugate in $G$. Together with the previous paragraph, this means we can assume that $(k_1,k_2) = (k_1,0)$. By considering the non-trivial action of $Z(L)$ on each summand in $Q$, we therefore deduce that there is at most one class of non-$G$-cr complements to $Q$ in $QX$.

Now consider $P_{23456}$. We claim that a non-$G$-cr subgroup $A_1$ with irreducible image in $L_{23456}$ must be $G$-conjugate to a subgroup of $P_{12345}$. It suffices to exhibit an element of $N_{G}(T)$ which sends the roots subgroups in $L'_{23456}$ to those of $L'_{12345}$, and also sends $U_{\alpha_1}$ to $U_{\alpha_6}$, since the image of $U_{\alpha_{1}}$ generates the $L_{23456}$-module of high weight $\lambda_{4}$ in $Q_{23456}/Q_{23456}(2)$. The following element fits the bill:
\[ n_7 n_6 n_5 n_4 n_3 n_2 n_4 n_5 n_6 n_7 n_1 n_3 n_4 n_5 n_6 n_2 n_4 n_5 n_3 n_4 n_1 n_3 n_2 n_4 n_5 n_6 n_7 \]
and hence up to $G$-conjugacy there is at most one non-$G$-cr subgroup $A_1$ in a $D_5$-parabolic subgroup of $G$, with irreducible image in the Levi factor.

Recall that $A_{1} A_{5}'$ denotes a subgroup of this type lying in a Levi $E_{6}$ subgroup of $G$, and let $Y \cong A_1 \hookrightarrow A_1 A_5'$ via $(1^{[s]},W(5)^{[r]})$, where $s \neq r + 1$. From Section \ref{sec:E6D5} we know that $Y$ is non-$E_6$-cr, lying in a $D_5$-parabolic subgroup of $E_6$, with irreducible image in the Levi factor. Hence by Lemma \ref{lem:BMR}, $Y$ is non-$G$-cr, and thus $Y$ is a representative of the unique class of non-$G$-cr subgroups of $QX$.

From Table \ref{E7p5tab}, we find that $\text{dim}(C_{L(G)}(Y)) = 3$, and since $Y < E_{6}$ we know that $C_G(Y)$ contains a 1-dimensional torus. We now construct a $2$-dimensional unipotent subgroup of $Q$ centralising a conjugate of $Y$; it then follows that $C_G(Y)^{\circ}$ has the form $U_{2}T_{1}$.

First note that $Z(Q)$ contains a trivial $X$-submodule, which therefore centralises $QX$, and in particular centralises every complement to $Q$ in $QX$. Moreover the root subgroup $U_{\alpha_7}$ commutes with $Q(2)$, since there is no root $\alpha$ of level $2$ such that $\alpha + \alpha_7$ is also a root. Therefore, if we can show that some non-$G$-cr complement $Z$ to $Q$ in $QX$ is contained in $Q(2)X$, it follows that $U_{\alpha_7}$ also lies in $C_{G}(Z)^{\circ}$.

The natural map $H^{0}(X,Q) \to H^{0}(X,Q/Q(2))$ is surjective since $H^{0}(X,Q) = U_{\alpha_{7}}$ has trivial intersection with $Q(2)$. From the long exact sequence of cohomology, it follows that we have an exact sequence $0 \to H^{1}(X,Q(2)) \to H^{1}(X,Q)$ of pointed sets, so a non-zero element of $H^{1}(X,Q(2))$ gives rise to a non-zero element of $H^{1}(X,Q)$. Thus there exists a non-$G$-cr complement $Z$ to $Q(2)$ in $Q(2)X$, and as described in the previous paragraph we have $C_{G}(Z)^{\circ} = U_{2}T_{1}$.

\subsection{$L'=A_1 D_4$}

Let $P = P_{23457} = QL$ be the unique standard $A_{1} D_{4}$-parabolic subgroup of $G$. Let $X,Y \cong A_1 \hookrightarrow A_1 D_4$ via $(1^{[r]},3^{[s]} \otimes 1^{[s+1]})$ $(rs=0)$ where $V_{D_{4}}(\lambda_4) \downarrow X = V_{D_{4}}(\lambda_3) \downarrow Y = 4^{[s]} + 2^{[s+1]}$. The levels of $Q$ are as follows:
\begin{align*}
Q/Q(2) \downarrow X &= (3^{[s]} \otimes 1^{[s+1]}) + (1^{[r]} \otimes 4^{[s]}) + (1^{[r]}  \otimes 2^{[s+1]}), \\
Q(2)/Q(3) \downarrow X &=  (1^{[r]} \otimes 3^{[s]} \otimes 1^{[s + 1]}) + 0, \\
Q(3)/Q(4) \downarrow X &= 3^{[s]} \otimes 1^{[s+1]}, \\
Q(4) \downarrow X &= 0.
\end{align*}

By Lemma \ref{lem:h1fora1}, for any $r$ and $s$ we have $\mathbb{V}_{X,Q} \cong K^{2}$. The trivial module in $Q(2)/Q(3)$ is generated by the image of the root subgroup $U_{0112221}$, which commutes with $X$, hence induces a homomorphism $Q/Q(2) \to Q(3)/Q(4)$ of $X$-modules, which is non-trivial since $x_{0112221}(1)$ does not commute with all of $Q/Q(2)$, for instance it does not commute with $x_{\alpha_{1}}(1)$. This must therefore induce an isomorphism from the summand $3^{[s]} \otimes 1^{[s+1]}$ to $Q(3)/Q(4)$.

Applying Corollary \ref{cor:trivs}, complements to $Q$ in $QX$ are parametrised by pairs $(k_1,k_2)$ of elements of $K$ with $k_1 k_2 = 0$. In addition, the element $n_{0112221}$ of $N_{G}(T)$ normalises each root subgroup in $L'$ and swaps $U_{\alpha_{1}}$ and $U_{1112221}$. It follows that any complements corresponding to $(k,0)$ and $(0,k)$ are $G$-conjugate. Lastly, the non-trivial action of $Z(L)$ reduces us to the cases $(k_1,k_2) = (1,0)$ and $(0,0)$, by Lemma \ref{lem:torus}. Hence we have at most one $G$-conjugacy class of non-$G$-cr complements to $Q$ in $QX$. 

With similar reasoning, there is at most one $G$-conjugacy class of non-$G$-cr complements to $Q$ in $QY$. Furthermore, the Weyl group element $n_{0112221} n_{2234321} n_{1} n_{3} n_{4} n_{2} n_{5} n_{4} n_{3} n_{1} \in N_{G}(T)$ normalises $U_{\alpha_{1}}$ and induces an outer automorphism of $L'$. Hence a lift of this element to $N_{G}(T)$ conjugates each complement to $Q$ in $QY$ to a subgroup of $QX$, and so up to conjugacy in $G$ there is at most one non-$G$-cr subgroup $A_{1}$ with irreducible image in a Levi factor $A_{1} D_{4}$.

Let $Z \cong A_1 \hookrightarrow A_1 D_5$ via $(1^{[r]},T(8)^{[s]})$. By Lemma \ref{lem:BMR}, $Z$ is non-$G$-cr. Since $Z$ is contained in an $A_1 D_4$-parabolic subgroup of $A_1 D_5$, it is contained in an $A_1 D_4$-parabolic subgroup of $G$. Hence $Z$ is a representative for the unique $G$-conjugacy class of non-$G$-cr subgroups $A_1$ in $P$.

For each $r$ and $s$, $\text{dim}(C_{L(G)}(Z)) = 1$, and so $C_G(Z)^\circ = C_G(A_1 D_5)^\circ = T_1$. 

\subsection{$L'=A_2 A_3$}

The three standard $A_{2} A_{3}$-parabolic subgroups of $G$ are $P_{13567}$, $P_{13467}$ and $P_{23467}$. Let $X \cong A_1 \hookrightarrow A_2 A_3$ via $(2,1 \otimes 1^{[1]})$. First let $P = P_{13567} = QL$. Then the actions of $X$ on the levels of $Q$ are as follows: 
\begin{align*}
Q/Q(2) \downarrow X &= 3 \otimes 1^{[1]} + 1 \otimes 1^{[1]} + 0, \\
Q(2)/Q(3) \downarrow X &= 3 \otimes 1^{[1]} + 1 \otimes 1^{[1]}, \\
Q(3)/Q(4) \downarrow X &= 4 + 2 + 2 \otimes 2^{[1]} + 0, \\
Q(4)/Q(5) \downarrow X &= 1 \otimes 1^{[1]}, \\
Q(5)/Q(6) \downarrow X & = 1 \otimes 1^{[1]}, \\
Q(6)/Q(7) \downarrow X & = 2. 
\end{align*}

By Lemma \ref{lem:h1fora1}, we have $\mathbb{V} \cong K^2$. The root group $U_{\alpha_{2}}$ gives rise to the trivial $L'$-submodule in $Q/Q(2)$, and does not commute with the root group $U_{\alpha_{4}}$. A non-trivial element of $U_{\alpha_{2}}$ therefore induces a non-trivial homomorphism of $L'$-modules $Q/Q(2) \to Q(2)/Q(3)$, and hence by Corollary \ref{cor:trivs}, parametrising complements to $Q$ in $QX$ by pairs $(k_1,k_2)$ of elements of $K$, we may assume that $k_{1} k_{2} = 0$. Further, using Lemma \ref{lem:torus}, the action of the torus $Z(L)$ reduces us to the cases $(k_1,k_2) = (1,0)$ or $(0,1)$. Lastly, the element $n_2 \in N_{G}(T)$ normalises the root subgroups in $L'$ and swaps $U_{\alpha_4}$ and $U_{\alpha_2 + \alpha_4}$; hence we may assume $(k_1,k_2) = (0,1)$ or $(0,0)$, and we have at most one conjugacy class of non-$G$-cr complements to $Q$ in $QX$.

Similar calculations hold for $P_{13467} = Q_{13467}L_{13467}$. We have $\mathbb{V} \cong K^{2}$, the modules with non-vanishing first cohomology group appearing in levels $1$ and $3$, respectively generated by the images of $U_{\alpha_5}$ and $U_{0112210}$. The subgroup $X < L_{13467}$ commutes with elements of $Q_{13467}/Q_{13467}(4)$ of the form
\[ x_{1112100}(t) x_{1111110}(t) x_{0112110}(t) x_{0111111}(t) Q_{13467}(4) \]
and these induce non-zero $X$-module homomorphisms $Q_{13467}/Q_{13467}(2) \to Q_{13467}(3)/Q_{13467}(4)$. An element of $N_{G}(T)$ normalising $L'$ and fusing complements corresponding to $(k_1,0)$ and $(0,k_1) \in \mathbb{V}$ is $n_{1223221}n_{1122100}n_{0101100}$. Moreover, the element $n_{4} n_{3} n_{1} n_{5} n_{4} n_{3} n_{6} n_{5} n_{4} n_{7} n_{6} n_{5}$ conjugates the root subgroups in $L_{13467}$ to those in $L_{13567}$, and also conjugates $U_{0112210}$ to a subgroup of $Q_{13567}$. It follows that each complement to $Q_{13467}$ in $Q_{13467}X$ is $G$-conjugate to a subgroup of $P_{13567}$.

In $P_{23467} = Q_{23467}L_{23467}$ we again have $\mathbb{V} \cong K^{2}$, the modules with non-vanishing first cohomology group appearing in levels $2$ and $3$, generated by the images of $U_{1011100}$ and $U_{1112210}$. The element $n_{1111000} \in N_{G}(T)$ conjugates the root subgroups in $L_{23467}$ to those in $L_{13467}$, and sends both $U_{1011100}$ and $U_{1112210}$ to subgroups of $Q_{13467}$. Thus each non-$G$-cr subgroup of $P_{23467}$ is $G$-conjugate to a subgroup of $P_{13467}$, hence is also conjugate to a subgroup of $P_{13567}$.

Up to $G$-conjugacy, there is therefore at most one non-$G$-cr subgroup $A_{1}$ contained in an $A_{2}A_{3}$-parabolic of $G$ with irreducible image in the Levi factor. Let $Y \cong A_1 < A_6$ via $W(6)$. Then $Y$ is non-$G$-cr by Lemma \ref{lem:BMR} and is contained in an $A_2 A_3$-parabolic subgroup of $A_6$ and hence of $G$. Therefore, $Y$ is a representative of this conjugacy class of subgroups.

Now, from Table \ref{E7p5tab}, $C_{L(G)}(Y)$ is $3$-dimensional. Also, $C_G(Y)^{\circ}$ contains the $1$-dimensional torus $C_G(A_{6})^{\circ}$. Moreover with similar considerations to those in Section \ref{E7p5D5} we see that a complement $X_{[0,1]} \le Q(2)L_{13567}$ commutes with $U_{\alpha_2}$ as well as with the 1-dimensional subgroup giving rise to the trivial submodule in $Q(3)/Q(4)$; hence $C_{G}(X_{[0,1]})^{\circ} = U_{2}T_{1}$

\subsection{$L'=A_1^2 A_3$}

The three standard $A_{1}^{2} A_{3}$-parabolic subgroups of $G$ are $P_{12567}$, $P_{12457}$ and $P_{23567}$. Let $X,Y \hookrightarrow A_1^2 A_3$ via $(1^{[r]},1^{[s+1]},3^{[s]})$ and $(1^{[s+1]},1^{[r]},3^{[s]})$ respectively ($rs=0$ in both cases), and let $Z \cong A_1 \hookrightarrow A_1^2 A_3$ via $(1,1,1 \otimes 1^{[1]})$. First consider $P_{12567} = P = QL$. Then $L'$ acts on the levels of $Q$ as follows: 
\begin{align*}
Q/Q(2) \downarrow L' &= (1,0,000) + (0,1,100), \\
Q(2)/Q(3) \downarrow L' &= (1,1,100), \\
Q(3)/Q(4) \downarrow L' &= (1,0,010), \\
Q(4)/Q(5) \downarrow L' &= (0,0,010), \\
Q(5)/Q(6) \downarrow L' & = (0,1,001), \\
Q(6)/Q(7) \downarrow L' & = (0,0,000), \\
Q(7) \downarrow L' &= (1,0,000). 
\end{align*}
From Lemma \ref{lem:h1fora1} it follows that $\mathbb{V}_{X,Q} \cong K^2$, $\mathbb{V}_{Y,Q} = 0$, and $\mathbb{V}_{Z,Q} \cong K$. Both $X$ and $Z$ have summands in some level with a non-vanishing second cohomology group. We first analyse complements to $Q$ in $QX$.

Note first that $Q(4)/Q(5) \downarrow X = 4^{[s]} + 0$. We can identify the following positive and negative root elements of $X$:
\begin{align*}
x_{+}(t) &= x_{\alpha_{1}}(t^{5^r}) x_{\alpha_{2}}(t^{5^{s+1}}) x_{\alpha_{5}}(t^{5^s}) x_{\alpha_{6}}(t^{5^{s}}) x_{\alpha_{7}}(t^{5^{s}})\\ &\quad \times x_{0000110}(2t^{2(5^{s})}) x_{0000011}(2t^{2(5^{s})}) x_{0000111}(3t^{3(5^{s})}),\\
x_{-}(t) &= x_{-\alpha_{1}}(t^{5^r}) x_{-\alpha_{2}}(t^{5^{s+1}}) x_{-\alpha_{5}}(3t^{5^s}) x_{-\alpha_{6}}(4t^{5^{s}}) x_{-\alpha_{7}}(3t^{5^{s}})\\&\quad \times x_{-0000110}(t^{2(5^{s})}) x_{-0000011}(t^{2(5^{s})}) x_{-0000111}(3t^{3(5^{s})}),
\end{align*}
and these elements commute with every element of $Q(4)$ of the form $x_{1122210}(u) x_{1122111}(u)$. If $u \neq 0$ then this element does not commute with $U_{\alpha_{7}}$ and therefore induces a non-trivial homomorphism of $X$-modules $Q/Q(2) \to Q(5)/Q(6)$. Therefore, applying Corollary \ref{cor:trivs} and Lemma \ref{lem:torus}, there are at most two conjugacy classes of non-$G$-cr complements to $Q$ in $QX$, corresponding to $(1,0)$ and $(0,1) \in \mathbb{V}_{X,Q}$. Furthermore, the element $n_{1122100}n_{1122221}n_{6}n_{5}n_{7}n_{6} \in N_{G}(T)$ stabilises the set of root subgroups in $L'$, inducing an outer automorphism on the $A_{3}$ factor. This element also swaps the root subgroups $U_{\alpha_{4}}$ and $U_{1123210}$. Since the images of these root groups generate the two $X$-modules with non-vanishing first cohomology group, a lift of this Weyl group element exchanges complements to $Q$ in $QX$ corresponding to $(1,0)$ and $(0,1) \in \mathbb{V}$.

Since $\mathbb{V}_{Y,Q} = 0$, all complements to $Q$ in $Q Y$ are $Q$-conjugate to $Y$. Also, with similar calculations to those of Section \ref{sec:E7A1D5p5}, we find that a cocycle $Z \to Q(2)/Q(3)$ lifts to a cocycle $Z \to Q(2)$ only if it is a coboundary. Hence every complement to $Q$ in $QZ$ is $Q$-conjugate to $Z$.

Now consider $P_{12457}$. With identical arguments, we find that there are no non-$G$-cr complements to $Q_{12457}$ in $Q_{12457}Y$ or $Q_{12457}Z$, and that there is exactly one $G$-conjugacy class of non-$G$-cr complements to $Q_{12457}$ in $Q_{12457}X$. Moreover, the element $n_{1234321} n_{1} n_{3} n_{4} n_{2} n_{5} n_{4} n_{3} n_{6} n_{5} n_{4} n_{2} n_{7} n_{6} n_{5} n_{4} n_{3} \in N_{G}(T)$ conjugates the root subgroups of $L_{12457}$ to those in $L_{12567}$, and sends the subgroup $X < L_{12457}$ to the corresponding subgroup  of $L_{12567}$. This element also conjugates $U_{\alpha_6}$, which gives rise to the module in $Q_{12457}$ with non-vanishing first cohomology group, to $U_{\alpha_{4}}$, and hence complements to $Q$ in $Q_{12457}X$ are $G$-conjugate to subgroups of $P_{12567}$.

An entirely similar argument shows that for $P_{23567}$ there are no non-$G$-cr complements to $Q_{23567}$ in $Q_{23567}X$ or $Q_{23567}Z$, and there is exactly one $G$-conjugacy class of non-$G$-cr complements to $Q_{23567}$ in $Q_{23567}Y$. Moreover, the element $n_{3} n_{1} \in N_{G}(T)$ conjugates $L'_{23567}$ to $L'_{12567}$, sending $Y < L'_{23567}$ to the subgroup $X < L'_{12567}$. This also fixes $U_{1123210}$, hence complements to $Q_{23567}$ in $Q_{23567}Y$ are also conjugate to a subgroup of $P_{12567}$.

Therefore, for each $r$ and $s$, there is at most one $G$-conjugacy class of non-$G$-cr subgroups $A_{1}$ with irreducible image in a Levi factor of type $A_{1}^{2} A_{3}$. Let $A \cong A_1 \hookrightarrow A_1 A_5$ via $(1^{[r]},W(5)^{[s]})$. Then $A$ is non-$G$-cr by Lemma \ref{lem:BMR} and is contained in an $A_1^2 A_3$-parabolic subgroup of the Levi $A_1 A_5$ and hence in an $A_1^2 A_3$-parabolic subgroup of $G$. It is also clear from Table \ref{E7p5tab} that distinct choices of $(r,s)$ lead to non-conjugate subgroups of $G$, and thus each $G$-conjugacy class above exists.

Finally, the complement $X_{[0,1]}$ to $Q$ in $QX < P_{12567}$ above commutes with the 2-dimensional subgroup $\left< U_{1224321},\ x_{1122210}(t) x_{1122111}(t) \, : \, t \in K \right>$. Since $A$ lies in a Levi subgroup $A_{1} A_{5}$, we see $C_G(A)^{\circ}$ contains a $1$-dimensional torus, hence $C_G(A)^\circ = U_2 T_1$. 

\subsection{$L'=A_1^3 A_2$}

Let $P = P_{12357} = QL$ be the unique standard $A_{1}^{3}A_{2}$-parabolic subgroup of $G$. Let $X \cong A_1 \hookrightarrow A_1^3 A_2$ via $(1^{[r]},1^{[s]},1^{[t]},2^{[u]})$. The action of $L'$ on the levels of $Q$ is as follows: 
\begin{align*}
Q/Q(2) \downarrow L' &= (1,1,0,01) + (0,1,1,00), \\
Q(2)/Q(3) \downarrow L' &= (1,0,1,01) + (0,0,0,10), \\
Q(3)/Q(4) \downarrow L' &= (0,1,1,10), \\
Q(4)/Q(5) \downarrow L' &= (0,0,0,10) + (1,0,1,00), \\
Q(5)/Q(6) \downarrow L' & = (1,1,0,00), \\
Q(6) \downarrow L' & = (0,0,0,01).
\end{align*}
We may assume $rstu = 0$. If $u = r = s - 1$ or $u = s = r - 1$ then $H^{1}(X,Q/Q(2)) \cong K$, and similarly if $u = r = t - 1$ or $u = t = r - 1$ then $H^{1}(X,Q(2)/Q(3)) \cong K$, and if $u = s = t - 1$ or $u = t = s - 1$ then $H^{1}(X,Q(3)/Q(4)) \cong K$. Each of these cohomology groups vanishes if neither of the two corresponding conditions holds, and no other $X$-module occurring can have non-zero first cohomology group. Hence $\mathbb{V} \neq 0$ if and only if $(r,s,t)$ is a permutation of $(u,u+1,a)$ for some $a \in K$. In this case, $\mathbb{V} \cong K$ unless $a = u$ or $a = u + 1$, in which case $\mathbb{V} \cong K^2$ and $u = 0$.

If $\mathbb{V} \cong K$ then by Lemma \ref{lem:torus} there is at most one $G$-conjugacy class of non-$G$-cr complements to $Q$ in $QX$.

Now suppose $\mathbb{V} \cong K^{2}$, so that $u = 0$ and $(r,s,t)$ is a permutation of $(1,0,0)$ or $(1,1,0)$. If $s = t$, or if $s \neq t$ but $r = s = 0$, then there is a trivial $X$-submodule in level $1$, inducing non-zero $X$-module homomorphisms $Q/Q(2) \to Q(2)/Q(3)$ and $Q(2)/Q(3) \to Q(3)/Q(4)$. If $s \neq t$ and $r = s = 1$, then level $5$ contains an $X$-direct summand $2^{[1]}$, which has non-vanishing second cohomology group. Finally, if $s \neq t$ and $r \neq s$ then either $r = t = 1$ and level $4$ contains an $X$-direct summand $2^{[1]}$, or $r = t = 0$ and level $2$ contains a trivial submodule, inducing a non-zero homomorphism of $X$-modules $Q/Q(2) \to Q(3)/Q(4)$.

With similar calculations to previous sections, using the above paragraph we find that whenever $\mathbb{V} \cong K^{2}$, so that complements are parametrised by $(k_1,k_2) \in K^{2}$, then either $k_1 k_2 = 0$ is necessary for the partial map $\rho \, : \, \mathbb{V} \to H^{1}(X,Q)$ to be defined, or Corollary \ref{cor:trivs} applies and a complement to $Q$ in $QX$ corresponding to $(k_1,k_2)$ is $Q$-conjugate to one corresponding to either $(k_1,0)$ or $(0,k_2)$. Applying Lemma \ref{lem:torus} reduces us to $(k_1,k_2) = (1,0)$ or $(0,1)$ or $(0,0)$.

Next, the element $n_{6} n_{5} n_{7} n_{6} \in N_{G}(T)$ stabilises the set of root subgroups contained in $L$, swapping those in the second and third $A_1$ factors, and also swaps $U_{\alpha_4}$ and $U_{0001110}$, which give rise to the modules in levels $1$ and $2$ with non-vanishing first cohomology group. Similarly, the element $n_{0112100} n_{1} n_{3} n_{4} n_{2} n_{5} n_{4} $ swaps the root subgroups in the first two $A_1$ factors and swaps $U_{0001110}$ and $U_{0112110}$ which give rise to the relevant modules in levels $2$ and $3$, and finally the element $n_{0112221} n_{1011110} n_{1111111} n_{1} n_{2} n_{3} n_{7}$ swaps the root subgroups in the first and third $A_{1}$ factors, and swaps the appropriate root subgroups in levels $1$ and $3$. It follows that for each $(r,s,t,u)$ there is at most one non-$G$-cr complement to $Q$ in $QX$, up to conjugacy in $G$, and furthermore, the six potential non-$G$-cr subgroups corresponding to permutations of $(r,s,t)$ are all conjugate in $G$. Hence there exists at most one non-$G$-cr subgroup $A_1$ of $P$ with irreducible image in the Levi factor, for each set of twists $(u,u+1,a,u)$ with $ua = 0$.

Let $Y \cong A_1 \hookrightarrow A_1 D_5$ via $(1^{[a]},T(6)^{[u]})$, where $ua = 0$. By Lemma \ref{lem:BMR}, $Y$ is non-$G$-cr. Moreover, $Y$ is contained in an $A_1^3 A_2$-parabolic subgroup of $A_1 D_5$ and hence of $G$. Therefore, $Y$ is a representative of the conjugacy class found in the above analysis.

From Table \ref{E7p5tab} we see that $C_{L(G)}(Y))$ has dimension $1$ if $a \neq u, u+1$; $2$ if $a = u+1$; and $3$ if $a = u$. In each case, $C_G(Y)^\circ \ge C_G(A_1 D_5)^\circ = T_1$.

If $X$ above corresponds to twists $(0,1,1,0)$ then we check that $X$ centralises the abelian group $\{ x_{0000110}(c) x_{0000011}(c) \, : \, c \in K \}$, which also commutes with the root groups giving rise to the module in level $2$ with non-vanishing first cohomolgy group. If instead $X$ corresponds to twists $(0,1,0,0)$ we check that $X$ centralises the 2-dimensional unipotent group
\[ 
\left\{
\begin{array}{l}
x_{1 0 1 1 1 1 0}(c) x_{0 1 1 1 1 1 0}(-c) x_{0 1 0 1 1 1 1}(3c) x_{0 0 1 1 1 1 1}(c) x_{1 1 1 2 2 2 1}(c^2)\\
\quad \times \ x_{1 2 2 3 2 1 0}(d) x_{1 1 2 3 2 1 1}(d)
\end{array}
\, : \, c,d \in K
\right\}
\]
which also commutes with the root subgroups giving rise to the $X$-module in level $3$ with non-zero first cohomology group.

It follows that $C_G(Y)^{\circ} = T_1$ if $a \neq u, u + 1$; $C_G(Y)^{\circ} = U_{1} T_{1}$ if $a = u + 1$; and $C_G(Y)^{\circ} = U_{2} T_{1}$ if $a = u$.

\subsection{$L'=D_4$}
Let $P = P_{2345} = QL$ be the unique standard $D_{4}$-parabolic subgroup of $G$. Let $X$, $Y$, $Z$ be representatives of the three $L'$-conjugacy classes of $L'$-irreducible subgroups, with $X$ (resp.\ $Y$, $Z$) acting on $V_{D_4}(\lambda_1)$ (resp.\ $V_{D_4}(\lambda_3)$, $V_{D_4}(\lambda_4)$) as $4 + 2^{[1]}$ (in this case, the remaining $8$-dimensional modules $V_{D_4}(\lambda_i)$ restrict as $3 \otimes 1^{[1]}$). The action of $D_4$ on the levels of $Q$ is as follows: 
\begin{align*}
Q/Q(2) \downarrow L' &= \lambda_3 + \lambda_4 + 0, \\
Q(2)/Q(3) \downarrow L' &= \lambda_1 + \lambda_4, \\
Q(3)/Q(4) \downarrow L' &= \lambda_1 + 0, \\
Q(4)/Q(5) \downarrow L' &= \lambda_3, \\
Q(5)/Q(6) \downarrow L' & = 0. \\
\end{align*}
Reading left-to-right and down the levels, the non-trivial $L'$-modules are generated by the images of the root groups $U_{\alpha_1}$, $U_{\alpha_{6}}$, $U_{1011110}$, $U_{0000011}$, $U_{1011111}$ and $U_{1112221}$. The element $w \stackrel{\textup{def}}{=} n_{6} n_{5} n_{4} n_{2} n_{3} n_{1} n_{4} n_{3} n_{5} n_{4} n_{2} n_{6} n_{5} n_{4} n_{3} n_{1} \in N_{G}(T)$ induces a triality automorphism on $L'$ and sends $U_{0000011} \to U_{1112221} \to U_{1011111} \to U_{0000011}$. This element will shortly be used to show that all non-$G$-cr subgroups $A_1$ of $QL$ are conjugate to a subgroup of $QX$.

Now, $\mathbb{V}_{X,Q} \cong \mathbb{V}_{Y,Q} \cong \mathbb{V}_{Z,Q} \cong  K^4$. The trivial summands in levels $1$ and $3$ are respectively generated by the root groups $U_{\alpha_{7}}$ and $U_{0112221}$, inducing isomorphisms between the $L'$-modules of high weight $\lambda_4$ in levels $1$ and $2$ and between the modules of high weight $\lambda_1$ in levels $2$ and $3$, and between the modules of high weight $\lambda_3$ in levels $1$ and $4$.

First consider $X$. Then $H^{1}(X,Q/Q(2)) \cong K^{2}$, and $H^{1}(X,Q(2)/Q(3)) \cong H^{1}(X,Q(4)/Q(5)) \cong K$, and also $H^{2}(X,Q(2)/Q(3)) \cong H^{2}(X,Q(3)/Q(4)) \cong K$. With similar calculations to Section \ref{sec:E6D4}, we find that a cocycle corresponding to $(k_1,k_2) \in H^{1}(X,Q/Q(2))$ lifts to a cocycle $X \to Q/Q(3)$ only if $k_1 k_2 = 0$, and also a cocycle $X \to Q/Q(3)$ corresponding to $(k_1,k_2,k_3)$ lifts to $Q/Q(4)$ only if $k_1 k_3 = 0$. Moreover, using the isomorphisms of $L'$-modules coming from $U_{\alpha_{7}}$ and $U_{0112221}$, if cocycles $X \to Q$ are parametrised by $(k_1,k_2,k_3,k_4) \in \mathbb{V}$ then we may assume that $k_2k_3 = k_1 k_4 = 0$.

Next, the element $n_7 \in N_{G}(T)$ normalises the root subgroups in $L'$ and the positive root groups giving rise to $L$-modules of high weight $\lambda_3$ in the filtration of $Q$, while exchanging the root groups giving rise to the two $L$-modules of high weight $\lambda_4$. Hence a complement to $Q$ in $QX$ corresponding to $(k_1,k_2,k_3,k_4) \in \mathbb{V}$ is $G$-conjugate to one corresponding to $(k_1,k_3,k_2,k_4)$. Similarly, the element $n_{0112221}$ normalises the root subgroups in $L'$ and swaps the positive root groups which give rise to the two $X$-modules of high weight $\lambda_3$, and so a complement corresponding to $(0,0,0,k_4)$ is $G$-conjugate to one corresponding to $(k_4,0,0,0)$. Finally, the element $n_{1122210} n_{1112110}$ normalises $L'$ and $U_{1011111}$, and sends $U_{0000011}$ to $U_{1122221}$ (level $4$) and sends $U_{1112221}$ to $U_{0000111}$ (level $2$). It follows that complements to $X$ in $QX$ corresponding to $(0,0,k_3,k_4)$ and $(0,0,k_4,k_3)$ are $G$-conjugate.

Putting this together and applying Lemma \ref{lem:torus}, it follows that there are at most two non-$G$-cr complements to $Q$ in $QX$ up to $G$-conjugacy, corresponding to $(k_1,k_2,k_3,k_4) = (0,0,0,1)$ and $(0,0,1,1)$. Similar reasoning gives the same conclusion for $QY$ and $QZ$, hence each of these contains at most two non-$G$-cr complements to $Q$ up to conjugacy in $G$. Moreover, these complements each correspond to 4-tuples $(k_1,k_2,k_3,k_4)$ of elements of $K$ with at most two non-zero entries, thus these complements lie in the subgroup $\left<U_{\alpha},U_{\beta},L\right>$ for two appropriate positive roots $\alpha, \beta \in \{ \alpha_1,\alpha_6,1011110,0000011,1011111,1112221 \}$. The element $w$ above now conjugates $X \to Z \to Y$ and also permutes the positive root groups appropriately, and hence each complement to $Q$ in $QY$ or $QZ$ is $G$-conjugate to a subgroup of $QX$. Thus there are at most two $G$-conjugacy classes of non-$G$-cr subgroups $A_1$ in a $D_{4}$-parabolic subgroup with irreducible image in the Levi factor.

Let $A \cong A_1 < D_5$ via $T(8)$ and $B \cong A_1 \hookrightarrow A_1 A_5'$ via $(1^{[1]},W(5))$ (recall that $A_1 A_5'$ is the conjugacy class of $A_1 A_5$ subgroups contained in $E_6$). Then $A$ and $B$ are non-$G$-cr by Lemma \ref{lem:BMR}. Moreover, $A$ is contained in a $D_4$-parabolic subgroup of $D_5$ and hence in a $D_{4}$-parabolic subgroup of $G$, and $B$ is contained in a $D_4$-parabolic subgroup of $E_6$ by Theorem \ref{thm:e6} and hence in a $D_{4}$-parabolic subgroup of $G$. From Table \ref{E7p5tab}, we find that $A$ and $B$ are not $GL_{56}$-conjugate and hence are not $G$-conjugate. Therefore, $A$ and $B$ are representatives of the two $G$-conjugacy classes of non-$G$-cr subgroups in $D_4$-parabolic subgroups.   

From Table \ref{E7p5tab}, $\text{dim}(C_{L(G)}(A)) = 6$ and $\textup{dim}(C_{L(G)}(B)) = 4$. Also $C_G(A)^{\circ} \ge C_G(D_5)^{\circ} = A_{1} T_{1}$, and $C_G(B)^{\circ} \ge C_G(A_{1}A_{5}') = T_{1}$. Now, in the above analysis we have shown that a complement to $Q$ in $QX$ is $G$-conjugate to a subgroup of $\left< U_{0000011},\ U_{1011111},\ X \right>$. But each of these three groups commutes with each of the positive root groups $U_{\alpha_7}$, $U_{0112221}$ and $U_{2234321}$, so each complement to $X$ in $QX$ must centralise a 3-dimensional unipotent group generated by positive root elements. It follows that $C_G(A)^{\circ} = U_{2} A_{1} T_{1}$ and $C_G(B)^{\circ} = U_{3} T_{1}$.

\subsection{$L'=A_1 A_3$}

There are $11$ standard $A_1 A_3$-parabolic subgroups of $G$. The corresponding Levi factors fall into two conjugacy classes; nine standard parabolic subgroups have Levi factor conjugate to $L_{1567}$, and two have Levi factor conjugate to $L_{2567}$. For each standard parabolic $P = QL$, we need to consider $X \hookrightarrow L$ via $(1^{[1]},3)$.

Let $P = QL$ be a standard parabolic subgroup whose Levi factor is $G$-conjugate to $L_{1567}$. Then for each of the $9$ choices of $P$, we get $\mathbb{V} \cong K^{2}$. In each case, one of two scenarios occurs:

{\bf Case 1)} The there exists a root element in $Q$, centralised by $L'$, inducing a non-zero isomorphism between the two modules $3 \otimes 1^{[1]}$ in the filtration of $Q$, and there also exists an element of $N_{G}(T)$ which fixes the root subgroups in $L'$ and swaps the root subgroups giving rise to these two $X$-modules. By Corollary \ref{cor:trivs} and Lemma \ref{lem:torus}, in this case, up to $G$-conjugacy there exists at most one non-$G$-cr complement to $Q$ in $QX$, corresponding to $(0,1) \in \mathbb{V}$.

{\bf Case 2)} There exists an $L'$-module $(0,\lambda_2)$ in some level $V_{j} = Q(j)/Q(j+1)$ of $Q$, which restricts to $X$ as $\bigwedge^{2}(3) = 4 + 0$. In this case, the $1$-dimensional trivial submodule lifts to an $X$-invariant subgroup $\{ x_{\alpha}(c)x_{\beta}(kc) \, : \, c \in K \} \subseteq Q^{X}$ for some pair of roots $\alpha$, $\beta$ and some fixed $k \in K$; again we get a non-zero homomorphism of $X$-modules between the two modules $3 \otimes 1^{[1]}$ occurring in the filtration of $Q$. Furthermore there exists an element of $N_{G}(T)$ stabilising the root subgroups in $L$ and swapping the appropriate root subgroups in $Q$, and so by Corollary \ref{cor:trivs} and Lemma \ref{lem:torus}, up to $G$-conjugacy there exists at most one non-$G$-cr complement to $Q$ in $QX$, corresponding to $(0,1) \in \mathbb{V}$.

Finally, for each of these nine parabolics, there exists an element of $N_{G}(T)$ sending root subgroups of the standard Levi factor to those of $L_{1567}$, and also sending the root subgroups generating one of the modules $3 \otimes 1^{[1]}$ to the appropriate root subgroups in $Q_{1567}$. It follows that all non-$G$-cr subgroups $A_{1}$ of these nine $A_{1}A_{3}$-parabolic subgroups of $G$, with irreducible image in the Levi factor, are $G$-conjugate.

Next consider the two remaining standard parabolics, $P_{2567}$ and $P_{2457}$, whose Levi factors are $G$-conjugate. Then for each, we get $\mathbb{V} \cong K^{4}$. In $P_{2567}$, we find that the subgroup $X < L$ given by the root elements
\begin{align*}
x_{+}(t) &= x_{\alpha_{2}}(t^5) x_{\alpha_{5}}(t) x_{\alpha_{6}}(2t) x_{\alpha_{7}}(3t) x_{0000110}(4t^2) x_{0000011}(2t^2) x_{0000111}(3t^3),\\
x_{-}(t) &= x_{-\alpha_{2}}(t^5) x_{-\alpha_{5}}(3t) x_{-\alpha_{6}}(2t) x_{-\alpha_{7}}(t) x_{-0000110}(3t^2)x_{-0000011}(t^2) x_{-0000111}(3t^3)
\end{align*}
commutes with the 9-dimensional unipotent subgroup of $Q$ given by
\begin{align*}
&\left\{
\begin{array}{c}
x_{0112210}(a)x_{0112111}(3a)x_{1112210}(b)x_{1112111}(3b)\\
\times \quad x_{1122210}(c)x_{1122111}(3c)
\end{array}
\, : \, a,b,c \in K\right\} \\
& \qquad \times \ U_{\alpha_{1}} U_{\alpha_{2}} U_{\alpha_{3}}U_{1224321}U_{1234321}U_{2234321},
\end{align*}
and this group induces non-zero $X$-module homomorphisms between each pair of modules in $Q$ with non-vanishing first cohomology group. By Proposition \ref{prop:trivs} and Lemma \ref{lem:torus}, every complement to $Q$ in $QX$ is $G$-conjugate to one corresponding to the vector $(k_1,k_2,k_3,k_4) \in \mathbb{V}$ where at most one coordinate is non-zero. Furthermore, for each pair of modules $3 \otimes 1^{[1]}$ occurring in the filtration of $Q$, there exists an element of $N_{G}(T)$ sending the root subgroups giving rise to one module to those giving rise to the other. It follows that up to $G$-conjugacy there is at most one non-$G$-cr subgroup $A_1$ in $P_{2567}$ with irreducible image in the Levi factor. Entirely similar calculations hold for the parabolic subgroup $P_{2457}$, and furthermore the element $n_{3} n_{4} n_{2} n_{5} n_{4} n_{3} n_{6} n_{5} n_{4} n_{2} n_{7} n_{6} n_{5} n_{4} n_{3} \in N_{G}(T)$ sends the root subgroups in $L_{2457}$ to those in $L_{2567}$, and also sends $U_{\alpha_6}$ to $U_{\alpha_4}$. Since these each give rise to an $X$-module $3 \otimes 1^{[1]}$ in the filtration of the corresponding unipotent radical, it follows that the classes of non-$G$-cr subgroups $A_1$ arising in each of these parabolic subgroups are $G$-conjugate.

Let $Y \cong A_{1} < A_{5}'$ via $W(5)$, and let $Z \cong A_{1} < A_{5}$ via $W(5)$. Then $Y$ and $Z$ are non-$G$-cr by Lemma \ref{lem:BMR}, and considering their composition factors on $L(G)$ and $V_{56}$ tells us that that $Y$ and $Z$ are not conjugate in $G$, and also that $Y$ lies in a parabolic subgroup whose Levi factor is $G$-conjugate to $L_{1567}$, and $Z$ lies in a parabolic subgroup whose Levi factor is $G$-conjugate to $L_{2567}$. Hence each of the two possible classes of non-$G$-cr subgroups above exists.

From Table \ref{E7p5tab}, $\text{dim}(C_{L(G)}(Y)) = 6$ and so $\text{dim}(C_G(Y)^\circ) \leq 6$. Also, $C_G(Y)^{\circ} \ge C_G(A_{5}')^{\circ} = A_{1} T_{1}$. Moreover in $P_{1567}$, the root subgroups giving rise to an $X$-module $3 \otimes 1^{[1]}$ in $Q(3)/Q(4)$ each commute with the elements $x_{\alpha_{2}}(c)$, $x_{1224321}(c)$ and $x_{1122210}(c)x_{1122111}(3c)$, for all $c \in K$, which generate the $3$-dimensional unipotent subgroup $Q^{X}$. Since $A_{1} < C_G(A_5')$ does not contain a $2$-dimensional unipotent subgroup, it follows that $C_G(Y)^{\circ} = U_{2} A_{1} T_{1}$.

Similarly, from Table \ref{E7p5tab}, $\text{dim}(C_{L(G)}(Z)) = 14$ and so $\text{dim}(C_G(Z)^\circ) \leq 14$. Also $C_G(Z)^{\circ} \ge C_G(A_5)^{\circ} = A_{2}$. Above, we have found a 9-dimensional subgroup of $Q_{2567}$ centralised by $X < L_{2567}$. On the other hand, let $l$ be maximal such that $Q_{2567}(l)$ contains a module $3 \otimes 1^{[1]}$ (in fact $l = 6$ here). Then any complement to $Q_{2567}(l)$ in $Q_{2567}(l)X$ must commute with all of $Q_{2567}^{X}$, since each element of this induces a homomorphism $Q_{2567}(l) \to Q_{2567}(l+j)$ for some $j$, which must be the zero map by maximality of $l$. Hence a non-$G$-cr complement to $Q_{2567}$ in $Q_{2567}X$, which we have shown is conjugate to a subgroup of $Q_{2567}(l)X$ by an element of $N_{G}(T)$, centralises a 9-dimensional unipotent subgroup generated by positive root elements. Since $A_{2}$ does not contain a $4$-dimensional unipotent subgroup, it follows that $C_G(Z)^{\circ} = U_{6} A_{2}$.

\subsection{$L'=A_1^2 A_2$}

There are 12 standard $A_1^{2} A_2$-parabolic subgroups of $G$. Let $X, Y \cong A_1 \hookrightarrow A_1^2 A_2$ via $(1,1^{[1]},2)$ and $(1^{[1]},1,2)$ respectively. In each such standard parabolic subgroup $P = QL$, there is precisely one level of $Q$ containing an $A_1^2 A_2$-module of the form $(1,1,10)$ or $(1,1,01)$, and it follows that $\mathbb{V} \cong K$ in each case. In the standard parabolic subgroup $P_{1235} = Q_{1235} L_{1235}$, the module $(1,1,01)$ is generated as an $X$-module by the image of the root group $U_{\alpha_4}$. Then the element $n_{0112221} n_{2234321} n_{1} n_{3} n_{4} n_{2} n_{5} n_{4} n_{3} n_{1} \in N_{G}(T)$ fixes $U_{\alpha_4}$ and swaps the $A_{1}$ factors of $L_{1235}$, hence each complement to $Q_{1235}$ in $Q_{1235}X$ is $G$-conjugate to a subgroup of $Q_{1235}Y$.

For $P$ equal to each of the other 11 standard parabolic subgroups, the following element of the Weyl group sends the roots of the standard Levi subgroup of $P$ to the roots in $L_{1235}$ and also sends a positive root $\alpha$, whose root group generates the $X$-module $(1,1,01)$ or $(1,1,10)$, to $\alpha_{4}$. Thus a lift of this element conjugates non-$G$-cr subgroups of these standard parabolics to subgroups of $P_{1235}$.

\begin{longtable}{c|c|c}
$P$ & Root $\alpha$ & Element of $N_{G}(T)$ \\ \hline
$P_{1457}$ & $0111110$ & $n_{6} n_{5} n_{4} n_{2} n_{3} n_{1} n_{4} n_{3} n_{7} n_{6} n_{5} n_{4} n_{2}$\\
$P_{2367}$ & $0001100$ & $n_{6} n_{5} n_{4} n_{2} n_{3} n_{1} n_{4} n_{3} n_{5} n_{4} n_{2} n_{6} n_{5} n_{4} n_{3} n_{1} n_{7} n_{6} n_{5} $ \\
$P_{1357}$ & $0112110$ & $n_{0112221} n_{1011111} n_{1111110} n_{1} n_{2} n_{3}$ \\
$P_{1256}$ & $0011000$ & $n_{6} n_{5} n_{4} n_{2} n_{3} n_{1} n_{4} n_{3} n_{5} n_{4} n_{2} n_{6} n_{5} n_{4} $ \\
$P_{1236}$ & $0001100$ & $n_6 n_5$ \\
$P_{1467}$ & $0111100$ & $ n_{0112221} n_{1011110} n_{0101110} n_{2} n_{3} n_{4} n_{5} $ \\
$P_{2356}$ & $\alpha_4$ & $n_{6} n_{5} n_{4} n_{2} n_{3} n_{1} n_{4} n_{3} n_{5} n_{4} n_{2} n_{6} n_{5} n_{4} n_{3} n_{1} $ \\
$P_{1267}$ & $0011100$ & $n_{6} n_{5} n_{4} n_{2} n_{3} n_{1} n_{4} n_{3} n_{5} n_{4} n_{2} n_{6} n_{5} n_{4} n_{7} n_{6} n_{5} $ \\
$P_{1237}$ & $0001110$ & $n_{6} n_{5} n_{7} n_{6}$ \\
$P_{1247}$ & $0011110$ & $n_{6} n_{5} n_{4} n_{2} n_{3} n_{1} n_{4} n_{3} n_{7} n_{6} $ \\
$P_{1246}$ & $0011100$ & $n_{6} n_{5} n_{4} n_{2} n_{3} n_{1} n_{4} n_{3} $

\end{longtable}
And therefore there is at most one $G$-conjugacy class of non-$G$-cr subgroups $A_{1}$ in an $A_{1}^{2}A_{2}$-parabolic subgroup with irreducible image in a Levi factor.

Let $Z \cong A_{1} < D_5$ via $T(6)$. Then $Z$ lies in an $A_{1}^{2} A_2$-parabolic subgroup of the Levi subgroup $D_{5}$ and is non-$D_5$-cr. Hence $Z$ lies in an $A_{1}^{2} A_2$-parabolic subgroup of $G$, and by Lemma \ref{lem:BMR} it is non-$G$-cr. From Table \ref{E7p5tab}, $\text{dim}(C_{L(G)}(Z)) = 4$, and so $C_G(Z)^\circ = C_G(D_5)^\circ = A_1 T_1$.

\subsection{$L' = D_5$ (Previously missing case)} \label{sec:E7D5again}

The two standard $D_5$-parabolic subgroups of $G$ are $P_{12345}$ and $P_{23456}$. Let $X \cong A_1 < D_5$ via $2 \otimes 2^{[1]} + 0$. For each standard parabolic subgroup $P = QL$, the unipotent radical $Q$ has three levels, of which levels $1$ and $2$ contain a summand $3 \otimes 1^{[1]}$ which have non-trivial cohomology. Arguing precisely as in Section~\ref{E7p5D5}, with $X$ having exactly the same trivial modules in the levels of $X$, we conclude that $X$ gives rise to at most one non-$G$-cr subgroup in each parabolic subgroup $P$. Moreover, again as argued in Section~\ref{E7p5D5}, the non-$G$-cr subgroups in each class of parabolics are in fact $G$-conjugate to one another. The non-$E_6$-cr subgroup $Y$ (or $Z$) exhibited in Section~\ref{sec:E6D5again} remains non-$G$-cr by Proposition~\ref{lem:BMR}, giving a representative of this class.

Inspecting Table~\ref{E6p5tab}, we see that $\dim C_{L(G)}(Y) = 3$. Once again arguing as Section~\ref{E7p5D5}, we conclude that $C_{G}(Y)^{\circ} = U_2 T_1$.


\section{Proof of Theorem \ref{thm:e7}: $G=E_7$, $p=7$} \label{sec:e7p7}

In this section we prove Theorem \ref{thm:e7} in the case $p=7$, in which case both $A_1$ and $G_2$ subgroups occur. The starting point is the following lemma. 

\begin{lemma} \label{lem:e7p7badlevi}
Let $L$ be a Levi subgroup of $G$ containing an $L$-irreducible subgroup $X$ of type $A_1$ or $G_2$. If there exists a parabolic subgroup $P$ of $G$ with Levi factor $L$ and unipotent radical $Q$, such that $H^{1}(X,M \downarrow X) \neq 0$ for some level $M$ of $Q$, then $X$ and the type of $L'$ appear in Table \ref{tab:e7p7badX}.
\end{lemma}

\begin{longtable}{>{\raggedright\arraybackslash}p{0.15\textwidth - 2\tabcolsep}>{\raggedright\arraybackslash}p{0.47\textwidth-\tabcolsep}@{}}
\caption{$L'$-irreducible $X$ with $H^1(X,M\downarrow X) \neq 0$\label{tab:e7p7badX}} \\ 
$L'$ & Embedding of $X$ \\ \hline
$E_6$ & $X \cong G_2$, see Proposition \ref{prop:irredG2} \\
& $X \cong A_1 \hookrightarrow A_1 A_5$ via $(1^{[1]},5)$ \\
$A_1 A_2 A_3$ & $X \cong A_1 \hookrightarrow A_1 A_2 A_3$ via $(1^{[1]},2,3)$ \\ \hline
\end{longtable}

\begin{proof}
This follows in the same manner as Lemma \ref{lem:e6badlevi}, where the classes of subgroups $X \cong G_2$ follow from Proposition \ref{prop:irredG2}. We note that $X \cong G_2 < A_{6}$ via $10$ has a composition factor $20$ on $V_{A_6}(\lambda_3)$. This module occurs in the filtration of the unipotent radical of an $A_{6}$-parabolic and has non-vanishing first cohomology group by Lemma \ref{lem:G2mods}\ref{G2mod20}. However, $V_{A_6}(\lambda_3) =\bigwedge^{3} (V_{A_6}(\lambda_1))$, hence $V_{A_6}(\lambda_3) \downarrow X = \bigwedge^{3} (V_{G_2}(10))$ is tilting, by Lemma \ref{lem:tiltingnoH1}(iii). By Lemma \ref{lem:tiltingnoH1}(iv), it follows that $H^1(X, V_{A_6}(\lambda_3) \downarrow X) = 0$.
\end{proof}

\subsection{$L'=E_6$} \label{sec:E7E6p7}
The unique standard $E_6$-parabolic subgroup of $G$ is $P = P_{123456} = QL$. In this case $Q$ is an irreducible $L'$-module of high weight $\lambda_1$. First, let $X \cong A_1 \hookrightarrow A_1 A_5'$ via $(1^{[1]},5)$. Since $V_{E_6}(\lambda_1) \downarrow A_1 A_5' = (1,\lambda_1) + (0,\lambda_4)$, it follows that $V_{E_6}(\lambda_1) \downarrow X = (1^{[1]} \otimes 5) + T(8) + 0$. Therefore, $H^1(X,Q) \cong K$ and by Lemma \ref{lem:torus} we find that there is exactly one $G$-conjugacy class of non-$G$-cr complements to $Q$ in $QX$. Let $A \cong A_1 < A_7$ via $W(7)$. Then $A$ is non-$G$-cr by Lemma \ref{lem:BMR} and is contained in an $A_1 A_5$-parabolic subgroup of $A_7$, and therefore in an $E_{6}$-parabolic subgroup of $G$. Thus $Y$ is a representative of this class of non-$G$-cr subgroups of $P$.

Still with $P = QL$ as above, we now let $Y$ be the $L'$-irreducible subgroup of $L'$ of type $G_{2}$. Then $Y$ is contained in a subgroup $F_4$ of $L'$, and hence $V_{E_{6}}(\lambda_1) \downarrow Y = 20 + 00$. By Lemma \ref{lem:G2mods}, $H^1(Y,Q) \cong K$ and thus by Lemma \ref{lem:torus} we have exactly one $G$-conjugacy class of non-$G$-cr complements to $Q$ in $QY$. We let $B$ be a representative of such a non-$G$-cr complement and now prove that $B$ is not properly contained in any proper connected reductive subgroup of $G$. Suppose that $H$ is maximal among proper reductive subgroups of $G$ containing $B$. Using Lemma \ref{lem:maximalexcep}, the restriction $V_{56} \downarrow B = T(20)^2$ as given in Table \ref{E7p7g2tab}, and the restrictions in \cite[Tables 8.2, 8.6]{MR1329942}, we see that $H$ must be simple of type $A_{7}$. But this implies that $B$ stabilises a 1-space on $V_{A_7}(\lambda_1)$, since the only non-trivial irreducible $G_2$-modules of dimension at most $8$ are Frobenius twists of the $7$-dimensional module $10$. Thus $B$ lies in an $A_6$-parabolic subgroup of $H$; but $V_{56} \downarrow A_{6} = \lambda_1/\lambda_2/\lambda_5/\lambda_6$ by \cite[Table 8.6]{MR1329942}, which implies that $V_{56} \downarrow B$ has a 7-dimensional section as a $B$-module; a contradiction. Thus $B$ does not lie in any proper reductive subgroup of $G$.

Now $\text{dim}(C_{L(G)}(A)) = \textup{dim}(C_{L(G)}(B)) = 1$; since $Q$ is abelian and contains a $1$-dimensional trivial $L'$-submodule, this submodule is a $1$-dimensional unipotent group centralised by $A$ and by $B$. Thus $C_G(A)^{\circ} = C_{G}(B)^{\circ} = U_1$.

For use later in Section \ref{sec:restrictions} when computing the restriction of $G$-modules $V_{56}$ and $L(G)$ to $B$, we now show that a conjugate of $A$ is in fact a subgroup of $B$. We let $Z$ be a maximal subgroup $A_{1}$ of $Y$, so that $V_{G_{2}}(10) \downarrow Z = 6$. Since $V_{G_2}(10)$ is tilting and $p > 2$, by considering weight multiplicities it follows that the symmetric square $S^{2}(V_{G_2}(10)) = T(20)$. Since the $Z$-module $6$ is tilting, it follows that $S^{2}(6) = T(12) + T(8)$, and therefore $Q \downarrow Z = 12 + T(8) + 0$. This shows two things: Firstly, $Z$ is $L'$-irreducible, since its composition factors on $Q = V_{E_{6}}(\lambda_1)$ are incompatible with the actions of every proper Levi subgroup of $L'$ (see \cite[Table 8.7]{MR1329942}), and furthermore $H^{1}(Z,Q) \cong K$. Thus $Z$ is $L'$-conjugate to the subgroup $X$ above, so without loss of generality we assume $Z = X$. Secondly, there exists a non-trivial extension of $Y$-modules $V = 00|20$ such that $V \downarrow Z = (0|12) + T(8)$ does not contain a trivial submodule. Since $H^{1}(Y,20)$ and $H^{1}(Z,12)$ are each $1$-dimensional, the restriction map $H^{1}(Y,20) \to H^{1}(Z,12)$ is an isomorphism of vector spaces. In particular each non-$G$-cr complement to $X$ in $QX$ lies in a non-$G$-cr complement to $Z$ in $QZ$ and it follows that a conjugate of $A$ is contained in $B$. 

\subsection{$L'=A_1 A_2 A_3$} \label{sec:E7A1A2A3p7}

Let $P = P_{123567} = QL$ by the unique standard $A_1 A_2 A_3$-parabolic subgroup of $G$. Let $X \cong A_1 \hookrightarrow A_1 A_2 A_3$ via $(1^{[1]},2,3)$. Using Lemma \ref{lem:h1fora1} as in previous calculations, we find that $H^{1}(X,Q/Q(2)) \cong K$, while the corresponding cohomology group for the other levels of $Q$ vanishes. By Lemma \ref{lem:torus}, there is at most one conjugacy class of non-$G$-cr subgroups in $QX$.

Let $Y \cong A_1 \hookrightarrow A_1 G_2$ via $(1,6)$ and $Z \cong A_1 \hookrightarrow G_2 C_3$ via $(6,5)$. Using the restrictions of $L(G)$ and $V_{G}(\lambda_7)$ to $A_{1}G_{2}$ and $G_{2}C_{3}$ as given in \cite[Table 10.1,\ 10.2]{MR2044850}, we find that the action of $Y$ and $Z$ on these modules is as given in Table \ref{E7p7a1tab}. In particular, $Y$ and $Z$ each fix a non-zero element of $V_{G}(\lambda_7)$, and therefore lie in a proper subgroup of dimension at least $133 - 56 = 77$; all such subgroups are contained in parabolics. Furthermore $\text{dim}(C_{L(G)}(Y)) = 0$, and so $Y$ and $Z$ cannot centralise a non-trivial torus of $G$, and thus do not lie in a proper Levi subgroup. Hence $Y$ and $Z$ are non-$G$-cr. Since the action of $Y$ and $Z$ on $L(G)$ does not agree with a non-$G$-cr subgroup $A_{1}$ in a parabolic subgroup of type $E_{6}$, we deduce that $Y$ and $Z$ lie in a parabolic subgroup of type $A_{1}A_{2}A_{3}$ and hence both are representatives of the unique class above.

Since $\text{dim}(C_{L(G)}(Y)) = 0$ we deduce that $C_G(Y)^\circ = 1$.


\section{Proof of Theorem \ref{thm:e8}: $G=E_8$, $p=7$} \label{sec:e8}

In this section we prove Theorem \ref{thm:e8}. As in the previous sections, the starting point is the following lemma, which determines parabolics $P = QL$ of $G$ and $L$-irreducible subgroups $X$ such that $H^{1}(X,Q)$ may be non-zero. The proof is identical to that of Lemma \ref{lem:e7p7badlevi}. 

\begin{lemma} \label{lem:e8badlevi}
Let $L$ be a Levi subgroup of $G$ containing an $L$-irreducible subgroup $X$ of type $A_1$ or $G_2$. If there exists a parabolic subgroup $P$ of $G$ with Levi factor $L$ and unipotent radical $Q$, such that $H^{1}(X,M \downarrow X) \neq 0$ for some level $M$ of $Q$, then $X$ and the type of $L'$ appear in Table \ref{tab:e8p7badX}.
\end{lemma}

\begin{longtable}{>{\raggedright\arraybackslash}p{0.13\textwidth - 2\tabcolsep}>{\raggedright\arraybackslash}p{0.87\textwidth-\tabcolsep}@{}}
\caption{$L'$-irreducible $X$ with $H^1(X,M\downarrow X) \neq 0$ \label{tab:e8p7badX}} \\
$L'$ & Embedding of $X$ \\ \hline
$D_7$ & $X \cong G_2 < D_7$ via $01$ (two $L'$-conjugacy classes) \\
$A_1 E_6$ & $X \cong A_1 \hookrightarrow A_1 A_1 A_5 < A_1 E_6$ via $(1^{[r]},1^{[s+1]},1^{[s]})$ $(rs=0)$ \\
$A_2 D_5$ & $X \cong A_1 \hookrightarrow A_2 D_5$ via $(2^{[r]},4^{[r]} + 1^{[r+1]} \otimes 1^{[s]})$ $(rs=0;$ $r+1 \neq s)$ \\
$A_3 A_4$ & $X \cong A_1 \hookrightarrow A_3 A_4$ via $(1 \otimes 1^{[1]},4)$ \\
$E_6$ & $X \cong A_1 \hookrightarrow A_1 A_5$ via $(1^{[1]},5)$ \\
& $X \cong G_2$, see Proposition \ref{prop:irredG2} \\
$D_6$ & $X \cong A_1 < D_6$ via $1^{[1]} \otimes 5$ (two $L'$-conjugacy classes) \\
$A_1 A_5$ & $X \cong A_1 \hookrightarrow A_1 A_5$ via $(1^{[1]},5)$ \\
$A_2 D_4$ & $X \cong A_1 \hookrightarrow A_2 D_4$ via $(2, 4+2^{[1]})$ \\
& $X \cong A_1 \hookrightarrow A_2 D_4$ via $(2, 3\otimes 1^{[1]})$  (two $L'$-conjugacy classes) \\
$A_1^2 A_4$ & $X \cong A_1 \hookrightarrow A_1^2 A_4$ via $(1,1^{[1]},4)$ \\
& $X \cong A_1 \hookrightarrow A_1^2 A_4$ via $(1^{[1]},1,4)$ \\
$A_1 A_2 A_3$ & $X \cong A_1 \hookrightarrow A_1 A_2 A_3$ via $(1^{[1]},2,3)$ \\ \hline 
\end{longtable}

\subsection{$L'=D_7$} \label{sec:E8D7}
Let $P = P_{2345678} = QL$ be the unique standard $D_7$-parabolic subgroup of $G$. Let $X \cong G_2 < L'$ with $V_{D_7}(\lambda_1) \downarrow X = 01$. Then there are two $L$-conjugacy classes of such subgroups in $D_7$, which are distinguished by their action on $V_{D_7}(\lambda_7)$. Indeed, as outlined in Section \ref{sec:restrictions}, if $X$ and $Y$ are representatives of these subgroup classes then we can verify computationally that $V_{D_7}(\lambda_7) \downarrow X$ and $V_{D_{7}}(\lambda_7) \downarrow Y$ are uniserial with two composition factors, of dimension $26$ and $38$, and we can therefore pick $V_{D_7}(\lambda_7) \downarrow X = 11 | 20 $ and $V_{D_7}(\lambda_7) \downarrow Y = 20 | 11$.

Now $Q$ has two levels, and as $L'$-modules we have $Q / Q(2) = \lambda_7$ and $Q(2) = \lambda_1$. By Lemma \ref{lem:G2mods}(i) and Lemma \ref{lem:tiltingnoH1}(iv) we have $H^{1}(X,Q(2)) = H^{1}(Y,Q(2)) = 0$. By Lemma \ref{lem:G2mods}(iii) and (iv) we have $W(20) = 20 | 00$ and $W(11) = 11|20$. This implies that $H^{1}(G_2,20) \cong K$ and $H^{1}(G_2,11) = H^{0}(G_2,11) = 0$. From the long exact sequence of cohomology induced from $20 \hookrightarrow 11|20 \twoheadrightarrow 11$, we deduce that $H^{1}(G_2,11|20) \cong H^{1}(G_2,20) \cong K$. Now assume that $H^{1}(G_2,20|11) \neq 0$. If $V$ is a corresponding indecomposable extension of $20|11$ by the trivial module, then $V^{\ast}$ has shape $11|(20/00)$. Since all high weights here are less than $11$, this module is an image of $W(11)$, which is absurd. Therefore $H^{1}(G_2,20|11) = 0$. We have just shown that $H^1(X,Q/Q(2)) \cong K$ and $H^1(Y,Q/Q(2)) = 0$. Hence $\mathbb{V}_{X,Q} \cong K$ and $\mathbb{V}_{Y,Q} = 0$, and by Lemma \ref{lem:torus} there exists at most one $G$-conjugacy class of non-$G$-cr complements to $Q$ in $QX$, and none in $QY$. 

Consider $Z = G_2 \hookrightarrow G_2 G_2 < G_2 F_4$ via $(10,10)$, where the second factor $G_2$ is maximal in $F_4$. By \cite[Lemma 7.13]{Tho1}, $Z$ is contained in a $D_7$-parabolic subgroup of $G$. Furthermore, from Table \ref{E8p7a1tab} we see that $L(G) \downarrow Z$ has no non-zero trivial submodules and hence $Z$ is not contained in any Levi subgroup of $G$. Therefore, $Z$ is a representative of the conjugacy class of non-$G$-cr subgroups in $QX$. Since $\text{dim}(C_{L(G)}(Z)) = 0$, it follows that $C_G(Z)^\circ = 1$.

\subsection{$L'=A_1 E_6$}
Let $P = P_{1234568} = QL$ be the unique standard $A_1 E_6$-parabolic subgroup of $G$. Let $X \cong A_1 \hookrightarrow A_1 A_1 A_5 < A_1 E_6$ via $(1^{[r]},1^{[s+1]},1^{[s]})$ $(rs=0)$. Then $Q$ has three levels, and $H^{1}(X,Q(2)/Q(3)) \cong K$, while the corresponding group for the other levels vanishes. Applying Lemma \ref{lem:torus}, there is at most one $G$-conjugacy class of non-$G$-cr subgroups in $QX$. 

Now consider $Y \cong A_1 \hookrightarrow A_1 A_7$ via $(1^{[r]},W(7)^{[s]})$, where $A_7$ is a maximal connected subgroup of an $E_7$ Levi subgroup of $G$. Then $Y$ is non-$G$-cr by Lemma \ref{lem:BMR} and $Y$ is contained in an $A_1 E_6$-parabolic subgroup of $A_1 E_7$ and hence of $G$. Therefore $Y$ is a representative of the conjugacy class of non-$G$-cr complements to $Q$ in $QX$.

From Table \ref{E8p7a1tab}, $\text{dim}(C_{L(G)}(Y)) = 1$. Since the image of $Y$ in $A_7$ is non-$E_7$-cr, by Table \ref{E7p7a1tab} this image centralises a 1-dimensional unipotent subgroup of $E_7$, and it follows that $C_G(Y)^\circ = U_1$.

\subsection{$L'=A_2 D_5$} \label{sec:E8p7A2D5}
Let $P = P_{1234578} = QL$ be the unique standard $A_2 D_5$-parabolic subgroup of $G$. Let $X \cong A_1 \hookrightarrow A_2 D_5$ via $(2^{[r]},4^{[r]} + 1^{[r+1]} \otimes 1^{[s]} + 0)$ $(rs=0$; $r+1 \neq s)$. The action of $X$ on the levels of $Q$ is as follows:
\begin{align*}
Q/Q(2) \downarrow X &= 5^{[r]} \otimes 1^{[r+1]} + 3^{[r]} \otimes 1^{[r+1]} + 1^{[r]} \otimes 1^{[r+1]} + 5^{[r]} \otimes 1^{[s]} + 3^{[r]} \otimes 1^{[s]} + 1^{[r]} \otimes 1^{[s]},\\
Q(2)/Q(3) \downarrow X &= 6^{[r]} + 4^{[r]} + 2^{[r]} \otimes 1^{[r+1]} \otimes 1^{[s]} + (2^{[r]})^2, \\
Q(3)/Q(4) \downarrow X &= 3^{[r]} \otimes 1^{[r+1]} + 3^{[r]} \otimes 1^{[s]}, \\
Q(4)/Q(5) \downarrow X &= 2^{[r]}.
\end{align*}
We see that $\mathbb{V}_{X,Q} \cong K$, the unique module with non-vanishing first cohomology group occurring in level $1$. Applying Lemma \ref{lem:torus} there is at most one non-$G$-cr complement to $Q$ in $QX$, for each $r$ and $s$, up to $G$-conjugacy.

Let $Y \cong A_1 \hookrightarrow A_1 A_1 G_2 < A_1 E_7$ via $(1^{[s]},1^{[r]},6^{[r]})$ ($rs=0$; $r \neq s+1$). From Section \ref{sec:E7A1A2A3p7} we know that the image of $Y$ in $E_7$ is non-$E_7$-cr. By \cite[Lemma 2.12]{MR2178661}, it follows that $Y$ is non-$A_1 E_7$-cr, and since $A_1 E_7$ is a subsystem subgroup of $G$, by Lemma \ref{lem:BMR} we conclude that $Y$ is non-$G$-cr. The image of $Y$ in $E_7$ lies in an $A_1 A_2 A_3$-parabolic subgroup of $E_7$, and so $Y$ lies in a parabolic subgroup of $G$ whose Levi factor contains a subgroup of type $A_1^{2} A_2 A_3$. The only such Levi subgroup is $A_{2} D_{5}$, hence $Y$ is a representative of the class of non-$G$-cr subgroups above.

From Table \ref{E8p7a1tab}, we have $\text{dim}(C_{L(G)}(Y)) = 0$ if $r \neq s$ and $\text{dim}(C_{L(G)}(Y)) = 1$ if $r=s$. It follows that $C_G(Y)^\circ$ is trivial if $r \neq s$. If $r = s = 0$ there is a module $1 \otimes 1 = 2 + 0$ occurring in $Q/Q(2)$. Identifying the root elements of $X$ as in previous calculations, we find that for $c \in K^*$, the following element generates a $1$-dimensional subgroup of $Q$ which commutes with $X$ and with the roots giving the module $5 \otimes 1^{[1]}$ in $Q/Q(2)$:
\begin{align*}
x_{00011100}&(c)x_{00001110}(2c)x_{00000111}(c)x_{12232100}(6c)x_{11232110}(5c)x_{11222111}(6c)x_{12233210}(6c^2)\\ & \times x_{12232211}(3c^2)x_{11233211}(3c^2)x_{11232221}(c^2)x_{12233321}(4c^3)x_{23464321}(6c^3)x_{23465431}(2c^4)
\end{align*}
and it follows that $C_{G}(Y)^{\circ} = U_{1}$.

\subsection{$L'=A_3 A_4$} \label{sec:E8A3A4}
Let $P = P_{1234678} = QL$ be the unique standard $A_3 A_4$-parabolic subgroup of $G$. We need to consider $X \cong A_1 \hookrightarrow A_3 A_4$ via $(1 \otimes 1^{[1]},4)$. With similar calculations to previous sections, we find that $H^{1}(X,Q(2)/Q(3)) \cong K$, while the corresponding cohomology groups for other levels vanishes. Applying Lemma \ref{lem:torus}, there is at most one class of non-$G$-cr complements to $Q$ in $QX$.

Let $Y = A_1 < A_8$ via $W(8)$. Then $Y$ is non-$A_8$-cr and thus non-$G$-cr by Lemma \ref{lem:BMR}, and comparing the composition factors of $Y$ on $L(G)$ with those of Levi subgroups of $G$ in \cite[Table 8.1]{MR1329942} shows that $Y$ can only lie in a parabolic subgroup of $G$ with Levi factor $A_{3} A_{4}$. Thus a conjugate of $Y$ is a representative of the class of non-$G$-cr complements to $Q$ in $QX$

From Table \ref{E8p7a1tab}, we have $\text{dim}(C_{L(G)}(Y)) = 0$ and hence $C_G(Y)^\circ = 1$. 

\subsection{$L'=E_6$} \label{sec:E8E6}

Let $P = P_{123456} = QL$ be the unique standard $E_{6}$-parabolic subgroup of $G$. Let $X \cong A_1 \hookrightarrow A_1 A_5 < E_6$ via $(1^{[1]},5)$. Then the actions of $X$ on the levels of $Q$ are as follows: 
\begin{align*}
Q/Q(2) \downarrow X &= 5 \otimes 1^{[1]} + T(8) + 0^2,\\
Q(2)/Q(3) \downarrow X &= 5 \otimes 1^{[1]} + T(8) + 0, \\
Q(3)/Q(4) \downarrow X &= 5 \otimes 1^{[1]} + T(8) + 0, \\
Q(4)/Q(5) \downarrow X &= 0, \\
Q(5)/Q(6) \downarrow X &= 0.
\end{align*}

From Lemma \ref{lem:h1fora1} it follows that $\mathbb{V}_{X,Q} \cong K^{3}$. As in previous sections, it is straightforward to determine root subgroups of $X$ in terms of root subgroups of $A_{1} A_{5}$, which are root subgroups of $G$. The trivial modules in levels $1$ to $5$ each lift to elements of $Q^{X}$ as follows:
\begin{center}
\begin{tabular}{c|c}
Level & Elements \\ \hline
$1$ & $x_{\alpha_{8}}(a)x_{11221110}(b)x_{11122110}(5b)x_{01122210}(3b)$\\
$2$ & $x_{11221111}(a)x_{11122111}(5a)x_{01122211}(3a)$ \\
$3$ & $x_{22343221}(a)x_{12343321}(2a)x_{12244321}(3a)$ \\
$4$ & $x_{23465431}(a)$ \\
$5$ & $x_{23465432}(a)$
\end{tabular}
\end{center}
If $a \neq 0$, the root element $x_{\alpha_{8}}(a)$ induces a non-trivial homomorphism $Q/Q(2) \to Q(2)/Q(3)$ of $L$-modules, and centralises each root group for a root of level 2. Similarly, each element of the form $x_{11221110}(a) x_{11122111}(5a) x_{01122210}(3a)$ induces a non-trivial $X$-module homomorphism $Q(2)/Q(3) \to Q(3)/Q(4)$, and centralises the root groups in $Q(1)$ besides $U_{\alpha_{8}}$. Finally, a lift of a non-trivial element of $Q(2)/Q(3)$ gives rise to a non-trivial $X$-module homomorphism $Q/Q(2) \to Q(3)/Q(4)$. Applying Proposition \ref{prop:trivs}, if complements to $Q$ in $QX$ are parametrised by $(k_1,k_2,k_3) \in \mathbb{V}_{X,Q}$, we may assume $k_1 k_2 = k_2 k_3 = k_1 k_3 = 0$. Next, the element $n_{8} \in N_{G}(T)$ normalises each root subgroup in $L$ and swaps the root subgroups occurring in level $1$, other than $\alpha_8$, with those of level $2$. In addition, the element $n_{7} n_{6} n_{5} n_{4} n_{2} n_{3} n_{1} n_{4} n_{3} n_{5} n_{4} n_{2} n_{6} n_{5} n_{4} n_{3} n_{1} n_{7} n_{6} n_{5} n_{4} n_{2} n_{3} n_{4} n_{5} n_{6} n_{7}$ stabilises the root subgroups in $L$ and exchanges the root subgroups in levels $2$ and $3$. Together with Lemma \ref{lem:torus} and the subsequent discussion, we conclude that there is at most one $G$-conjugacy class of non-$G$-cr complements to $Q$ in $QX$.

Recall that $A_7'$ denotes a subgroup $A_7$ of $G$ which lies in a Levi subgroup $E_7$. Let $A \cong A_1 < A_7'$ via $W(7)$. Then from Section \ref{sec:E7E6p7} we know that $A$ is non-$E_7$-cr and contained in an $E_6$-parabolic subgroup of $E_7$, with irreducible image in the Levi factor. Hence $A$ is non-$G$-cr by Lemma \ref{lem:BMR} and $A$ is a representative of the class of non-$G$-cr complements to $Q$ in $QX$.

Now let $Y \cong G_2$ be an $E_6$-irreducible subgroup (see Proposition \ref{prop:irredG2}). Then the actions of $Y$ on the levels of $Q$ are as follows: 
\begin{align*}
Q/Q(2) \downarrow Y &= 20 + 00^2,\\
Q(2)/Q(3) \downarrow Y &= 20 + 00, \\
Q(3)/Q(4) \downarrow Y &= 20 + 00, \\
Q(4)/Q(5) \downarrow Y &= 00, \\
Q(5)/Q(6) \downarrow Y &= 00.
\end{align*}

By Lemma \ref{lem:G2mods}, we have $H^1(G_2,20) \cong K$, hence complements are parametrised by $(k_1,k_2,k_3) \in \mathbb{V}_{Y,Q}$. An entirely similar argument to the above shows that there is at most one $G$-conjugacy class of non-$G$-cr complements to $Q$ in $QY$. Moreover, a representative subgroup may be taken to lie in $Q(3) Y$. 

Let $B$ be the non-$E_7$-cr subgroup $G_2$ of $E_7$ given by Theorem \ref{thm:e7}. Then $B$ is non-$G$-cr by Lemma \ref{lem:BMR} and is contained in an $E_6$-parabolic subgroup $E_{7}$ and hence of $G$. Hence $B$ is a representative of the conjugacy class of non-$G$-cr subgroups contained in $QY$.

From Table \ref{E8p7a1tab}, we have $\text{dim}(C_{L(G)}(A)) = \text{dim}(C_{L(G)}(B)) = 8$. We claim that $C_{G}(B)^\circ = U_5 \bar{A}_1$. Since $A$ is conjugate to a subgroup of $B$, as proved in Section \ref{sec:E7E6p7}, it follows from the claim that $C_{G}(A)^\circ = U_5 \bar{A}_1$. To prove the claim we start by noting that the elements of $Q^X$ in the above table are all centralised by $Y$ and so $Q^Y$ is a $6$-dimensional unipotent subgroup of $C_G(Y)$. We let $R$ be the subgroup of $Q$ generated by $Q^Y$ and $Q(3)$. We know that $B$ is conjugate to a non-$G$-cr subgroup of $R Y$, call it $Z$. It now follows that $C_{R Y}(Z)$ contains $Q^Y$. Moreover, since $\bar{A}_1 = C_{G}(E_7) < C_G(B)$, it follows that $C_{G}(Z)$ contains a subgroup of type $A_1$ generated by root subgroups of $G$. The intersection of this subgroup $A_1$ with $Q^Y < Q$ can be at most $1$-dimensional. Therefore, we have found a subgroup $U_5 \bar{A}_1 < C_{G}(Z)$. Since $\text{dim}(C_G(Z)) \leq 8$ we have proved that $C_{G}(Z)^\circ = U_5 \bar{A}_1$ and the claim immediately follows.

\subsection{$L'=D_6$}
Let $P = P_{234567} = QL$ be the unique standard $D_{6}$-parabolic subgroup of $G$. Let $X$ and $Y$ be representatives of the two conjugacy classes of $A_1$ subgroups in $D_6$ which act as $5 \otimes 1^{[1]}$ on $V_{D_6}(\lambda_1)$, with $V_{D_6}(\lambda_5) \downarrow X \cong V_{D_6}(\lambda_6) \downarrow Y = T(9) + 5 \otimes 2^{[1]}$ and $V_{D_6}(\lambda_6) \downarrow X \cong V_{D_6}(\lambda_5) \downarrow Y = T(8) \otimes 1^{[1]} + 3^{[1]}$. Then the actions of $X$ and $Y$ on the levels of $Q$ are as follows:

\begin{align*}
Q/Q(2) \downarrow X &= 5 \otimes 1^{[1]} + T(8) \otimes 1^{[1]} + 3^{[1]},\\
Q(2)/Q(3) \downarrow X &= T(9) + 5 \otimes 2^{[1]} + 0, \\
Q(3)/Q(4) \downarrow X &= 5 \otimes 1^{[1]}, \\
Q(4)/Q(5) \downarrow X &= 0,\\
\ \\
Q/Q(2) \downarrow Y &= 5 \otimes 1^{[1]} + T(9) + 5 \otimes 2^{[1]},\\
Q(2)/Q(3) \downarrow Y &= T(8) \otimes 1^{[1]} + 3^{[1]} + 0, \\
Q(3)/Q(4) \downarrow Y &= 5 \otimes 1^{[1]}, \\
Q(4)/Q(5) \downarrow Y &= 0.
\end{align*}
By Lemma \ref{lem:h1fora1} we have $\mathbb{V}_{X,Q} \cong \mathbb{V}_{Y,Q} \cong K^2$. For both $X$ and $Y$, the root groups $U_{\alpha_8}$ and $U_{22343211}$ give rise to the modules $5 \otimes 1^{[1]}$ in levels $1$ and $3$, and a non-trivial element of $U_{22343210}$ is fixed by $L'$ and induces a non-trivial $L'$-module homomorphism $Q/Q(2) \to Q(3)/Q(4)$. Thus if complements to $Q$ in $QX$ or $QY$ are parametrised by $(k_1,k_2) \in \mathbb{V}_{X,Q}$ or $\mathbb{V}_{Y,Q}$ respectively, by Corollary \ref{cor:trivs} we may assume that $k_1 k_2 = 0$. In addition, the element $n_{22343210} \in N_{G}(T)$ normalises the root subgroups in $L'$, and swaps $U_{\alpha_8}$ and $U_{22343211}$. Thus, applying Lemma \ref{lem:torus}, up to $G$-conjugacy there exists at most one non-$G$-cr complement to $Q$ in $QX$ and at most one non-$G$-cr complement to $Q$ in $QY$. Finally, the element $n_{8} n_{7} n_{6} n_{5} n_{4} n_{2} n_{3} n_{4} n_{5} n_{6} n_{7} n_{8}$ induces a graph automorphism on $L'$, swapping the class of $X$ and $Y$, while centralising the root group $U_{22343211}$. Thus a non-$G$-cr complement to $Q$ in $QX$ is conjugate to a subgroup of $QY$.

Let $Z \cong A_1 < D_7$ via $T(12)$. Then $Z$ is non-$G$-cr by Lemma \ref{lem:BMR} and is contained in a $D_6$-parabolic subgroup of the Levi subgroup $D_7$. Hence $Z$ lies in a $D_7$-parabolic subgroup of $G$, and is thus conjugate to a non-$G$-cr complement to $Q$ in $QX$.

From Table \ref{E8p7a1tab}, $\text{dim}(C_{L(G)}(Z)) = 3$, and $C_G(Z) \ge C_G(D_7) = T_1$. Now, the two root groups $U_{22343210}$ and $U_{23465432}$ centralise $L$ and $Q(2)$, thus a non-$G$-cr complement to $Q$ in $QX$ corresponding to $(0,1) \in \mathbb{V}_{X,Q}$, which lies in $Q(2)X$, also centralises this 2-dimensional unipotent subgroup generated by positive root elements. It follows that $C_G(Z)^{\circ} = U_{2} T_{1}$.

\subsection{$L'=A_1 A_5$}

The three standard $A_{1} A_{5}$-parabolic subgroups of $G$ are $P_{124567}$, $P_{145678}$ and $P_{134568}$. Let $X \cong A_1 \hookrightarrow A_1 A_5$ via $(1^{[1]},5)$. First consider $P_{124567}$. Then the actions of $X$ on the levels of $Q = Q_{124567}$ are as follows:

\begin{align*}
Q/Q(2) \downarrow X &= T(8) \otimes 1^{[1]} + 5 + 1^{[1]},\\
Q(2)/Q(3) \downarrow X &= T(8) + 5 \otimes 1^{[1]} + 0, \\
Q(3)/Q(4) \downarrow X &= T(9) + 5 + 1^{[1]}, \\
Q(4)/Q(5) \downarrow X &= 5 \otimes 1^{[1]}, \\
Q(5)/Q(6) \downarrow X & = 5, \\
Q(6)/Q(7) \downarrow X &= 0.
\end{align*} 

By Lemma \ref{lem:h1fora1} we have $\mathbb{V}_{X,Q} \cong K^2$. The modules of high weight $5 \otimes 1^{[1]}$ are generated by $U_{00111111}$ and $U_{12343211}$, and the element $n_{12343210} n_{1} n_{3} n_{4} n_{2} n_{5} n_{4} n_{3} n_{6} n_{5} n_{4} n_{2} n_{7} n_{6} n_{5} n_{4} n_{3} \in N_{G}(T)$ stabilises the set of root subgroups in $L'$ while swapping $U_{00111111}$ and $U_{12343211}$. Now, $X$ lies in a subgroup $A_{1}C_{3}$ of $L'$. This subgroup acts on $Q(2)/Q(3)$ as $(1,\lambda_1) + (0,\bigwedge^{2}\lambda_1)$. Since $\bigwedge^{2}\lambda_1$ is self-dual it follows that $\bigwedge^{2}\lambda_1 = \lambda_2 + 0$. This trivial direct summand is the trivial $X$-module summand in $Q(2)/Q(3)$. Using Lemma \ref{lem:tiltingnoH1}, the other $A_{1}C_{3}$-summands occurring in the filtration of $Q$ are tilting, and therefore have zero first cohomology group. Thus the trivial summand in $Q(2)/Q(3)$ lifts to a $1$-dimensional subgroup of $Q^{A_{1}C_{3}} \le Q^{X}$, and therefore gives rise to a non-trivial $X$-module homomorphism $Q(2)/Q(3) \to Q(4)/Q(5)$. By Corollary \ref{cor:trivs}, if complements to $Q$ in $QX$ are parametrised by $(k_1, k_2) \in \mathbb{V}_{X,Q}$, then we may assume $k_1 k_2 = 0$. Applying Lemma \ref{lem:torus} and the above Weyl group elements, there is at most one $G$-conjugacy class of non-$G$-cr complements to $Q$ in $QX$. 

Entirely similar arguments hold for $P_{134568}$ and $P_{145678}$. In the filtration of $Q_{134568}$, one of the two modules $5 \otimes 1^{[1]}$ is generated as an $X$-module by the image of $U_{\alpha_{7}}$, and in $Q_{145678}$ one such module is generated by $U_{\alpha_3}$. The element $n_{13354321} n_{22343210} n_{10111111} n_{2} n_{4} n_{3} n_{5} n_{4} n_{2} n_{6} n_{5} n_{4} n_{3} n_{8} \in N_{G}(T)$ sends the root subgroups in $L_{124567}$ to those in $L_{134568}$ and sends $U_{00111111}$ to $U_{\alpha_7}$, while the element $n_{2} n_{4} n_{5} n_{6} n_{7} n_{8}$ sends the root subgroups in $L_{124567}$ to those of $L_{145678}$ and sends $U_{00111111}$ to $U_{\alpha_3}$. It follows that up to $G$-conjugacy there is at most one non-$G$-cr subgroup contained in an $A_{1} A_{5}$-parabolic with irreducible image in the Levi factor.

Let $Y \cong A_1 < A_7$ via $W(7)$. Then $Y$ is non-$G$-cr by Lemma \ref{lem:BMR} and contained in an $A_1 A_5$-parabolic subgroup of the Levi subgroup $A_7$ and hence in an $A_1 A_5$-parabolic subgroup of $G$. Hence $Y$ is a representative of the class of non-$G$-cr subgroups above.

From Table \ref{E8p7a1tab}, $\text{dim}(C_{L(G)}(Y)) = 3$, and $C_G(Y)^{\circ} \ge C_G(A_7) = T_{1}$. Moreover, a non-$G$-cr subgroup of $QX$ corresponding to $(0,1) \in \mathbb{V}_{X,Q}$ lies in $Q(4)X$, and thus centralises both $Z(Q) = Q(6)$ and the 1-dimensional subgroup of fixed points giving rise to the trivial $X$-module in $Q(2)/Q(3)$. Thus $C_G(Y)$ contains a 2-dimensional unipotent subgroup containing only positive root elements, and it follows that $C_G(Y)^\circ = U_2 T_1$.

\subsection{$L'=A_2 D_4$}
Let $P = P_{234578} = QL$ be the unique standard $A_2 D_4$-parabolic subgroup of $G$. Let $X$, $Y$ and $Z$ be representatives of the three $L'$-conjugacy classes of $L'$-irreducible subgroups, with the natural $A_2$-module $\lambda_1$ restricting to $X$, $Y$ and $Z$ with high weight $2$, and $V_{D_4}(\lambda_3) \downarrow X \cong V_{D_4}(\lambda_4) \downarrow Y \cong V_{D_4}(\lambda_1) \downarrow Z \cong 4 + 2^{[1]}$, so that the other two of these 8-dimensional $D_4$-modules restricts as $3 \otimes 1^{[1]}$. The actions of $L'$ on the levels of $Q$ are as follows:
\begin{align*}
Q/Q(2) \downarrow L' &= (00,\lambda_3) + (10,\lambda_4),\\
Q(2)/Q(3) \downarrow L' &= (10,\lambda_1) + (01,0), \\
Q(3)/Q(4) \downarrow L' &= (01,\lambda_3), \\
Q(4)/Q(5) \downarrow L' &= (00,\lambda_1) + (01,0), \\
Q(5)/Q(6) \downarrow L' & = (00,\lambda_4), \\
Q(6)/Q(7) \downarrow L' &= (10,0).
\end{align*} 
It follows that $\mathbb{V}_{X,Q} \cong \mathbb{V}_{Y,Q} \cong \mathbb{V}_{Z,Q} \cong K^{2}$. The two $L'$-modules occurring in the filtration of $Q$ with non-vanishing first cohomology group for $X$ (resp. $Y$ and $Z$) are generated as $L'$-modules by the images of the root groups $U_{\alpha_6}$ and $U_{10111100}$ (resp.\ $U_{10111100}$ and $U_{11122210}$; and $U_{\alpha_6}$ and $U_{11122210}$).

Whenever the partial map $\rho \, : \, \mathbb{V}_{X,Q} \to H^{1}(X,Q)$ is defined, let $X_{[a,b]}$ denote a complement to $Q$ in $QX$ corresponding to $(a,b) \in \mathbb{V}_{X,Q}$, and similarly for $Y_{[a,b]}$ and $Z_{[a,b]}$.

Note that $Q(5)/Q(6)$ contains a $Y$-module direct summand $2^{[1]}$, and $Q(4)/Q(5)$ contains a $Z$-module summand $2^{[1]}$. By Lemma \ref{lem:h2fora1} these modules have non-vanishing second cohomology group, and so the map $\rho$ is not necessarily defined everywhere. If the basis  of $\mathbb{V}_{Y,Q}$ is chosen to consist of a non-zero element from each group $H^{1}(Y,Q(2)/Q(3))$ and $H^{1}(Y,Q(3)/Q(4))$, and similarly for $Z$, with similar calculations to those of Section \ref{sec:E6D4} we find that the condition $ab = 0$ is necessary for the complements $Y_{[a,b]}$ or $Z_{[a,b]}$ to exist.

The typical element of $Z(L)$:
\[ h(u,t) \stackrel{\textup{def}}{=} h_1(t^2u^{-6})h_2(t)h_3(t^2u^{-3})h_4(t^2)h_5(tu^3)h_6(u^6)h_7(u^4)h_8(u^2) \]
acts as the scalar $u^{5}t^{-1}$ on the module generated by the image of $U_{\alpha_6}$, and as $tu^{-4}$ on the module generated by the image of $U_{10111100}$, hence Lemma \ref{lem:torus} applies and complements to $Q$ in $QX$ are $G$-conjugate to one of $X_{[0,0]}$, $X_{[0,1]}$, $X_{[1,0]}$ or $X_{[1,1]}$, while complements to $Q$ in $QY$ are conjugate to $Y_{[0,0]}$, $Y_{[1,0]}$ or $Y_{[0,1]}$, and complements to $Q$ in $QZ$ are conjugate to $Z_{[0,0]}$, $Z_{[0,1]}$ or $Z_{[1,0]}$.

The element $n_1 n_3 n_4 n_2 n_5 n_4 n_3 n_1 \in N_{G}(T)$ induces a non-trivial graph automorphism of the $D_4$ factor of $L'$. This fixes the $L'$-conjugacy class of $X$ whilst swapping $Y$ and $Z$, and also fixes the root group $U_{11122210}$ whilst swapping $U_{\alpha_6}$ and $U_{10111100}$. Therefore, $X_{[0,1]}$ and $X_{[1,0]}$ are conjugate in $G$; as are $Y_{[1,0]}$ and $Z_{[1,0]}$; and $Y_{[0,1]}$ and $Z_{[0,1]}$. Similarly, the element $n_{01122210} n_8 n_7 n_6 n_5 n_4 n_2 n_3 n_4 n_5 n_6$ swaps $U_{10111100}$ and $U_{11122210}$, and also stabilises $Y$ and swaps $X$ and $Z$. Therefore $X_{[0,1]}$ is $G$-conjugate to $Z_{[0,1]}$, and $Y_{[1,0]}$ is $G$-conjugate to $Y_{[0,1]}$. Finally, the element
\[ n_{22343210} n_{23465432} n_8 n_7 n_6 n_5 n_4 n_2 n_3 n_4 n_5 n_6 n_7 n_8 \]
normalises $Y$, swaps $X$ and $Z$, and stabilises $U_{\alpha_6}$. Hence $X_{[1,0]}$ and $Z_{[1,0]}$ are conjugate in $G$. It follows that there are at most two $G$-conjugacy classes of non-$G$-cr subgroups $A_{1}$ with irreducible image in a Levi factor of type $A_{2} D_{4}$, represented by $X_{[1,1]}$ and $X_{[0,1]}$.

Consider $A \cong A_1 < D_7$ via $T(10)$ and $B \cong A_1 \hookrightarrow A_1 A_1 G_2 < A_1 E_7$ via $(1,1^{[1]},6^{[1]})$. Then $A$ lies in an $A_2 D_4$-parabolic subgroup of the Levi subgroup $D_7$, and hence in an $A_2 D_4$-parabolic subgroup of $G$. Also, by Lemma~\ref{lem:BMR}, $B$ is non-$G$-cr if and only if it is non-$A_1 E_7$-cr. By \cite[Lemma 2.12]{MR2178661}, this is the case if and only if the image of $B$ in $E_7$ is non-$E_7$-cr, which we know is true from Theorem~\ref{thm:e7}. Furthermore, we know that the image of $B$ in $E_7$ lies in an $A_1 A_2 A_3$-parabolic subgroup of $E_7$ (cf.\ Section \ref{sec:E7A1A2A3p7}), and so $B$ lies in a parabolic subgroup of $G$ whose Levi factor contains a subgroup $A_1^{2} A_2 A_3$. The only possibility for this is $A_{2} D_{5}$; then by consideration of composition factors of $B$ on $L(G)$ we must have $V_{D_5}(\lambda_1) \downarrow B = 4 + 1^{[1]} \otimes 1^{[1]} + 0 = 4 + 2^{[1]} + 0^{2}$. In particular the image of $B$ in $D_{5}$ lies in a subgroup $D_{4}$, so $B$ is contained in an $A_{2} D_{4}$-parabolic subgroup of $G$, with irreducible image in the Levi factor.

Thus $A$ and $B$ are non-$G$-cr subgroups of $G$, each contained in an $A_{2} D_{4}$-parabolic subgroup with irreducible image in the Levi factor. Since $\textup{dim}(C_{L(G)}(B)) = 0$ we have $C_G(B)^{\circ} = 1$, and since $\textup{dim}(C_{L(G)}(A)) = 1$ we have $C_{G}(A)^{\circ} = C_G(D_7)^{\circ} = T_{1}$. This also shows that $A$ and $B$ are not $G$-conjugate.

\subsection{$L'=A_1^2 A_4$}

The four standard $A_1^{2} A_4$-parabolic subgroups of $G$ are $P_{235678}$, $P_{125678}$, $P_{124568}$ and $P_{123468}$. Let $X, Y \cong A_1 < A_1^2 A_4$ via $(1,1^{[1]},4)$ and $(1^{[1]},1,4)$, respectively. If $P = QL$ is one of these four parabolics and $X < L$, then in the filtration of $Q$ there is a unique $L$-module direct summand $(1,1,4)$, generated as an $L$-module by the image of $U_{\alpha}$ as in the table below. This gives rise to a unique indecomposable summand in the filtration of $Q$ for $X$ (resp.\ $Y$) with non-vanishing first cohomology group. Then the element of $N_{G}(T)$ given in the table below sends the root subgroups of $L$ to those of $L_{235678}$ and sends $U_{\alpha}$ to $U_{\alpha_4}$. By Lemma \ref{lem:torus}, there is at most one $G$-conjugacy class of non-$G$-cr complements to $Q$ in $QX$ and at most one $G$-conjugacy class of non-$G$-cr complements to $Q$ in $QY$ for each choice of parabolic, and such complements arising for different choices of parabolic subgroup are all conjugate in $G$.
\begin{center}
\begin{tabular}{c|c|c}
$P$ & Root $\alpha$ & Element of $N_{G}(T)$ \\ \hline
$P_{235678}$ & $\alpha_4$ & $1$ \\
$P_{125678}$ & $\alpha_{3} + \alpha_{4}$ & $n_{1} n_{3}$ \\
$P_{124568}$ & $00111110$ & $n_{1} n_{3} n_{4} n_{2} n_{5} n_{4} n_{6} n_{5} n_{7} n_{6} n_{8} n_{7} $ \\
$P_{123468}$ & $01122110$ & $n_{22343210} n_{10111111} n_{11110000} n_{2} n_{5} n_{4} n_{3} n_{6} n_{5} n_{4} n_{7} n_{6} n_{5} $
\end{tabular}
\end{center}
Furthermore, the element $n_{22343210} n_{23465432} n_{8} n_{7} n_{6} n_{5} n_{4} n_{2} n_{3} n_{4} n_{5} n_{6} n_{7} n_{8}$ stabilises the set of root subgroups in $L_{235678}$, swapping the two $A_1$ factors, while fixing $U_{\alpha_4}$. Hence a non-$G$-cr complement to $Q$ in $QY$ is $G$-conjugate to a subgroup of $QX$, so up to $G$-conjugacy there exists at most one non-$G$-cr subgroup $A_1$ with irreducible image in a Levi factor of type $A_1^2 A_4$.

Let $Z \cong A_1 < D_7$ via $T(8)$. Then $Z$ is non-$G$-cr by Lemma \ref{lem:BMR} and is contained in a $A_1^2 A_4$-parabolic subgroup of the Levi subgroup $D_7$, and hence in an $A_1^2 A_4$-parabolic subgroup of $G$. Therefore, $Z$ is a representative for the class of non-$G$-cr subgroups above. 

From Table \ref{E8p7a1tab}, $\text{dim}(C_{L(G)}(Z)) = 1$ and so $C_G(Z)^\circ = C_G(D_7)^\circ = T_1$.

\subsection{$L'=A_1 A_2 A_3$}

The four standard $A_1 A_2 A_3$-parabolic subgroups of $G$ are $P_{123567}$, $P_{123678}$, $P_{124678}$ and $P_{124578}$. We need to consider $X \cong A_1 \hookrightarrow A_1 A_2 A_3$ via $(1^{[1]},2,3)$. If $P = QL$ is one of these four parabolics, in the filtration of $Q$ there is a unique $L$-module direct summand on which each simple factor of $L$ acts non-trivially. This gives rise to a unique indecomposable $X$-module direct summand with non-vanishing first cohomology group, generated as an $X$-module by the image of $U_{\alpha}$ as in the table below. Then the corresponding element of $N_{G}(T)$ sends the root subgroups of $L$ to those of $L_{123567}$ and sends $U_{\alpha}$ to $U_{\alpha_{4}}$. By Lemma \ref{lem:torus}, there is at most one $G$-conjugacy class of non-$G$-cr complements to $Q$ in $QX$ for each choice of parabolic, and such non-$G$-cr complements arising for different choices of parabolic subgroup are all conjugate in $G$.

\begin{center}
\begin{tabular}{c|c|c}
$P$ & Root $\alpha$ & Element of $N_{G}(T)$ \\ \hline
$P_{123567}$ & $\alpha_4$ & $1$ \\
$P_{123678}$ & $00011000$ & $n_{8} n_{7} n_{6} n_{5} $ \\
$P_{124678}$ & $00111000$ & $n_{8} n_{7} n_{6} n_{5} n_{4} n_{2} n_{3} n_{1} n_{4} n_{3} $ \\
$P_{124578}$ & $00111100$ & $n_{8} n_{7} n_{6} n_{5} n_{4} n_{2} n_{3} n_{1} n_{4} n_{3} n_{5} n_{4} n_{2} n_{6} n_{5} n_{4} n_{7} n_{6} n_{5} n_{8} n_{7} n_{6}$
\end{tabular}
\end{center}

Let $Y \cong A_1 \hookrightarrow A_1 G_2 < E_7$ via $(1,6)$. Then $Y$ is non-$E_7$-cr by Theorem \ref{thm:e7}, hence is non-$G$-cr by Lemma \ref{lem:BMR}. In addition, from Section \ref{sec:E7A1A2A3p7} we know that $Y$ is contained in an $A_1 A_2 A_3$-parabolic subgroup of $E_7$, hence $Y$ lies in an $A_1 A_2 A_3$-parabolic subgroup of $G$, and is therefore a representative of the unique class of non-$G$-cr subgroups arising above.

From Table \ref{E8p7a1tab}, $\text{dim}(C_{L(G)}(Y)) = 3$ and so $C_G(Y)^\circ = C_G(E_7)^\circ = A_1$.


\section{Proof of Theorem \ref{thm:reductive} and Corollaries \ref{cor:variations}--\ref{cor:countable}} \label{sec:corollaries}

Having now proved Theorems \ref{thm:e6}--\ref{thm:e8}, in this section we prove Theorem \ref{thm:reductive} and all the corollaries stated in the introduction. Throughout, $G$ denotes an exceptional simple algebraic group, over an algebraically closed field of characteristic $p = 5$ or $7$.

\proof[Proof of Theorem \ref{thm:reductive}.] Let $X$ be a non-$G$-cr connected reductive subgroup of $G$. If $X$ is not simple, then since $p$ is good for $G$ it follows from \cite[Theorem 1.3]{MR2431255} that some simple factor of the derived subgroup $X'$ is non-$G$-cr. Conversely if some such simple factor of $X'$ is non-$G$-cr, then as this factor is normal in $X$ it follows from \cite[Theorem 3.10]{MR2178661} that $X$ is non-$G$-cr.

Thus it suffices to enumerate reductive subgroups of $XC_G(X)^{\circ}$ containing $X$, for each non-$G$-cr simple subgroup $X$ in Tables \ref{E6p5tab}--\ref{E8p7g2tab}. The results are precisely Tables \ref{tab:semisimple} and \ref{tab:reductive}. Where we have written ``$\infty$-many classes'', the distinct classes arise since a group of type $A_{1}T_{1}$ contains infinitely many pairwise non-conjugate 1-dimensional tori. For instance, $GL_{2}(K)$ contains the 1-dimensional tori $\left\{ \left(\begin{smallmatrix}t & 0 \\ 0 & t^{n} \end{smallmatrix}\right) \, : \, t \in K^* \right\}$ for $n \ge 1$. It remains to prove that no two subgroups given in Tables \ref{tab:semisimple} and \ref{tab:reductive} are $\textup{Aut}(G)$-conjugate. Firstly, notice that each reductive subgroup $X$ contains a \emph{unique} non-$G$-cr normal simple subgroup, so if $X_{1}$ and $X_{2}$ are non-$G$-cr reductive subgroups then we may assume $X_1$, $X_2 \le YC_G(Y)^{\circ}$ for some non-$G$-cr simple subgroup $Y$. Since $X_1$ and $X_2$ are isomorphic, inspection now shows that this is impossible unless either $X_1 = X_2$, or $C_G(Y) = A_{1}T_{1}$ and $X_1$, $X_2$ are each equal to the product of $Y$ with a 1-dimensional torus of $C_G(Y)$. Consideration of $L(G) \downarrow X_1$, $L(G)\downarrow X_2$ then shows that $X_1$ and $X_2$ can only be conjugate if the corresponding tori are $N_G(Y) = YC_G(Y)$-conjugate, which gives the result. \qed

\proof[Proof of Corollary \ref{cor:variations}] For $G$-cr subgroups of $G$, this follows from \cite[Theorems 3, 4]{MR1973585}. For non-$G$-cr subgroups, by Theorems \ref{thm:e6}--\ref{thm:e8} it suffices to inspect $L(G) \downarrow X$ for each subgroup $X$ in Tables~\ref{E6p5tab}--\ref{E8p7g2tab}. \qed

\proof[Proof of Corollary \ref{cor:overgroups}] Let $X$ be a non-$G$-cr reductive subgroup of $G$, and suppose that $X'$ is not contained in a proper subsystem subgroup of $G$. Then $Z(X)^{\circ} = 1$ and $C_G(X)^{\circ}$ is unipotent, otherwise $X$ centralises a non-trivial torus and lies in the corresponding Levi subgroup. Now, the first column of Tables \ref{E6p5tab}--\ref{tab:semisimple} gives a proper subsystem subgroup of $G$ containing $X$, for each non-$G$-cr subgroup $X$, with the exception of the simple groups given in (ii)--(v) of this corollary. These four subgroups cannot lie in a proper subsystem subgroup of the relevant group $G$, since all non-$G$-cr subgroups of subsystem subgroups appear elsewhere in Tables \ref{E6p5tab}--\ref{tab:semisimple} and therefore represent different conjugacy classes of subgroups. This proves that exactly one of (i)--(v) holds; furthermore each of the subgroups in (ii)--(iv) is uniquely determined up to conjugacy in the ambient group $G$, and the subgroups in (v) are determined up to $\operatorname{Aut}(G)$-conjugacy.

Now let $M$ be connected and maximal among reductive subgroups of $G$. If $M$ is $G$-reducible, then either $M$ is $G$-cr and therefore maximal among proper Levi subgroups of $G$, or $M$ is non-$G$-cr. In the latter case, either $G = E_{7}$, $p = 7$ or $G = E_6$, $p = 5$, with $M$ respectively conjugate to a non-$G$-cr subgroup from (iii) or (v) above. If instead $M$ is $G$-irreducible, then an application of the Borel-Tits Theorem shows that $M$ is in fact a maximal connected subgroup of $G$, as required. \qed

\proof[Proof of Corollary \ref{cor:max}] Let $M$ be maximal among proper reductive subgroups of $G$ and let $X \le M$ be $M$-irreducible and non-$G$-cr. If $Z(M)^{\circ} \neq 1$ then $M$ is a Levi subgroup of $G$, and so all $M$-cr subgroups of $M$ are $G$-cr by \cite[Proposition 3.2]{MR2167207}. Thus $M$ is semisimple. Also $X$ is semisimple as $X$ cannot centralise a non-trivial torus of $M$.

If $M$ is $G$-cr then either $M$ is the derived subgroup of a maximal Levi subgroup, or $M$ is $G$-irreducible and therefore a maximal connected subgroup of $G$. In the former case, $M$ is a subsystem subgroup of $G$ and so every $M$-irreducible subgroup is $G$-cr. We will consider the case that $M$ is a non-subsystem maximal connected subgroup of $G$ shortly.

If instead $M$ is non-$G$-cr, then from Corollary \ref{cor:overgroups} we know that either $G = E_{7}$, $p=7$ with $M = G_2$ as in part (iii); or $G = E_6$, $p = 5$ with $M = A_1$ as in part (v); or $X$ is an $M$-irreducible proper subgroup of $M$. The maximal $M$-irreducible subgroups of $M$ are of type $A_{1}$, $A_{1}\tilde{A}_{1}$ and $A_{2}$; we have proved in Section \ref{sec:E7E6p7} that the former of these is non-$G$-cr, while the latter two are centralisers (in $M$) of non-central semisimple elements, of order $2$ and $3$ respectively, and therefore lie in a proper subsystem subgroup of $G$. By Theorem \ref{thm:e7}, the only non-$G$-cr subgroups of $G$ are conjugate either to $M$ or its maximal subgroup $A_{1}$, or to a second class of subgroups $A_{1}$, which by Corollary \ref{cor:overgroups} do not lie in a proper subsystem subgroup of $G$, and hence do not lie in $M$. This justifies the entries in Table \ref{tab:max} with $M$ non-$G$-cr.

We may now assume that $M$ is $G$-irreducible, and therefore a maximal connected subgroup of $G$. Thus $M$ is one of the subgroups given by Lemma \ref{lem:maximalexcep}, and by hypothesis $M$ is not a subsystem subgroup. Furthermore if $G = E_{6}$ and $M = F_{4}$ or $C_{4}$, then $M$ is the centraliser of an involutary automorphism of $G$, and since $p \neq 2$ it follows from \cite[Corollary 3.21]{MR2178661} that every $M$-irreducible subgroup of $M$ is $G$-cr. Since $(G,p) = (E_6,5)$, $(E_7,5)$, $(E_7,7)$ or $(E_8,7)$, the possibilities for $M$ are therefore as follows: $M = A_{2}G_{2}$ or $A_{2}$ (two classes) when $G = E_{6}$; $M = G_{2}C_{3}$, $A_{1}F_{4}$, $A_{1}G_{2}$, $A_{1}A_{1}$ or $A_{2}$ when $G = E_{7}$; and $M = G_{2}F_{4}$, $B_{2}$ or $A_{1}A_{2}$ when $G = E_{8}$.

If $G = E_{6}$ then by Theorems \ref{thm:e6} and \ref{thm:reductive}, $X$ has type $A_{1}$ or $A_{1}A_{1}$. If $M = A_{2}G_{2}$ then the image of $X$ in the $G_{2}$ factor lies in a maximal subgroup $A_{1}A_{1}$ or $A_{2}$ of this factor. Now $L(G) \downarrow A_{2}G_{2} = L(A_2) + L(G_2) + (11,10)$ by \cite[p.\ 193]{MR1048074}. Since $L(A_2)$ restricted to a maximal subgroup $A_{1}$ has shape $4 + 2$, and $L(G_2) \downarrow A_1A_1 = (2,0)+(0,2) + (1,3)$ and $L(G_2) \downarrow A_2 = L(A_2) + 01 + 10$, it follows that $L(G) \downarrow X$ has at least three $3$-dimensional direct summands and two $4$-dimensional summands, with pairwise zero intersections. No possible $X$ in Table \ref{E6p5tab} or \ref{tab:semisimple} satisfies this.

On the other hand, if $M$ is a maximal subgroup of type $A_{2}$ and $X$ is a maximal subgroup $A_{1}$ of $M$, then using the restrictions $V_{G}(\lambda_1) \downarrow M = 22|11$ or $11|22$ and $L(G) \downarrow M = 41 + 14 + 11$ given in \cite[Tables 10.1,10.2]{MR2044850}, it follows that $V_{G}(\lambda_1) \downarrow X = T(8) + W(6) + 4^{2}$ or $T(8) + W(6)^{*} + 4^{2}$, and $L(G) \downarrow X = T(10)^{2}+T(6)^{2}+4^{3}+2$. In particular $X$ fixes a $1$-space on $V_{G}(\lambda_1)$ and lies in a subgroup of dimension at least $78 - \textup{dim}(V_{G}(\lambda_1)) = 51$. By Lemma \ref{lem:maximalexcep} the only such subgroups are contained in parabolic subgroups, hence $X$ is $G$-reducible. Since $X$ is fixed-point-free on $L(G)$, $X$ cannot centralise a non-trivial torus of $G$, and hence is non-$G$-cr. Finally, using Theorem \ref{thm:e6} and comparing the composition factors of $X$ on $L(G)$ and $V_{G}(\lambda_1)$ with Table \ref{E6p5tab}, we see that $X$ is $\textup{Aut}(G)$-conjugate to a subgroup $A_1 < A_{1}A_{5}$ via $(1,W(5))$ as in Table \ref{tab:max}.

For $G = E_{7}$ or $E_{8}$ we proceed similarly, considering each possible subgroup $M$ and its action on $L(G)$ and $V_{G}(\lambda_7)$ when $G = E_{7}$; these are given by \cite[Tables 10.1,\ 10.2]{MR2044850}. In each case, we find that no simple subgroup of $M$ of the appropriate isomorphism type can act on $L(G)$ or $V_{G}(\lambda_7)$ in the manner given in Tables \ref{E7p5tab}--\ref{E8p7g2tab}, except for the subgroups of $M$ given in the final column of Table \ref{tab:max}. In these cases, the restriction is compatible with the non-$G$-cr subgroups $X$ given in the fourth column there. Thus if we show that these subgroups of $M$ are non-$G$-cr, then as their conjugacy class is determined by their action on $L(G)$, it will follow that they are indeed conjugate to the non-$G$-cr subgroups $X$ in column 4. Furthermore since none of the subgroups listed in column 4 centralise a non-trivial reductive subgroup of $G$, it follows that no semisimple, non-simple, non-$G$-cr subgroup of $G$ lies in any such $M$.

It remains to prove that each subgroup $X$ of $M$ given in the final column of Table \ref{tab:max} is in fact non-$G$-cr. Firstly, using the fact that the action of $X$ on $L(G)$ agrees with the appropriate non-$G$-cr subgroup in the fourth column of Table \ref{tab:max}, we see that $L(G) \downarrow X$ has no trivial direct summands. If $X$ centralises a non-trivial torus $S$ of $G$, then $L(C_{G}(S)) = L(S) \oplus L(C_{G}(S)')$ is a direct summand of $L(G) \downarrow C_G(S)$; a complement is given by the sum of those root spaces of $L(G)$ not centralised by the action of $S$ (see also \cite[Lemma 3.9]{MR2608407}). Since $X < C_{G}(S)$, this implies that $L(S)$ is a trivial direct summand of $L(G) \downarrow X$, a contradiction. Hence each possible subgroup $X$ is $G$-indecomposable.

We now show that every subgroup $X$ of each subgroup $M$ lies in a parabolic subgroup of $G$. For $G = E_{7}$ it suffices to note that $X$ has a non-zero fixed point on $V_{G}(\lambda_7)$; the corresponding stabiliser has dimension at least $133 - 56 = 77$, and every subgroup of $G$ of such a dimension lies in a parabolic subgroup by Lemma \ref{lem:maximalexcep}. Thus $X$ is indeed $G$-reducible. For $G = E_{8}$ it is shown in \cite[Section 3.3]{MR2044850} that a subgroup of type $A_{1}$ with the same composition factors on $L(G)$ as the subgroup $A_{1} < A_{8}$ via $W(8)$ must in fact be conjugate to a subgroup of $A_{8}$; it follows that $X$ is conjugate to this non-$G$-cr subgroup. Finally the subgroup $G_{2} \hookrightarrow G_{2}G_{2} < G_{2}F_{4}$ via $(10,10)$ is shown to be $G$-reducible in \cite[Lemma 7.13]{Tho1}, completing the proof. \qed

\proof[Proof of Corollary \ref{cor:restricted}] This is a matter of inspecting Tables \ref{E6p5tab}--\ref{tab:semisimple}, noting that only the non-$G$-cr subgroups of type $G_{2}$ and $A_{1}G_{2}$ are restricted; the $A_{1}$ factor of the non-$G$-cr subgroup $A_{1}G_{2}$ of $G = E_{8}$ is restricted because it is the centraliser of a subsystem subgroup $E_{7}$, and $L(G) \downarrow A_{1}E_{7} = L(A_{1}) + L(E_{7}) + (1,\lambda_7)$ (see \cite[1.8]{MR1048074}), so the high weights of $L(G) \downarrow A_{1}$ are at most 2. \qed

\proof[Proof of Corollary \ref{cor:trivcentraliser}] If $X$ is a reductive subgroup of $G$ with $C_{G}(X)^{\circ} = 1$ then $Z(X)^{\circ} = 1$ and $X$ lies in no proper Levi subgroup of $G$, since then $X$ would centralise a non-trivial torus. Thus either $X$ is $G$-irreducible or $X$ is non-$G$-cr. Inspecting Tables \ref{E6p5tab}--\ref{tab:semisimple}, we see that each non-$G$-cr subgroup $X$ with $C_{G}(X)^{\circ} = 1$ appears in Table \ref{tab:centraliser}, as required. \qed

\proof[Proof of Corollaries \ref{cor:fincomps} and \ref{cor:countable}] Theorems \ref{thm:e6}--\ref{thm:reductive} show that for each exceptional simple algebraic group $G$, there are only a finite number of non-$G$-cr reductive subgroups of $G$ having a specified set of composition factors on $L(G)$, which immediately implies each of these corollaries. \qed


\section{Restrictions of $G$-modules to non-$G$-cr subgroups} \label{sec:restrictions}
Let $V$ be either $L(G)$, or $V_{27}$ or $V_{56}$ when $G$ is of type $E_{6}$ or $E_{7}$ respectively, and let $X$ be one of the non-$G$-cr simple subgroups of $G$ listed in Tables \ref{E6p5tab}--\ref{E8p7g2tab}. Here, we justify the restrictions $V \downarrow X$ given in those tables, which have been used in proving aspects of Theorems \ref{thm:e6}--\ref{thm:e8}.

We begin by considering the cases where $X$ is properly contained in a proper semisimple subgroup $H$ of $G$, which is then given in the tables. For example, in the first line of Table \ref{E6p5tab}, $X$ is contained in a Levi subgroup $H$ of type $A_{5}$ in $G = E_{6}$. Now, each subgroup $H$ occurring is either a subsystem subgroup of $G$, or is contained in one of a small number of known maximal connected subgroups (cf.\ Lemma \ref{lem:maximalexcep}). In the latter case, the action of the maximal connected subgroup on the low-dimensional $G$-modules is given explicitly in \cite[Table 10.1]{MR2044850}, and the action of $H$ is straightforward to determine from this.

If $H$ is a subsystem subgroup occurring, and if either $H$ has maximal rank or $G = E_{6}$ and $H$ is a subgroup of type $D_{5}$, then $V \downarrow H$ is given in \cite[Lemmas 11.2, 11.8, 11.10]{MR2883501}. Moreover, if $G = E_{7}$ and $H = E_{6}$ then $V \downarrow H$ is given in \cite[Tables 8.2, 8.6]{MR1329942}. From these it is straightforward to determine $V \downarrow H$ for all other type of subgroup $H$ occurring.

If a factor $H_{0}$ of $H$ is classical or of type $G_{2}$, then the given embedding determines the action of $X$ on $V_{H_0}(\lambda_1)$ (considered as an $H$-module in the obvious way). From this it is straightforward to determine the high weights of related modules, such as symmetric and alternating powers, and tensor products. Moreover as discussed in the proof of Lemma \ref{lem:tiltingnoH1}, since $p > 3$ each of the modules $S^{2}(V)$, $\bigwedge^{2}(V)$ and $\bigwedge^{3}(V)$ occurs as a direct summand in a tensor power of $V$. It then follows from \cite[Corollary 1.3]{DonkinFiltration} that if $V$ has a filtration by Weyl modules (respectively, dual Weyl modules), so too do these symmetric and alternating powers. This is sufficient information to determine $V \downarrow X$, unless some factor $H_{0}$ has type $D_{n}$ and $V \downarrow H_{0}$ involves a spin module; we will address this problem shortly.

As an example, consider the final line of Table \ref{E6p5tab}, where $G = E_{6}$, $p = 5$ and where $X \hookrightarrow H = A_{1} A_{5}$ via $(1^{[r]},W(5)^{[s]})$ with $rs = 0$. From \cite[Propositions 2.1, 2.3]{MR1329942} it follows that
\begin{align*}
L(G) \downarrow H &= L(A_1) + L(A_{5}) + (1,\lambda_{3}),\\
V_{27} \downarrow H &= (0,\lambda_{4}) + (1,\lambda_1)
\end{align*}
where the direct-sum decompositions follow as no factors here can extend another indecomposably. Now, the $A_{5}$-module $\lambda_{4} = \bigwedge^{2}(\lambda_1)^{*}$ restricts to the image of $X$ as $\bigwedge^{2}(W(5)^{[s]})^{*}$. Since $(W(5)^{[s]})^{*}$ has high weights $5^{s+1}$ and $3 (5^{s})$, it follows easily that $\lambda_{4} \downarrow X$ has high weights $8 (5^s)$, $4 (5^{s})$, $0$ and $0$. Since $\bigwedge^{2}(W(5))^{\ast}$ has a filtration by duals of Weyl modules, and since $W_{X}(4)$ and $W_{X}(0)$ are irreducible, while $W_{X}(8) = V_{X}(8)|0$, it follows that $\bigwedge^{2}(W(5)^{[s]})^{\ast} = 4^{[s]} \oplus T(8)^{[s]}$. The other $H$-direct summands of $L(G)$ and $V_{27}$ are also constructed as tensor products and alternating powers of $(0,\lambda_1)$ and $(1,0)$, and we proceed in an entirely similar manner. The restriction in Table \ref{E6p5tab} follows.

Now consider spin modules for $D_{n}$. Given a sub-torus of a simply-connected group of type $D_{n}$, the weights of this torus on the spin modules $V_{D_n}(\lambda_{n-1})$ and $V_{D_n}(\lambda_n)$ are straightforward to determine from the weights of this torus on the natural $2n$-dimensional module; an example of such a calculation is given in \cite[pp.\ 195--197]{MR1057341}. For the exact module structure, we make use of \textsc{Magma}'s functionality to work within groups of Lie type. In particular, \textsc{Magma} allows us to explicitly construct a finite quasi-simple subgroup $X(q) < X$, as an irreducible matrix group of degree $2n$ over a field of size $q$, preserving an explicit quadratic form. We can then take a pre-image of $X(q)$ under the natural $2n$-dimensional representation of a group of Lie type $D_{n}$, and then take the image of this under a representation of high weight $\lambda_{n-1}$ or $\lambda_n$.

As an example, let $G$ be simple of type $E_{8}$, let $L'$ be a simple subgroup of type $D_{7}$ as in Section \ref{sec:E8D7}, and let $X$ be a simple subgroup of type $G_{2}$, with $V_{D_{7}}(\lambda_{1}) \downarrow X = 01$. Using \textsc{Magma} as above, we find that the spin modules $V_{D_{7}}(\lambda_{6})$ and $V_{D_{7}}(\lambda_7)$ restrict to a finite subgroup $X(7^{2})$ as uniserial modules with composition factor dimensions $38$ and $26$. Since the composition factors of $X$ on $L(G)$ must agree with those of the image of $X$ in $D_{7}$, which are given in \cite[Table 8.1]{MR1329942}, we deduce that $V_{D_7}(\lambda_6) \downarrow X = 11|20$ and $V_{D_7}(\lambda_6) \downarrow X = 20|11$, or vice-versa.

We now need to consider the case where $X$ is contained in no proper reductive overgroup of $G$. First consider the case where $X$ type $G_2$ and $G = E_7$, as in Table \ref{E7p7g2tab}. The only possible parabolic subgroups of $E_{7}$ containing $X$ have Levi factor of type $E_{6}$, and the image of $X$ under projection to the Levi factor lies in a proper subgroup $F_{4}$. Using \cite[Tables 8.4, 8.6]{MR1329942} and \cite[Tables 10.1, 10.2]{MR2044850}, we deduce that the restriction of $V$ to $X$ has the composition factors as given in Table \ref{E7p7g2tab}. To determine the exact structure of the restriction of $V$ to $X$, we consider the maximal subgroup $A_1$ of $X$, call it $Y$. From Section \ref{sec:E7E6p7} we know that $Y$ is conjugate to a subgroup $A_1 < A_7$ via $W(7)$ and therefore $V \downarrow Y$ is given in Table \ref{E7p7a1tab}. We claim that this enough to determine $V \downarrow X$. To prove the claim we first take $V = V_{56}$. We know that the $X$-composition factors of $V$ are $20^2 / 00^4$. Since $V$ is self-dual and $H^1(G_2,20) \cong K$, the possibilities for $V \downarrow X$ are thus $T(20)^2$, $T(20) + 20 + 00^2$, $W(20) + W(20)^{*} + 00^2$ or $20^2 + 00^4$. Moreover, we know that $V \downarrow Y = T(12)^2 + T(8)^2 = (0 | 12 | 0)^2 + (4|8|4)^2$. Thus the fixed-point space of $X$ on $V$ is at most $2$-dimensional, and therefore $V \downarrow X = T(20)^2$ as claimed in Table \ref{E7p7g2tab}. 

Now suppose $V = L(G)$. As given in Table \ref{E7p7a1tab}, the fixed-point space of $Y$ on $V$ is $1$-dimensional. The $X$-composition factors of $V$ are $11 / 20^3 / 01 / 00^3$. Suppose $M, N$ are two such composition factors. Then $\text{Ext}^1_X(M,N) \neq 0$ precisely when $\{M,N\} = \{11,20\}$ or $\{20,00\}$, in which case $\text{Ext}^1_X(M,N) \cong K$. This leaves a number of possibilities for $V \downarrow X$ but only one, namely $T(11) + T(20) + 01$, has fixed-point space of dimension at most $1$.

Finally, consider the case where $X$ has type $A_1$ and $G = E_6$, $p = 5$, with $X$ having no proper reductive overgroups in $G$. Similarly to the above, considering the action of $A_1 < D_5$ via $2 \otimes 2^{[1]} + 0$, the composition factors of $X$ on $V_{27}$ are $16/12/8/0^2$ and those on $L(G)$ are $22/16^2/14/12/10/8^2/2/0$. By \cite{MR925847}, $X$ cannot fix a $2$-space on $V_{27}$ or on its dual, which implies that $T(8)$ occurs as a direct summand, since no other composition factors extend $0$ indecomposably. This in turn implies that unipotent elements of $X$ must act on $V_{27}$ with at least two Jordan blocks of size $5$, by \cite[Lemma 2.1.7(iii)]{MR2044850}. Since these elements also have order $5$, by \cite[Table 5]{MR1351124} they lie in one of the classes $A_4+A_1$, $A_4$, $D_4(a_1)$ or $A_3 + A_1$. If $X$ were to act on $V_{27}$ as $16 + 12 + T(8)$, calculations in \textsc{Magma} show that the Jordan blocks of the unipotent elements would be $5^4/3^2/1$, which is impossible by [loc.~cit.]. Therefore $V_{27} \downarrow X = W(16) + T(8)$ or the dual of this. Note that the two possible actions of $X$ on $V_{27}$, being duals of one another and not self-dual, give two distinct $G$-classes of non-$G$-cr subgroups. The structure of $V_{27} \downarrow X$ now tells us that the unipotent elements have Jordan block structure $5^5/2$, hence lie in class $A_4 + A_1$. This in turn tells us that these have Jordan blocks $5^{15}/3$ on $L(G)$. Since $L(G)$ is self-dual and $L(X)$ gives the unique composition factor of $L(G)$ having high weight $2$, it is a direct summand of $L(G)$. The composition factor of high weight $14$ does not extend any other factor indecomposably, hence is also a direct summand. Finally, none of the factors or Weyl modules of high weight $22$, $16$, $12$, $10$, $8$ or $0$ have dimension a multiple of $5$, hence none of these can be a direct summand of $L(G)$. The only remaining possibility is $L(G) \downarrow X = T(22) + T(10) + 14 + 2$, as in Table~\ref{E6p5tab}.

\section{{Tables for Theorems \ref{thm:e6}--\ref{thm:reductive}}}
\label{sec:tables}

In this section, we give the tables referred to in Theorems \ref{thm:e6}--\ref{thm:reductive}. Each line of a table corresponding to the exceptional Lie type $G$, characteristic $p$, and subgroup type $A_{1}$ or $G_{2}$ gives a representative of an $\textup{Aut}(G)$-conjugacy class of non-$G$-cr subgroups of that type. For each simple non-$G$-cr subgroup $X$ we also give the action of $X$ on the adjoint module $L(G)$, and on the module $V_{27} = V_{G}(\lambda_1)$ when $G$ has type $E_{6}$, and on $V_{56} = V_{G}(\lambda_7)$ when $G$ has type $E_{7}$. Finally, we also give the connected centraliser $C_{G}(X)^{\circ}$ for each $X$. The notation for modules and embeddings, e.g. `$X < A_{5}$ via $W(5)$' is explained in Section \ref{sec:notation}.

In Table \ref{E6p5tab} each given $\textup{Aut}(G)$-class of subgroups splits into two $G$-conjugacy classes, interchanged by the graph automorphism of $G$.

\begin{longtable}{>{\raggedright\arraybackslash}p{0.25\textwidth - 2\tabcolsep}>{\raggedright\arraybackslash}p{0.25\textwidth - 2\tabcolsep}>{\raggedright\arraybackslash}p{0.4\textwidth - 2\tabcolsep}>{\raggedright\arraybackslash}p{0.1\textwidth-\tabcolsep}@{}}
\caption{Non-$G$-cr subgroups of type $A_1$ in $G = E_6$, $p=5$. \label{E6p5tab}} \\
\hline
Non-$G$-cr subgroup $X \cong A_1$ & $V_{27} \downarrow X$ & $L(G) \downarrow X$ & $C_{G}(X)^{\circ}$ \\
\hline
$X < A_5$ via $W(5)$ & $T(8) + W(5)^2 + 4$ & $T(10) + 9^2 + T(6) + T(5)^2 + 4 + 0^3$ & $\bar{A}_1$ \\
$X < D_5$ via $T(6)$ & $T(7) + T(6) + W(5) + 0$ & $T(10) + T(7)^2 + T(6)^2 + W(5) + W(5)^* + 4 + 0$ & $T_1$ \\
$X < D_5$ via $T(8)$ & $W(10) + T(8) +  4 + 0$ & $14 + T(10) + W(10) + W(10)^* + T(6) + 4^2 + 0$ & $T_1$ \\
$X \hookrightarrow \bar{A}_1 A_5$ via $(1^{[r]},W(5)^{[s]})$ $(rs=0)$ & $1^{[r]} \otimes W(5)^{[s]} + T(8)^{[s]} + 4^{[s]}$ & $1^{[r]} \otimes 9^{[s]} + 1^{[r]} \otimes T(5)^{[s]} + 2^{[r]} + T(10)^{[s]} + T(6)^{[s]} + 4^{[s]}$ & 1 \\
$X$ in no proper reductive overgroup (see Section \ref{sec:E6D5again}) & $W(16) + T(8)$ & $T(22) + 14 + T(10) + 2$ & 1 \\
\hline
\end{longtable}

\begin{longtable}{>{\raggedright\arraybackslash}p{0.205\textwidth - 2\tabcolsep}>{\raggedright\arraybackslash}p{0.3\textwidth - 2\tabcolsep}>{\raggedright\arraybackslash}p{0.35\textwidth - 2\tabcolsep}>{\raggedright\arraybackslash}p{0.145\textwidth-\tabcolsep}@{}}
\caption{Non-$G$-cr subgroups of type $A_1$ in $G = E_7$, $p=5$. \label{E7p5tab}} \\
\hline
Non-$G$-cr subgroup $X \cong A_1$ & $V_{56} \downarrow X$ & $L(G) \downarrow X$ & $C_{G}(X)^{\circ}$ \\
\hline
$X < A_5$ via $W(5)$ & $9 + T(5) + W(5)^3 + (W(5)^*)^3$ & $T(10) + T(8)^6 + T(6) + 4^7 + 0^8$ & $U_6 A_2$ \\
$X < A_5'$ via $W(5)$ & $T(8)^2 + W(5)^2 + (W(5)^*)^2 + 4^2 + 0^2$ & $T(10) + 9^2 +  T(8)^2 +  T(6) +  T(5)^2 + W(5)^2 + (W(5)^*)^2 + 4^3 + 0^4$ & $U_2 \bar{A}_1 T_1$ \\
$X < A_6$ via $W(6)$ & $W(10) + W(10)^* + T(6)^2 + W(6) + W(6)^*$ & $T(12)^3 + T(10) + T(8)^2 + W(6) + W(6)^* + 4^3 + 2 + 0$ & $U_2 T_1$ \\
$X < A_7$ via $W(7)$ & $W(12) + W(12)^* + T(8)^2 + 4^2$ & $T(16) + 14 + T(12) + T(10)^2 + T(8) + 4^3 + 2$ & $U_1$ \\
$X \hookrightarrow \bar{A}_1 A_5$ via $(1^{[r]},W(5)^{[s]})$ $(rs=0)$ & ${1^{[r]} \otimes W(5)^{[s]}} + {1^{[r]} \otimes (W(5)^*)^{[s]}} + 9^{[s]} + T(5)^{[s]} + W(5)^{[s]} +  (W(5)^*)^{[s]}$ & $ 2^{[r]} + {(1^{[r]} \otimes T(8)^{[s]})^2} + {(1^{[r]} \otimes 4^{[s]})^2} + (1^{[r]})^2 + T(10)^{[s]} + (T(8)^{[s]})^2 + T(6)^{[s]} + (4^{[s]})^3 + 0$ & $U_2 T_1$ \\
$X \hookrightarrow \bar{A}_1 A_5'$ via $(1^{[r]},W(5)^{[s]})$ $(rs=0)$ & ${1^{[r]} \otimes W(5)^{[s]}} + {1^{[r]} \otimes (W(5)^*)^{[s]}} + (T(8)^{[s]})^2 + (4^{[s]})^2 + 0^2$ & $2^{[r]} + {1^{[r]} \otimes 9^{[s]}} + {1^{[r]} \otimes T(5)^{[s]}} + {1^{[r]} \otimes W(5)^{[s]}} + {1^{[r]} \otimes (W(5)^*)^{[s]}} + T(10)^{[s]} + (T(8)^{[s]})^2 + T(6)^{[s]} + (4^{[s]})^3 + 0$ & $U_2 T_1$ $(r + 1 \neq s)$, $U_3 T_1$ $(r + 1 = s)$ \\
$X \hookrightarrow A_2 A_5$ via $(2^{[r]},W(5)^{[s]})$ $(rs=0)$ & ${2^{[r]} \otimes W(5)^{[s]}} + {2^{[r]} \otimes (W(5)^*)^{[s]}} + 9^{[s]} + T(5)^{[s]}$ & $4^{[r]} + {(2^{[r]} \otimes T(8)^{[s]})^2} + {(2^{[r]} \otimes 4^{[s]})^2} + 2^{[r]} + T(10)^{[s]} + T(6)^{[s]} + 4^{[s]}$ & $1$ \\
$X \hookrightarrow D_5$ via $T(6)$ & $T(7)^2 + T(6)^2 + W(5) + W(5)^* + 0^4$ & $T(10) + T(7)^4 + T(6)^4 + W(5)^2 + (W(5)^*)^2 + 4 + 0^4$ & $\bar{A}_1 T_1$\\
$X \hookrightarrow D_5$ via $T(8)$ & $W(10) + W(10)^* + T(8)^2 + 4^2 +0^4$ & $14 + T(10) + W(10)^2 + (W(10)^*)^2 + T(8)^2 + T(6) + 4^4 + 0^4$ & $U_2 \bar{A}_1 T_1$ \\
$X \hookrightarrow \bar{A}_1 D_5$ via $(1^{[r]},T(6)^{[s]})$ $(rs=0)$ & ${1^{[r]} \otimes T(6)^{[s]}} + (1^{[r]})^2 + (T(7)^{[s]})^2 + W(5)^{[s]} + (W(5)^*)^{[s]} $ & $2^{[r]} + {(1^{[r]} \otimes T(7)^{[s]})^2} + {1^{[r]} \otimes T(6)^{[s]}} + {1^{[r]} \otimes W(5)^{[s]}} + {1^{[r]} \otimes (W(5)^*)^{[s]}} + T(10)^{[s]} + (T(6)^{[s]})^2 + 4^{[s]} + 0$ & $T_1$ $(s \neq r,r+1)$, $U_{1}T_{1}$ ($s = r+1$), $U_2 T_1$ ($s=r$) \\
$X \hookrightarrow \bar{A}_1 D_5$ via $(1^{[r]},T(8)^{[s]})$ $(rs=0)$ & ${1^{[r]} \otimes T(8)^{[s]}} + (1^{[r]})^2 + W(10)^{[s]} + (W(10)^*)^{[s]} + (4^{[s]})^2$ & $2^{[r]} + {1^{[r]} \otimes W(10)^{[s]}} + {1^{[r]} \otimes (W(10)^*)^{[s]}} + {1^{[r]} \otimes T(8)^{[s]}} + {(1^{[r]} \otimes 4^{[s]})^2} + 14^{[s]} + T(10)^{[s]} + T(6)^{[s]} + 0$ & $T_1$ \\
$X < E_6$ (see Section \ref{sec:E7D5again}) & $W(16) + W(16)^* + T(8)^2 + 0^2$ & $T(22) + W(16) + W(16)^* + 14 + T(10) + T(8)^2 + 2 + 0$ & $U_2 T_1$ \\
\hline
\end{longtable}

\begin{longtable}{>{\raggedright\arraybackslash}p{0.31\textwidth - 2\tabcolsep}>{\raggedright\arraybackslash}p{0.24\textwidth - 2\tabcolsep}>{\raggedright\arraybackslash}p{0.36\textwidth - 2\tabcolsep}>{\raggedright\arraybackslash}p{0.09\textwidth-\tabcolsep}@{}}
\caption{Non-$G$-cr subgroups of type $A_1$ in $G = E_7$, $p=7$. \label{E7p7a1tab}} \\
\hline
Non-$G$-cr subgroup $X \cong A_1$ & $V_{56} \downarrow X$ & $L(G) \downarrow X$ & $C_{G}(X)^{\circ}$ \\
\hline
$X < A_7$ via $W(7)$ & $T(12)^2 + T(8)^2$ & $T(16) + T(14) + T(12) + T(10) + T(8)^3 + 6$ & $U_1$ \\
$X \hookrightarrow A_1 G_2$ via $(1,6)$ & $T(11) + T(9)^2 + T(7)$ & $T(14) + T(10)^4 + T(8)^2 + 6^3$ & $1$ \\
\hline
\end{longtable}

\begin{longtable}{>{\raggedright\arraybackslash}p{0.56\textwidth - 2\tabcolsep}>{\raggedright\arraybackslash}p{0.12\textwidth - 2\tabcolsep}>{\raggedright\arraybackslash}p{0.22\textwidth - 2\tabcolsep}>{\raggedright\arraybackslash}p{0.1\textwidth-\tabcolsep}@{}}
\caption{Non-$G$-cr subgroups of type $G_2$ in $G = E_7$, $p=7$. \label{E7p7g2tab}} \\
\hline
Non-$G$-cr subgroup $X \cong G_2$ & $V_{56} \downarrow X$ & $L(G) \downarrow X$ & $C_{G}(X)^{\circ}$ \\
\hline
$X$ in no proper reductive overgroup (see Section \ref{sec:E7E6p7}) & $T(20)^2$ & $T(11) + T(20) + 01$ & $U_1$ \\
\hline
\end{longtable}

\begin{longtable}{>{\raggedright\arraybackslash}p{0.33\textwidth - 2\tabcolsep}>{\raggedright\arraybackslash}p{0.55\textwidth - 2\tabcolsep}>{\raggedright\arraybackslash}>{\raggedright\arraybackslash}p{0.12\textwidth-\tabcolsep}@{}}
\caption{Non-$G$-cr subgroups of type $A_1$ in $G = E_8$, $p=7$. \label{E8p7a1tab}} \\
\hline
Non-$G$-cr subgroup $X \cong A_1$ & $L(G) \downarrow X$ & $C_{G}(X)^{\circ}$ \\
\hline
$X < A_7$ via $W(7)$ & $T(15)^2 + T(14) + T(12)^2 + T(10)  + T(9)^2 + T(8)^3 + T(7)^2 + W(7) + W(7)^* + 6 + 0$ & $U_2 T_1$ \\
$X < A_7'$ via $W(7)$ & $T(16) + T(14) + T(12)^5 + T(10) + T(8)^7 + 6 + 0^3$ & $U_5 \bar{A}_1$\\
$X < A_8$ via $W(8)$ & $T(18)^2 + T(16) + T(14)^3 + T(10)^2 + T(8) + 6^5 + 2 $ & $1$ \\
$X \hookrightarrow \bar{A}_1 A_7'$ via $(1^{[r]},W(7)^{[s]})$ $(rs=0)$ & $2^{[r]} + {(1^{[r]} \otimes T(12)^{[s]})^2} + {(1^{[r]} \otimes T(8)^{[s]})^2} + T(16)^{[s]} + T(14)^{[s]} + T(12)^{[s]} + T(10)^{[s]} + (T(8)^{[s]})^3 + 6^{[s]}$ & $U_1$ \\
$X < D_7$ via $T(8)$ & $T(14) + 13^2 + T(11)^2 + T(10)^2 + T(9)^2 + T(8)^3 + T(7)^2 + W(7) + W(7)^* + 6^3 + 0$ & $T_1$ \\
$X < D_7$ via $T(10)$ & $T(18) + T(16)^2 + T(14) + W(14) + W(14)^* + T(10)^4 + T(8)^2 + 6^3 + 0$ & $T_1$ \\
$X < D_7$ via $T(12)$ & $T(22) + W(21) + W(21)^* + T(15)^2 + T(14) + T(12)^2 + T(10) + T(9)^2 + 6 + 0$ & $U_{2}T_1$ \\
$X \hookrightarrow A_1 G_2 < E_7$ via $(1,6)$ & $T(14) +T(11)^2 + T(10)^4 + T(9)^4 + T(8)^2 + T(7)^2 + 6^3 + 0^3$ & $\bar{A}_1$ \\
$X \hookrightarrow \bar{A}_1 A_1 G_2 < \bar{A}_1 E_7$ via $(1^{[r]},1^{[s]},6^{[s]})$ $(rs=0)$ & $2^{[r]} + {1^{[r]} \otimes T(11)^{[s]}} + {(1^{[r]} \otimes T(9)^{[s]})^2} + {1^{[r]} \otimes T(7)^{[s]}} + T(14)^{[s]} + (T(10)^{[s]})^4 + (T(8)^{[s]})^2 + (6^{[s]})^3$ & $1$ $(r \neq s)$, $U_1$ $(r = s)$  \\
\hline
\end{longtable}

\vspace{-3pt}
 
\begin{longtable}{>{\raggedright\arraybackslash}p{0.4\textwidth - 2\tabcolsep}>{\raggedright\arraybackslash}p{0.35\textwidth - 2\tabcolsep}>{\raggedright\arraybackslash}>{\raggedright\arraybackslash}p{0.25\textwidth-\tabcolsep}@{}}
\caption{Non-$G$-cr subgroups of type $G_2$ in $G = E_8$, $p=7$. \label{E8p7g2tab}} \\
\hline
Non-$G$-cr subgroup $X \cong G_2$ & $L(G) \downarrow X$ & $C_{G}(X)^{\circ}$ \\
\hline
$G_2 < E_7$ & $T(11) + T(20)^5 + 01 + 00^3$ &$U_5 \bar{A}_1$ \\
$G_2 \hookrightarrow G_2 G_2 < G_2 F_4$ via $(10,10)$ & $30 + T(11) + 11 + 01^3$ & $1$ \\
\hline
\end{longtable}

\vspace{-3pt}

\begin{longtable}{>{\raggedright\arraybackslash}p{0.06\textwidth - 2\tabcolsep}>{\raggedright\arraybackslash}p{0.06\textwidth - 2\tabcolsep}>{\raggedright\arraybackslash}p{0.62\textwidth-\tabcolsep}@{}}
\caption{Non-simple, non-$G$-cr semisimple subgroups} \label{tab:semisimple} \\
\hline
$G$ & $p$ & Non-$G$-cr semisimple subgroup  \\
\hline
$E_6$ & 5 & $\bar{A}_1 A_1 < \bar{A}_1 A_5$ where $A_1 < A_5$ via $W(5)$ \\
$E_7$ & 5 & $A_2 A_1 < A_2 A_5$ where $A_1 < A_5$ via $W(5)$ \\
& 5 & $A_1 A_1 < A_2 A_5$ where $A_1 < A_2$ via $2$ and $A_1 < A_5$ via $W(5)$ \\
& 5 & $\bar{A}_1 A_1 < \bar{A}_1 A_5$ where $A_1 < A_5$ via $W(5)$ \\
& 5 & $\bar{A}_1 A_1 < \bar{A}_1 A_5'$ where $A_1 < A_5'$ via $W(5)$ \\
& 5 & $\bar{A}_1 A_1 < \bar{A}_1 D_5$ where $A_1 < D_5$ via $T(6)$ \\
& 5 & $\bar{A}_1 A_1 < \bar{A}_1 D_5$ where $A_1 < D_5$ via $T(8)$ \\
$E_8$ & 7 & $\bar{A}_1 A_1 < \bar{A}_1 A_7$ where $A_1 < A_7'$ via $W(7)$ \\
& 7 & $\bar{A}_1 A_1 < \bar{A}_1 E_7$ where $A_1 \hookrightarrow A_1 G_2  < E_7$ via $(1,6)$ \\
& 7 & $\bar{A}_1 G_2 < \bar{A}_1 E_7$ where $G_2 < E_7$ is non-$E_7$-cr \\
\hline
\end{longtable}

In the following table, recall from Section \ref{sec:corollaries} that ``$\infty$-many classes'' refers to the fact that a group of type $A_{1}T_{1}$ contains infinitely many pairwise non-conjugate 1-dimensional tori, and so a non-$G$-cr simple subgroup with such a centraliser gives rise to infinitely many pairwise non-conjugate, non-$G$-cr reductive subgroups.

\begin{longtable}{>{\raggedright\arraybackslash}p{0.06\textwidth - 2\tabcolsep}>{\raggedright\arraybackslash}p{0.06\textwidth - 2\tabcolsep}>{\raggedright\arraybackslash}p{0.53\textwidth-\tabcolsep}@{}}
\caption{Non-$G$-cr reductive subgroups $X$ with $Z(X)^{\circ} \neq 1$} \label{tab:reductive} \\
\hline
$G$ & $p$ & Non-$G$-cr reductive subgroup $X$  \\
\hline
$E_6$ & 5 & $A_1 T_1$ where $A_{1} < A_{5}$ via $W(5)$ \\
& & $A_1 T_1$ where $A_{1} < D_{5}$ via $T(6)$ \\
& & $A_1 T_1$ where $A_{1} < D_{5}$ via $T(8)$ \\
$E_{7}$ & $5$ & $A_{1} T_{2}$ where $A_{1} < A_{5}$ via $W(5)$ \\
& & $A_{1} \bar{A}_{1} T_{1}$ where $A_{1} < A_{5}$ via $W(5)$ \\
& & $A_{1} T_{1}$ where $A_{1} < A_{5}$ via $W(5)$ \\
& & $A_{1} T_{2}$ where $A_{1} < A_{5}'$ via $W(5)$ \\
& & $A_{1} \bar{A}_{1} T_{1}$ where $A_{1} < A_{5}'$ via $W(5)$ \\
& & $A_{1} T_{1}$ where $A_{1} < A_{5}'$ via $W(5)$ ($\infty$-many classes) \\
& & $A_{1}T_{1}$ where $A_{1} < A_{6}$ via $W(6)$ \\
& & $A_{1}T_{1}$ where $A_{1} < \bar{A}_{1} A_{5}$ via $(1^{[r]},W(5)^{[s]})$ ($rs = 0)$ \\
& & $A_{1}T_{1}$ where $A_{1} < \bar{A}_{1} A_{5}'$ via $(1^{[r]},W(5)^{[s]})$ ($rs = 0)$ \\
& & $A_{1}T_{2}$ where $A_{1} < D_{5}$ via $T(6)$ \\
& & $A_{1}\bar{A}_{1}T_{1}$ where $A_{1} < D_{5}$ via $T(6)$ \\
& & $A_{1}T_{1}$ where $A_{1} < D_{5}$ via $T(6)$ ($\infty$-many classes) \\
& & $A_{1}T_{2}$ where $A_{1} < D_{5}$ via $T(8)$ \\
& & $A_{1}\bar{A}_{1}T_{1}$ where $A_{1} < D_{5}$ via $T(8)$ \\
& & $A_{1}T_{1}$ where $A_{1} < D_{5}$ via $T(8)$ ($\infty$-many classes) \\
& & $A_{1}T_{1}$ where $A_{1} < \bar{A}_{1} D_{5}$ via $(1^{[r]},T(6)^{[s]})$\\
& & $A_{1}T_{1}$ where $A_{1} < \bar{A}_{1} D_{5}$ via $(1^{[r]},T(8)^{[s]})$\\
$E_{8}$ & $7$ & $A_{1} T_{1}$ where $A_{1} < A_{7}$ via $W(7)$\\
& & $A_{1}T_{1}$ where $A_{1} < A_{7}'$ via $W(7)$\\
& & $A_{1}T_{1}$ where $A_{1} < D_{7}$ via $T(8)$\\
& & $A_{1}T_{1}$ where $A_{1} < D_{7}$ via $T(10)$\\
& & $A_{1}T_{1}$ where $A_{1} < D_{7}$ via $T(12)$\\
& & $A_{1}T_{1}$ where $A_{1} < A_{1}G_{2}$ via $(1,6)$\\
& & $G_{2}T_{1}$ where $G_{2} < E_{7}$ is non-$E_{7}$-cr\\
\hline
\end{longtable}

\section{Further module decompositions} \label{tabs:misc}

The following tables give the restrictions of certain $H$-modules to $X$ for $p = 5$ or $7$, when $X$ is of type $A_{1}$ and $H$ is a certain semisimple subgroup of $G$ containing $X$. These modules have been used implicitly in Sections \ref{sec:e6}--\ref{sec:e8} and in calculating the given actions in Section \ref{sec:tables}. The given structure has been calculated, and can be verified, in the manner described in Section \ref{sec:restrictions}. 

\begin{longtable}{>{\raggedright\arraybackslash}p{0.3\textwidth - 2\tabcolsep}>{\raggedright\arraybackslash}p{0.35\textwidth - 2\tabcolsep}>{\raggedright\arraybackslash}p{0.35\textwidth-\tabcolsep}@{}}
\caption{Alternating Powers of Certain $A_1$-modules \label{tab:extpowers}} \\ \hline
$V$ & $\bigwedge^2(V)$ & $\bigwedge^3(V)$ \\ \hline
$1 \otimes 1^{[r]}$ & $2 + 2^{[r]}$ & $1 \otimes 1^{[r]}$\\
$2^{[r]} \otimes 1^{[s]}$ & $4^{[r]} + 2^{[r]} \otimes 2^{[s]} + 0$ & $4^{[r]}  \otimes 1^{[s]} + 3^{[s]} + 2^{[r]} \otimes 1^{[s]}$ \\
\hline
\end{longtable}

\begin{longtable}{>{\raggedright\arraybackslash}p{0.3\textwidth - 2\tabcolsep}>{\raggedright\arraybackslash}p{0.3\textwidth - 2\tabcolsep}>{\raggedright\arraybackslash}p{0.4\textwidth-\tabcolsep}@{}}
\caption{Spin modules for $D_{4}$ restricted to irreducible subgroups of type $A_{1}$}\\
\hline
$V_{D_4}(\lambda_1) \downarrow X$ & $V_{D_4}(\lambda_3) \downarrow X$ & $V_{D_4}(\lambda_4) \downarrow X$ \\
\hline
$1^{[r]} \otimes 1^{[s]} + 1^{[t]} \otimes 1^{[u]}$ & $1^{[r]} \otimes 1^{[u]} + 1^{[s]} \otimes 1^{[t]}$ & $1^{[r]} \otimes 1^{[t]} + 1^{[s]} \otimes 1^{[u]}$ (or vice versa) \\
$3^{[r]} \otimes 1^{[s]}$ & $3^{[r]} \otimes 1^{[s]}$ & $4^{[r]} + 2^{[s]}$ (or vice versa)\\
$4^{[r]} + 2^{[s]}$ & $3^{[r]} \otimes 1^{[s]}$ & $3^{[r]} \otimes 1^{[s]}$ \\
$6 + 0$ $(p=7)$ & $6 + 0$ & $6 + 0$  \\
\hline
\end{longtable}

\begin{longtable}{>{\raggedright\arraybackslash}p{0.35\textwidth - 2\tabcolsep}>{\raggedright\arraybackslash}p{0.35\textwidth-\tabcolsep}@{}}
\caption{Spin modules for $D_{5}$ restricted to various subgroups of type $A_{1}$ \label{tab:spinD5}}\\
\hline
$V_{D_5}(\lambda_1) \downarrow X$ & $V_{D_5}(\lambda_4) \downarrow X$ \\
\hline
$2 + 2^{[r]} + 2^{[s]} + 0$ & $(1 \otimes 1^{[r]} \otimes 1^{[s]})^2$ \\
$2^{[r]} + 2^{[s]} + 1^{[t]} \otimes 1^{[u]}$ & $1^{[r]} \otimes 1^{[s]} \otimes 1^{[t]} + 1^{[r]} \otimes 1^{[s]} \otimes 1^{[u]}$ \\
$4^{[r]} + 1^{[s]} \otimes 1^{[t]} + 0$ & $3^{[r]} \otimes 1^{[s]} + 3^{[r]} \otimes 1^{[t]}$ \\
$4 + 4^{[r]}$ & $3 \otimes 3^{[r]}$ \\
$6^{[r]} + 2^{[s]}$ $(p=7)$ & $6^{[r]} \otimes 1^{[s]} + 1^{[s]}$ \\
$T(6)$ $(p=5)$ & $T(7) + W(5)$ \\
$T(8)$ $(p=5)$ & $W(10) + 4$ \\
$2 \otimes 2^{[1]} + 0$ & $3 \otimes 1^{[1]} + 3^{[1]} \otimes 1$ \\
\hline
\end{longtable}

\begin{longtable}{>{\raggedright\arraybackslash}p{0.3\textwidth - 2\tabcolsep}>{\raggedright\arraybackslash}p{0.35\textwidth - 2\tabcolsep}>{\raggedright\arraybackslash}p{0.35\textwidth-\tabcolsep}@{}}
\caption{Spin modules for $D_{6}$ restricted to irreducible subgroups of type $A_{1}$}\\
\hline
$V_{D_6}(\lambda_1) \downarrow X$ & $V_{D_6}(\lambda_5) \downarrow X$ & $V_{D_6}(\lambda_6) \downarrow X$ \\
\hline
$1^{[r]} \otimes 1^{[s]} + 1^{[t]} \otimes 1^{[u]} + 1^{[v]} \otimes 1^{[w]}$ & $1^{[r]} \otimes 1^{[t]} \otimes 1^{[v]} + 1^{[r]} \otimes 1^{[u]} \otimes 1^{[w]} + 1^{[s]} \otimes 1^{[t]} \otimes 1^{[w]} + 1^{[s]} \otimes 1^{[u]} \otimes 1^{[v]}$ & $1^{[r]} \otimes 1^{[t]} \otimes 1^{[w]} + 1^{[r]} \otimes 1^{[u]} \otimes 1^{[v]} + 1^{[s]} \otimes 1^{[t]} \otimes 1^{[v]} + 1^{[s]} \otimes 1^{[u]} \otimes 1^{[w]}$ \\
$2 + 2^{[r]} + 2^{[s]} + 2^{[t]}$ & $(1 \otimes 1^{[r]} \otimes 1^{[s]} \otimes 1^{[t]})^2$ & same \\
$2^{[r]} \otimes 1^{[s]} \otimes 1^{[t]}$ & $4^{[r]} \otimes 1^{[s]} + 3^{[s]} + 2^{[r]} \otimes 2^{[t]} \otimes 1^{[s]}$ & $4^{[r]} \otimes 1^{[t]} + 3^{[t]} + 2^{[r]} \otimes 2^{[s]} \otimes 1^{[t]}$ (or vice versa) \\
$2^{[r]} \otimes 2^{[s]} + 2^{[t]}$ & $3^{[r]} \otimes 1^{[s]} \otimes 1^{[t]} + 3^{[s]} \otimes 1^{[r]} \otimes 1^{[t]}$ & same \\
$3^{[r]} \otimes 1^{[s]} + 1^{[t]} \otimes 1^{[u]}$ & $3^{[r]} \otimes 1^{[s]} \otimes 1^{[t]} + 4^{[r]} \otimes 1^{[u]} + 2^{[s]} \otimes 1^{[u]}$ & $3^{[r]} \otimes 1^{[s]} \otimes 1^{[u]} + 4^{[r]} \otimes 1^{[t]} + 2^{[s]} \otimes 1^{[t]}$ (or vice versa) \\
$4^{[r]} + 2^{[s]} + 1^{[t]} \otimes 1^{[u]}$ & $3^{[r]} \otimes 1^{[s]} \otimes 1^{[t]} + 3^{[r]} \otimes 1^{[s]} \otimes 1^{[u]}$  & same  \\
$5^{[r]} \otimes 1^{[s]}$ $(p=7)$ & $T(8)^{[r]} \otimes 1^{[s]} + 3^{[s]}$  &  $T(9)^{[r]} + 5^{[r]} \otimes 2^{[s]}$ (or vice versa) \\
$6^{[r]} + 4^{[s]}$ $(p=7)$ & $6^{[r]} \otimes 3^{[s]} + 3^{[s]}$ & same \\
$6^{[r]} + 1^{[s]} \otimes 1^{[t]} + 0$ $(p=7)$ & $6^{[r]} \otimes 1^{[s]} + 6^{[r]} \otimes 1^{[t]} + 1^{[s]} \otimes 1^{[t]}$ &  same \\
\hline
\end{longtable}

\begin{longtable}{>{\raggedright\arraybackslash}p{0.34\textwidth - 2\tabcolsep}>{\raggedright\arraybackslash}p{0.66\textwidth-\tabcolsep}@{}}
\caption{Spin modules for $D_{7}$ restricted to various subgroups of type $A_{1}$ when $p=7$}\\
\hline
$V_{D_7}(\lambda_1) \downarrow X$ & $V_{D_7}(\lambda_6) \downarrow X$ \\
\hline
$2^{[r]} + 2^{[s]} + 1^{[t]} \otimes 1^{[u]} + 1^{[v]} \otimes 1^{[w]}$ & ${1^{[r]} \otimes 1^{[s]} \otimes 1^{[t]} \otimes 1^{[v]}} + {1^{[r]} \otimes 1^{[s]} \otimes 1^{[u]} \otimes 1^{[w]}} + {1^{[r]} \otimes 1^{[s]} \otimes 1^{[t]} \otimes 1^{[v]}} + {1^{[r]} \otimes 1^{[s]} \otimes 1^{[u]} \otimes 1^{[v]}}$ \\
$3^{[r]} \otimes 1^{[s]} + 2^{[t]} + 2^{[u]}$ & $3^{[r]} \otimes 1^{[s]} \otimes 1^{[t]} \otimes 1^{[u]} + 4^{[r]} \otimes 1^{[t]} \otimes 1^{[u]} + 2^{[s]} \otimes 1^{[t]} \otimes 1^{[u]}$ \\
$4^{[r]} + 1^{[s]} \otimes 1^{[t]} + 1^{[u]} \otimes 1^{[v]} + 0$ & $3^{[r]} \otimes 1^{[s]} \otimes 1^{[u]} + 3^{[r]} \otimes 1^{[t]} \otimes 1^{[v]} + 3^{[r]} \otimes 1^{[s]} \otimes 1^{[v]} + 3^{[r]} \otimes 1^{[t]} \otimes 1^{[u]}$ \\
$4^{[r]} + 2^{[s]} + 2^{[t]} + 2^{[u]}$ & $(3^{[r]} \otimes 1^{[s]} \otimes 1^{[t]} \otimes 1^{[u]})^2$ \\
$4^{[r]} + 4^{[s]} + 2^{[t]} + 0$ & $(3^{[r]} \otimes 3^{[s]} \otimes 1^{[t]})^2$ \\
$4^{[r]} + 3^{[s]} \otimes  1^{[t]} + 0$ & $3^{[r]} \otimes 3^{[s]} \otimes 1^{[t]} + 3^{[r]} \otimes 4^{[s]} + 3^{[r]} \otimes 2^{[t]}$ \\
$4^{[r]} + 4^{[s]} + 1^{[t]} \otimes 1^{[u]}$ & $3^{[r]} \otimes 3^{[s]} \otimes 1^{[t]} + 3^{[r]} \otimes 3^{[s]} \otimes 1^{[u]}$ \\
$6^{[r]} + 2^{[s]} + 1^{[t]} \otimes 1^{[u]}$ & $6^{[r]} \otimes 1^{[s]} \otimes 1^{[t]} + 6^{[r]} \otimes 1^{[s]} \otimes 1^{[u]} + 1^{[s]} \otimes 1^{[t]} + 1^{[s]} \otimes 1^{[u]}$  \\
$6 + 6^{[r]}$ & ${6 \otimes 6^{[r]}} + 6 + 6^{[r]} + 0$ \\
$T(8)$ $(p=7)$ & $13 + T(11) + T(9) + T(7) + W(7)$ \\
$T(10)$ $(p=7)$ & $T(16) + W(14) + T(8) + 6$ \\
$T(12)$ $(p=7)$ & $W(21) + T(15) + T(9)$ \\
\hline
\end{longtable}

\begin{longtable}{>{\raggedright\arraybackslash}p{0.36\textwidth - 2\tabcolsep}>{\raggedright\arraybackslash}p{0.38\textwidth-\tabcolsep}@{}}
\caption{$V_{27} = V_{E_6}(\lambda_1)$ restricted to various subgroups of type $A_{1}$ when $p=5,7$}\\
\hline
$X$ & $V_{27} \downarrow X$ \\
\hline
$X < A_1 A_5$ via $(1^{[r]},5^{[s]})$ $(p=7)$ & $1^{[r]} \otimes 5^{[s]} + T(8)^{[s]} + 0$ \\
$X < A_1 A_5$ via $(1^{[r]},2^{[s]} \otimes 1^{[t]})$ & $1^{[r]} \otimes 2^{[s]} \otimes 1^{[t]} + 4^{[s]} + 2^{[s]} \otimes 2^{[t]} + 0$ \\
$X < A_2 G_2$ via $(2^{[r]},6^{[s]})$ $(p=7)$ & $4^{[r]} + 2^{[r]} \otimes 6^{[s]} + 0$ \\
$X < A_2^3$ via $(2,2^{[r]}, 2^{[s]})$ & $2 \otimes 2^{[r]} + 2 \otimes 2^{[s]} + 2^{[r]} \otimes 2^{[s]}$ \\
\hline
\end{longtable}

\begin{longtable}{>{\raggedright\arraybackslash}p{0.35\textwidth - 2\tabcolsep}>{\raggedright\arraybackslash}p{0.35\textwidth-\tabcolsep}@{}}
\caption{$V_{56} = V_{E_7}(\lambda_7)$ restricted to various subgroups of type $A_{1}$ when $p=7$}\\
\hline
$X$ & $V_{56} \downarrow X$ \\
\hline
$X < A_1 D_6$ & $V_{56} \downarrow A_1 D_6 = (1,\lambda_1) + (0,\lambda_6)$ \\
$X < A_1 A_1$ via $(1^{[r]},1^{[s]})$ & $6^{[r]} \otimes 3^{[s]} + 4^{[r]} \otimes 1^{[s]} + 2^{[r]} \otimes 5^{[s]}$ \\
$X < A_1 G_2$ via $(1^{[r]},6^{[s]})$ & ${3^{[r]} \otimes 6^{[s]}} + {1^{[r]} \otimes T(10)^{[s]}}$\\
$X < G_2 C_3$ via $(6^{[r]},5^{[s]})$ & $6^{[r]} \otimes 5^{[s]} + T(9)^{[s]}$ \\
$X < G_2 C_3$ via $(6^{[r]},2^{[s]} \otimes 1^{[t]})$ & $6^{[r]} \otimes 2^{[s]} \otimes 1^{[t]} + 4^{[s]} \otimes 1^{[t]} + 3^{[t]}$ \\
\hline
\end{longtable}

\section{Acknowledgements} 
The second author acknowledges financial support from the EPSRC and the hospitality at the University of Auckland, where much of the work for this paper was undertaken. We would like to thank Dr.\ D.\ Stewart for helpful conversations as well as Prof.\ M.\ Liebeck for bringing the problem to our attention and for his suggestions on an early version of the paper. Finally, we are extremely grateful to our anonymous referee for their careful reading of the paper and numerous helpful suggestions. 

\vspace{-3pt}

\bibliographystyle{amsplain}
\bibliography{biblio}

\end{document}